\documentclass{article}

\usepackage{lipsum} 

\usepackage[numbers]{natbib}

\title{An Introduction to Higher Categorical Algebra}
\author{David Gepner}

\usepackage[pdfborder=0,pagebackref,colorlinks,citecolor=darkblue,linkcolor=darkgreen,urlcolor=darkblue]{hyperref}

\usepackage{amsmath}
\usepackage{amscd}
\usepackage{amsbsy}
\usepackage{amssymb}
\usepackage{verbatim}
\usepackage{eufrak}
\usepackage{eucal}
\usepackage{microtype}
\usepackage{hyperref}
\usepackage{mathrsfs}
\usepackage{amsthm}
\usepackage{stmaryrd}
\usepackage{xy}
\input xy
\xyoption{all}
\usepackage{tikz}
\usetikzlibrary{matrix,arrows}
\usepackage{dsfont}
\usepackage[T1]{fontenc}
\usepackage{mathbbol}

\usepackage{makeidx}
\makeindex

\usepackage{color}
\definecolor{todo}{rgb}{1,0,0}
\definecolor{conditional}{rgb}{0,1,0}
\definecolor{e-mail}{rgb}{0,.40,.80}
\definecolor{reference}{rgb}{.20,.60,.22}
\definecolor{mrnumber}{rgb}{.80,.40,0}
\definecolor{citation}{rgb}{0,.40,.80}

\renewcommand{\rm}{\mdseries}

\let\oldmarginpar\marginpar
\renewcommand\marginpar[1]{\-\oldmarginpar[\raggedleft\footnotesize #1]%
{\raggedright\footnotesize #1}}

\usepackage{makeidx}

\usepackage{enumerate}
\usepackage{sseq}
\usepackage{amsmath}
\usepackage{amssymb}
\usepackage{stmaryrd}
\usepackage{epsfig}
\usepackage{psfrag}
\usepackage{bbm}     

\usepackage{eucal}
\usepackage{calligra}
\usepackage{amsthm}
\usepackage{graphicx}
\usepackage{mathrsfs}
\usepackage{comment}

\usepackage{xcolor}

\definecolor{darkgreen}{rgb}{0,0.30,0} 
\definecolor{darkred}{rgb}{0.75,0,0}
\definecolor{darkblue}{rgb}{0,0,0.6}

\def\makeautorefname#1#2{\expandafter\def\csname#1autorefname\endcsname{#2}}

\makeautorefname{equation}{Equation}%
\makeautorefname{footnote}{footnote}%
\makeautorefname{item}{item}%
\makeautorefname{figure}{Figure}%
\makeautorefname{table}{Table}%
\makeautorefname{part}{Part}%
\makeautorefname{appendix}{Appendix}%
\makeautorefname{chapter}{Chapter}%
\makeautorefname{section}{Section}%
\makeautorefname{subsection}{Section}%
\makeautorefname{subsubsection}{Section}%
\makeautorefname{paragraph}{Paragraph}%
\makeautorefname{subparagraph}{Paragraph}%
\makeautorefname{theorem}{Theorem}%
\makeautorefname{thm}{Theorem}%
\makeautorefname{addm}{Addendum}%
\makeautorefname{mainthm}{Main theorem}%
\makeautorefname{corollary}{Corollary}%
\makeautorefname{cor}{Corollary}%
\makeautorefname{lemma}{Lemma}%
\makeautorefname{lem}{Lemma}%
\makeautorefname{sublemma}{Sublemma}%
\makeautorefname{sublem}{Sublemma}%
\makeautorefname{subl}{Sublemma}%
\makeautorefname{proposition}{Proposition}%
\makeautorefname{property}{Property}
\makeautorefname{pro}{Property}
\makeautorefname{sch}{Scholium}%
\makeautorefname{step}{Step}%
\makeautorefname{conject}{Conjecture}%
\makeautorefname{conj}{Conjecture}%
\makeautorefname{questn}{Question}
\makeautorefname{quest}{Question}
\makeautorefname{qn}{Question}
\makeautorefname{definition}{Definition}%
\makeautorefname{defn}{Definition}%
\makeautorefname{defi}{Definition}%
\makeautorefname{def}{Definition}%
\makeautorefname{dfn}{Definition}%
\makeautorefname{notation}{Notation}
\makeautorefname{notn}{Notation}
\makeautorefname{rem}{Remark}%
\makeautorefname{remark}{Remark}
\makeautorefname{rems}{Remarks}%
\makeautorefname{rmk}{Remark}%
\makeautorefname{rk}{Remark}%
\makeautorefname{remarks}{Remarks}%
\makeautorefname{rems}{Remarks}%
\makeautorefname{rmks}{Remarks}%
\makeautorefname{rks}{Remarks}%
\makeautorefname{example}{Example}%
\makeautorefname{examp}{Example}%
\makeautorefname{exmp}{Example}%
\makeautorefname{exam}{Example}%
\makeautorefname{exa}{Example}%
\makeautorefname{axiom}{Axiom}%
\makeautorefname{axi}{Axiom}%
\makeautorefname{ax}{Axiom}%
\makeautorefname{case}{Case}%
\makeautorefname{claim}{Claim}%
\makeautorefname{clm}{Claim}%
\makeautorefname{assumpt}{Assumption}%
\makeautorefname{asses}{Assumptions}%
\makeautorefname{conclusion}{Conclusion}%
\makeautorefname{concl}{Conclusion}%
\makeautorefname{conc}{Conclusion}%
\makeautorefname{cond}{Condition}%
\makeautorefname{const}{Construction}%
\makeautorefname{con}{Construction}%
\makeautorefname{criterion}{Criterion}%
\makeautorefname{criter}{Criterion}%
\makeautorefname{crit}{Criterion}%
\makeautorefname{exercise}{Exercise}%
\makeautorefname{exer}{Exercise}%
\makeautorefname{exe}{Exercise}%
\makeautorefname{problem}{Problem}%
\makeautorefname{problm}{Problem}%
\makeautorefname{prob}{Problem}%
\makeautorefname{prob}{Problem}%
\makeautorefname{soln}{Solution}%
\makeautorefname{sol}{Solution}%
\makeautorefname{sum}{Summary}%
\makeautorefname{oper}{Operation}%
\makeautorefname{obs}{Observation}%
\makeautorefname{ob}{Observation}%
\makeautorefname{conv}{Convention}%
\makeautorefname{cvn}{Convention}%
\makeautorefname{warn}{Warning}%
\makeautorefname{note}{Note}%
\makeautorefname{fact}{Fact}%
\makeautorefname{ouch0}{Counterexample}%
\makeautorefname{introthm}{Theorem}%

\newtheorem{theorem}{Theorem}[subsection]

\newtheorem{corollary}{Corollary}[subsection]
\newtheorem{proposition}{Proposition}[subsection]

\theoremstyle{definition}

\newtheorem{definition}{Definition}[subsection]

\newtheorem{example}{Example}[subsection]

\newtheorem{remark}{Remark}[subsection]

\makeatletter

\let\c@corollary=\c@thm

\let\c@proposition=\c@thm

\let\c@lemma=\c@thm

\let\c@conj=\c@thm

\let\c@definition=\c@thm

\let\c@df=\c@thm

\let\c@example=\c@thm

\let\c@remark=\c@thm

\let\c@sch=\c@thm

\let\c@equation\c@thm

\makeatother

\newcommand{\A}{\mathcal{A}}
\newcommand{\B}{\mathcal{B}}
\newcommand{\C}{\mathcal{C}}
\newcommand{\D}{\mathcal{D}}
\newcommand{\E}{\mathcal{E}}
\newcommand{\F}{\mathcal{F}}

\newcommand{\I}{\mathcal{I}}

\newcommand{\M}{\mathcal{M}}
\newcommand{\N}{\mathcal{N}}
\renewcommand{\O}{\mathcal{O}}
\renewcommand{\P}{\mathcal{P}}

\renewcommand{\S}{\mathcal{S}}
\newcommand{\T}{\mathcal{T}}

\renewcommand{\AA}{\mathbf{A}}

\newcommand{\EE}{\mathbf{E}}
\newcommand{\FF}{\mathbf{F}}
\newcommand{\GG}{\mathbf{G}}

\newcommand{\NN}{\mathbf{N}}

\newcommand{\PP}{\mathbf{P}}
\newcommand{\QQ}{\mathbf{Q}}

\renewcommand{\SS}{\mathbf{S}}

\newcommand{\ZZ}{\mathbf{Z}}

\newcommand{\CCC}{\mathfrak{C}}

\renewcommand{\i}{\infty}

\newcommand{\too}{\longrightarrow}
\newcommand{\from}{\longleftarrow}
\newcommand{\op}{\mathrm{op}}
\DeclareMathOperator{\et}{\acute{e}t}

\newcommand{\n}{\langle n\rangle}
\newcommand{\m}{\langle m\rangle}

\DeclareMathOperator{\LMod}{LMod}
\DeclareMathOperator{\RMod}{RMod}
\newcommand{\Dual}{\mathrm{D}}

\DeclareMathOperator{\Mod}{Mod}

\DeclareMathOperator{\Alg}{Alg}
\DeclareMathOperator{\CAlg}{CAlg}

\DeclareMathOperator{\Sub}{Sub}

\DeclareMathOperator{\Aut}{Aut}

\DeclareMathOperator{\cocart}{cocart}
\DeclareMathOperator{\cart}{cart}
\DeclareMathOperator{\Cat}{Cat}
\DeclareMathOperator{\CAT}{CAT}
\DeclareMathOperator{\cof}{cof}
\DeclareMathOperator{\Der}{Der}

\DeclareMathOperator{\End}{End}

\DeclareMathOperator{\Exc}{Exc}
\DeclareMathOperator{\Ext}{Ext}

\DeclareMathOperator{\Fin}{Fin}
\DeclareMathOperator{\f}{fin}

\DeclareMathOperator{\Fun}{Fun}
\newcommand{\Funl}{\mathrm{LFun}}
\newcommand{\Funr}{\mathrm{RFun}}

\DeclareMathOperator{\Gpd}{Gpd}
\DeclareMathOperator{\Ho}{Ho}
\DeclareMathOperator{\Hom}{Hom}
\DeclareMathOperator{\id}{id}

\DeclareMathOperator{\loc}{loc}

\DeclareMathOperator{\nil}{nil}

\DeclareMathOperator{\Spec}{Spec}

\DeclareMathOperator{\proj}{proj}

\DeclareMathOperator{\Map}{Map}

\DeclareMathOperator{\Mon}{Mon}

\DeclareMathOperator{\colim}{colim}

\DeclareMathOperator{\PrL}{LPr}
\DeclareMathOperator{\Prl}{LPr}
\DeclareMathOperator{\PrR}{RPr}
\DeclareMathOperator{\Prr}{RPr}

\DeclareMathOperator{\PSp}{PSp}

\DeclareMathOperator{\Set}{Set}

\DeclareMathOperator{\Sym}{Sym}

\DeclareMathOperator{\fib}{fib}

\DeclareMathOperator{\Sp}{Sp}

\DeclareMathOperator{\GL}{GL}

\DeclareMathOperator{\Ten}{Ten}

\newcommand{\Ab}{\mathrm{Ab}}

\newcommand{\cn}{\mathrm{cn}}
\newcommand{\st}{\mathrm{st}}

\newcommand{\pt}{\mathrm{pt}}

\DeclareMathOperator{\fin}{fin}

\def\fibdown{\ar@{->>}[d]}
\def\hookdown{\ar@<-.5ex>[d]|{\phantom{a}}|<<{\put(-.7,2){$\scriptstyle\cap$}}}

\def\llarrow{  \hspace{.05cm}\mbox{\,\put(0,-2){$\leftarrow$}\put(0,2){$\leftarrow$}\hspace{.45cm}}}

\def\lllarrow{ \hspace{.05cm}\mbox{\,\put(0,-3){$\leftarrow$}\put(0,1){$\leftarrow$}\put(0,5){$\leftarrow$}\hspace{.45cm}}}

\def\llllarrow{\hspace{.05cm}\mbox{\,\put(0,-3){$\leftarrow$}\put(0,.5){$\leftarrow$}\put(0,4){$\leftarrow$}\put(0,7.5){$\leftarrow$}
               \hspace{.45cm}}}

\renewcommand{\epsilon}{\varepsilon}

\DeclareMathOperator{\Ind}{Ind}

\newcommand{\pb}{\ar@{}[dr]|{\text{\pigpenfont J}}}
\newcommand{\po}{\ar@{}[dr]|{\text{\pigpenfont R}}}

\newcommand{\CMon}{\mathrm{CMon}}

\author{David Gepner\thanks{David Gepner was supported by NSF Grant DMS-1406529 and the Mathematical Sciences Research Institute Program ``Derived Algebraic Geometry''.}}

\makeindex

\begin{document}

\maketitle

\tableofcontents

%

\section{Introduction}\label{sec:it}

\subsection{Higher structures}

Higher algebra is the study of algebraic structures which arise in the setting of higher category theory.
Higher algebra generalizes ordinary algebra, or algebra in the setting of ordinary category theory.
Ordinary categories have sets of morphisms between objects, and elements of a set are either equal or not.
Higher categories, on the other hand, have homotopy types of morphisms between objects, typically called {\em mapping spaces}. Sets are examples of homotopy types, namely the discrete ones, but in general it doesn't quite make sense to ask whether or not two ``elements'' of a homotopy type are ``equal''; rather, they are equivalent if they are represented by points which can be connected by a path in some suitable model for the homotopy type. But then any two such paths might form a nontrivial loop, leading to higher automorphisms, and so on.
The notion of equality only makes sense after passing to discrete invariants such as homotopy groups.

Since the higher categorical analogue of a set is a space, the higher categorical analogue of an abelian group ought to be something like a space equipped with a multiplication operation which is associative, commutative, and invertible up to coherent homotopy.
While invertibility is a {\em property} of an associative operation, commutativity is not; rather, it is {\em structure}.
This is because, in higher categorical contexts, it is not enough to simply permute a sequence of elements; instead, the permutation itself is recorded as a morphism.
A commutative multiplication operation must also act on morphisms, so that they may also be permuted, and so on and so forth, provided we keep track of these permutations as still higher morphisms.
There are a number of formalisms which make this precise, all of which are equivalent to (or obvious variations on) the notion of {\em spectrum} in the sense of algebraic topology.\footnote{The overuse of the term ``spectrum'' in mathematics is perhaps a potential cause for confusion; fortunately, it is almost always clear from context what is meant.}

On the one hand, a spectrum is an infinite delooping of a pointed space, thereby providing an abelian group structure on all its homotopy groups (positive and negative); on the other, a spectrum represents a cohomology theory, which is to say a graded family of contravariant abelian group valued functors on pointed spaces satisfying a suspension relation and certain exactness conditions.  These two notions are equivalent: the  functors which comprise the cohomology theory are represented by the spaces of the infinite delooping.
The real starting point of higher algebra is the observation that there is a symmetric monoidal structure on the $\i$-category of spectra which refines the tensor product of abelian groups, and that there are many important examples of algebras for this tensor product.

\subsection{Overview}

Ordinary algebra is set based, meaning that it is carried out in the language of ordinary categories.
As mentioned, the higher categorical analogue of sets are spaces, or $n$-truncated spaces if one chooses to work in an $(n+1)$-categorical context, and the truncation functors allow us to switch back and forth between categorical levels.
In the category of sets, limits and colimits reduce to intersections and unions in some ambient set; in higher category theory, however, these operations must be interpreted invariantly, which implies that a homotopy colimit of sets (viewed as spaces) need not be discrete.

We begin our exposition in \autoref{sec:ct} with some background on the behavior of colimits in higher categories, especially in $\i$-categories of presheaves or certain full subcategories thereof (the presentable $\i$-categories).
We then turn to the Grothendieck construction, which establishes a correspondence between fibrations and functors. The fibration perspective allows for an efficient approach to the theory of (symmetric) monoidal $\i$-categories and (commutative) algebras and modules therein, our main objects of interest, at least in the stable setting: a stable $\i$-category is a higher categorical analogue of an abelian category, an analogy which we make precise by comparing derived categories of abelian categories with stable $\i$-categories via t-structures.

Having equipped ourselves with the basic structures and language, in \autoref{sec:rt} we turn to a more detailed study of spectra and the smash product.
An associative (respectively, commutative) ring spectrum is defined as an algebra (respectively, commutative algebra) object in the stable $\i$-category of spectra, the universal stable $\i$-category.
The theory also allows for a notion of (left or right) module object of an $\i$-category which is (left or right) tensored over a monoidal $\i$-category.
We conclude this section with some remarks on localizations of ring spectra, which mirrors the ordinary theory save for the fact that ideals must be interpreted on the level of homotopy groups.

\autoref{sec:mt} is devoted to module theory.
In particular, monads appear as an instance of modules, allowing us to address monadicity, which plays a much more important role higher categorically due do the difficultly of ad hoc constructions.
Simplicial objects and their colimits, geometric realizations, feature in the construction of the relative tensor product as well as the definition of projective module.
We also study the more general class of perfect modules, which are colimits of shifted projective modules in a sense made rigorous by the theory of tor-amplitude, which acts as a substitute for projective resolutions.
We also consider free algebras and isolate various finiteness properties of modules and algebras which play important roles in higher algebra.

The final \autoref{sec:dt} deals with deformations of commutative algebras.
The formalism of the tangent bundle allows for an elegant construction of the cotangent complex, which governs derivations and square-zero extensions.
The Postnikov tower of truncations of a connective commutative algebra is comprised entirely of square-zero extensions, a crucial fact which implies an obstruction theory for computing the space of maps between commutative algebra spectra.
The main theorem in that all obstructions vanish in the \'etale case, from which it follows that the $\i$-category of \'etale commutative algebras over of fixed commutative ring spectrum $R$ is equivalent to an ordinary category, namely that of \'etale commutative algebras over its underlying ordinary ring $\pi_0 R$, a strong version of the topological invariance of \'etale morphisms.

We conclude this introductory section with some background on various approaches to homotopy coherence and some remarks on higher categorical ``set theory'': small and large spaces, categories, universes, etc.
The reader who is already familiar with these notions is encouraged to \autoref{sec:ct}, or even \autoref{sec:rt} if they are already familiar with stability, presentability, and symmetric monoidal structures.

\subsection{Homotopy coherence}

To lessen the prerequisites we avoid the use of operads or $\infty$-operads in this article altogether.
Nevertheless, for the sake of putting the theory of higher categorical algebra into historical context, and explaining some of the standard terminology, a few remarks are in order.

A space equipped with a homotopy coherently associative, or homotopy coherently associative and commutative, multiplication operation is traditionally referred to as an {\em $\AA_\i$-monoid}, or {\em $\EE_\i$-monoid}, meaning that it admits an action by an $\AA_\i$ (infinitely homotopy coherently associative) or $\EE_\i$ (infinitely homotopy coherently ``everything'', i.e. associative and commutative) operad. If the multiplication operation is invertible up to homotopy, the $\AA_\i$-monoid is said to be {\em grouplike}.
An $\EE_\i$-monoid is grouplike when viewed as an $\AA_\i$-monoid.
These operads were originally constructed out of geometric objects like associahedra, configuration spaces, or spaces of linear isometries.

The $\i$-categorical approach prefers to use small combinatorial models for associativity and commutativity, as in Segal's treatment \cite{Seg74}, by incorporating homotopy coherence into the language itself.
The result is a significantly more streamlined approach to homotopy coherent algebraic structures, as anticipated in the now extremely influential book of the same name of Boardmann--Vogt \cite{BV}, which contained the original definition of $\i$-category, well before the theory was systematically developed by Joyal and then Lurie.
Nevertheless, there is a rich interplay between geometry, topology, and higher category theory, as evidenced by the remarkable {\em cobordism hypothesis}, among other things.
Even the motivating example of the theory of operads, as originally developed by May \cite{May}, namely the {\em little $n$-cubes} operad, $n\in\NN$, collectively form the most important family of $\i$-operads, the so-called $\EE_n$ operads.

As $\i$-operads, $\AA_\i\simeq\EE_1$ and $\EE_1^{\otimes n}\simeq\EE_n$; the former equivalence is easy but the latter is equivalence is hard and requires both the {\em tensor product} of Boardmann--Vogt \cite{BV} and the  {\em additivity theorem} of Dunn \cite{Dunn}.
Since we will only be concerned with $\AA_\i$-algebras and $\EE_\i$-algebras in (symmetric) monoidal $\i$-categories, we choose to emphasize the analogy with ordinary algebra by referring to these objects as {\em associative and commutative algebras}, respectively.
This drastically simplifies the terminology and also allows us to reformulate actions by the associative or commutative $\i$-operad, respectively, in terms of functors from a category of ordinals or cardinals (arguably the must basic mathematical objects of all).
And while the abstract theory is quite powerful, one should keep in mind that many of the most important examples come from geometry and topology via these more classical constructions.

Heuristically, a $\infty$-categories are generalizations of (ordinary) categories in which there is a space, instead of a set, of morphisms between any pair of objects.
There are a number of ways of making this precise, but suffice it to say that {\em simplicially (or topologically) enriched categories} and {\em quasi-categories} yield equivalent models:
any $\infty$-category is equivalent to the {\em homotopy coherent nerve} $\mathrm{N}(\C)$ of a simplicially enriched category $\C$.
We will follow the usual notational conventions and sometimes refer to simplically enriched categories as {\em simplicial categories}, though this should not be confused with the more general notion of simplicial object in the category of categories.
\begin{remark}
The theory of quasi-categories \cite{Joy} has the distinct advantage that it allows for an easy construction of the $\i$-category $\Fun(\D,\C)$ of functors from an $\infty$-category $\D$ to an $\i$-category $\C$ as the exponential
\index{$\Fun(\D,\C)$}
\[
\Fun(\D,\C):=\C^{\D}
\]
in the cartesian closed category of simplicial sets.
This is a completely combinatorial object: if $X$ and $Y$ are simplicial sets, an $m$-simplex of $X^Y$ is natural transformation $\Delta^m\times Y\to X$ of functors $\Delta^{\op}\to\Set$, so it is completely determined by a compatible family of {\em functions} $\Delta^m_n\times Y_n\to X_n$, $n\in\NN$.
\end{remark}

\begin{remark}
The chief issue which arises when working with simplicial categories is that the simplicial category of simplicial functors from $\D$ to $\C$ is not in general invariant under weak equivalence (simplicial functors which are fully faithful and essentially surjective up to weak homotopy equivalence).
To obtain the homotopically correct simplicial category of functors we must replace $\D$ with a sufficiently ``free'' version $\D'\to\D$ of itself.
\end{remark}
\begin{remark}
In practice it is usually easier to apply the {\em homotopy coherent nerve} functor\index{$\mathrm{N}$}\index{homotopy coherent nerve}
\[
\mathrm{N}:\Cat_\Delta\too\Fun(\Delta^{\op},\Set)
\]
and work in simplicial sets.
Here $\Cat_\Delta$ denotes the category of simplicially enriched categories and $\mathrm{N}$ is the right adjoint of the colimit preserving functor $\CCC:\Fun(\Delta^{\op},\Set)\to\Cat_\Delta$\index{$\CCC$} determined by defining
$
\Map_{\CCC[\Delta^n]}(i,j)
$
to be the simplicial set of partially ordered subsets of $\{i,i+1,\ldots,j-1,j\}$.
\end{remark}
\begin{theorem}{\em \cite[Theorem 2.2.5.1]{HTT}}
The homotopy coherent nerve functor $\mathrm{N}:\Cat_\Delta\to\Fun(\Delta^{\op},\Set)$ is a right Quillen equivalence.
In particular, the homotopy coherent nerve $\mathrm{N}(\C)$ of a simplically enriched category $\C$ is a quasicategory provided $\C$ is enriched in Kan complexes (every $\Map_\C(A,B)$ is a Kan complex).
\end{theorem}

Higher categorical algebra is truly {\em homotopical} and not just {\em homological} in nature, meaning that many  of its most important objects simply do not exist within the world of chain complexes or derived categories.
The portion of the theory that can be formulated in these terms is {\em differential graded algebra}, the abstract study of which employs the language of {\em differential graded categories}.

\begin{definition}
A {\em differential graded category}
\index{differential graded category}
is a category enriched over the category $\mathrm{Ch}=\mathrm{Ch}(\Ab)$ of chain complexes of abelian groups.
We write $\Cat_\mathrm{dg}$\index{$\Cat_\mathrm{dg}$} for the category of small differential graded categories
\end{definition}

\begin{example}
Let ${\mathcal{A}}$ be an abelian category, or an additive subcategory of an abelian category.
Given complexes ${A, B \in \mathrm{Ch}(\mathcal{A})}$, there is a natural chain complex of abelian groups of homomorphisms $\underline{\Hom}(A, B)$, which in degree $n$ is the abelian group
$
{\underline{\Hom}(A, B)_n = \prod_{k\in\ZZ} \Hom_{\mathcal{A}}(A_k, B_{k+n})},
$
with differential given by the formula
$
(d\varphi)(a)=d(\varphi(a))-(-1)^n\varphi(d(a)).
$
In this way, ${\mathrm{Ch}(\mathcal{A})}$ acquires a canonical structure of a differential graded category: there is a chain complex of maps between any two objects, suitably compatible with composition.
Moreover, the abelian group of maps $\Hom(A,B)\cong \mathrm{H}_0\underline{\Hom}(A,B)$ falls out as the zeroth homology of this complex, so that $\mathrm{Ch}(\A)$ is the homotopy category of this differential graded category.
\end{example}

\begin{remark}\label{DK}
A differential graded category determines a simplicially enriched category via the {\em Dold-Kan correspondence}
\index{Dold-Kan correspondence} \cite{Wei94}, which asserts that there is a suitably monoidal equivalece of categories between connective chain complexes and simplicial abelian groups.
Said differently, we may regard a chain complex, provided it is identically zero in negative degrees,\footnote{We will use homological grading consistently throughout this article.}
as a simplicial abelian group.
This is done via the functor $\mathrm{Ch}_{\geq 0}( {\Ab}) \to\Fun(\Delta^{\op},\Ab)$ which associates to such a chain complex $A$ and nonempty finite ordinal $[n]$ the abelian group $\underset{[n]\twoheadrightarrow [m]}{\bigoplus}A_m$.
\end{remark}
\begin{remark}\label{dgn}
The functor $\mathrm{Ch}_{\geq 0}(\Ab)\to\Ab_\Delta$ alluded to above admits a lax monoidal structure \cite{Wei94}, and lax monoidal functors between enriching monoidal categories allow us to functorially change enrichment.
Since the (good) truncation functor $\tau_{\geq 0}:\mathrm{Ch}(\Ab)\to\mathrm{Ch}_{\geq 0}(\Ab)$ and underlying set functor $\Ab\to\Set$ also admit lax monoidal structures, we obtain a composite lax monoidal functor
\[
\mathrm{Ch}(\Ab)\too\mathrm{Ch}_{\geq 0}(\Ab)\too\Fun(\Delta^{\op},\Ab)\too\Fun(\Delta^{\op},\Set),
\]
and hence a functor from differential graded to simplicially enriched categories.
We write $(-)_\Delta:\Cat_\mathrm{dg}\to\Cat_\Delta$ for this functor and $\mathrm{N}(\C):=\mathrm{N}(\C_\Delta)$ for the homotopy coherent nerve of the differential graded category $\C$.
One can define a more explicit differential graded nerve functor $\mathrm{N_{dg}}:\Cat_\mathrm{dg}\to\Fun(\Delta^{\op},\Set)$, as in \cite[Construction 1.3.1.6]{HA}, but it is equivalent to the homotopy coherent nerve of the associated simplicially enriched category \cite[Proposition 1.3.1.17]{HA}).
\end{remark}

\subsection{Terminology and notation}

For the purposes of this article, we will always consider spaces from the point of view of $\infty$-category theory.
This means that we will tend to think of a topological space $X$ in terms of its associated singular simplicial set (its {\em fundamental $\infty$-groupoid})
\index{fundamental $\infty$-groupoid}
rather than as a point-set object, and that we reserve the right to replace $X$ with a homotopy equivalent space (or not even choose a representative of the homotopy type of $X$). 
If $X$ has the homotopy type of a cell complex then the realization of the singular complex of $X$ recovers $X$ up to homotopy equivalence; otherwise, the realization of the singular complex of $X$ recovers $X$ only up to weak homotopy equivalence.
As we will never need actual topological spaces (up to homoemorphism) at all anywhere in this article, ``space'' will always mean ``homotopy type'' in this sense. 
\begin{definition}
The $\infty$-category $\S$
\index{$\S$}
of {\em spaces} is the homotopy coherent nerve of the simplicially enriched category of Kan complexes, $\S=\mathrm{N}(\mathrm{Kan})$.
\end{definition}
\begin{remark}
As stipulated by Grothendieck's {\em homotopy hypothesis}, the $\infty$-category of spaces is equivalent to the $\infty$-category $\Gpd_\i$ of $\infty$-groupoids.
Here $\Gpd_\i\subset\Cat_\i$
\index{$\Gpd_\i$}
denotes the full subcategory consisting of those $\infty$-categories $\C$ for which any morphism $\Delta^1\to\C$ is an equivalence.
\end{remark}
We will also be interested in the $\infty$-category $\S_*$ of {\em pointed spaces},
\index{$\S_*$}
which can be modeled either internally as pointed spaces $\S_*=\S_{\pt/}$ or as the homotopy coherent nerve of the simplicial category of pointed Kan complexes.
\begin{remark}
The coproduct of a pair of objects $X$ and $Y$ of $\S_*$ is computed as the wedge product $X\lor Y$, the quotient of the disjoint union of $X$ and $Y$ obtained by identifying their basepoints.
Via the colimit preserving functor $\S\to\S_*$
which freely adjoins a basepoint, the cartesian product on $\S$ extends to the {\em smash product}
\index{smash product!of spaces}
on $\S_*$.
That is, the smash product sits in a cofiber sequence
\[
X\lor Y\too X\times Y\too X\land Y
\]
and has the property that if $X\simeq X'_+$ and $Y\simeq Y'_+$ then $X\land Y\simeq (X'\times Y')_+$.
\end{remark}
In practice, it is sometimes necessary to be precise about size by bounding various classes of objects by sets of cardinality less than some infinite regular cardinal.
For the purposes of this survey article, however, we will gloss over most of these distinctions and employ a very basic version of the theory of Grothendieck universes:
objects which exist in our first universe will be called {\em small}, objects which exist in the next universe will be called {\em large}, and objects which exist only in a still higher universe will be called {\em very large}.

As one ascends the higher categorical ladder, the size of the mathematical objects under consideration tends to increase.
While a space will always be assumed to be small, an $\infty$-category will typically be assumed to be large, unless otherwise mentioned.
 We write $\Cat_\i$
\index{$\Cat_\i$}
for the large $\i$-category of small $\i$-categories and $\CAT_\i$
\index{$\CAT_\i$}
for the very large $\i$-category of large $\i$-categories.

While some $\infty$-categorical statements will be formulated in the model of quasicategories, we will often employ the theory without reference to any particular model.
From this perspective, a full subcategory is assumed to be closed under equivalences, a cartesian fibration is any functor that is equivalent to a cartesian fibration (see \cite[Section 2.4]{HTT}, or \autoref{sec:GC} for an overview), and ordinary categories are $\infty$-categories with discrete mapping spaces.

\subsection{Acknowledgements}

The modern approach to algebra in higher categorical contexts surveyed in this article is based almost entirely on the groundbreaking work of Jacob Lurie.
It is a great pleasure to thank him for his enormous contributions to the subject.
Collectively, they comprise new and firm foundations for homotopy coherent mathematics.

This article is essentially a summary of the basics of Lurie's {\em Higher Algebra}, with results from {\em Higher Topos Theory} and {\em Spectral Algebraic Geometry} thrown in as necessary, simplified as much as possible in order to be accessible to the neophyte.
Nevertheless, we have tried our best to state definitions and theorems as precisely as possible and provide accurate references for these results (though in the interest of keeping the article as short as possible we have decided not to provide any proofs whatsoever, or even citations for most of the definitions).
Of course, Lurie's work builds on the combined efforts of a great many mathematicians --- far too many to list here, and any attempt to do so will inevitably leave omit many valuable contributions.

It is also a pleasure to thank Benjamin Antieau, Tobias Barthel, Jeremiah Heller, Lars Hesselholt, Achim Krause, Tyler Lawson, Haynes Miller, Thomas Nikolaus, Charles Rezk, Markus Spitzweck, and Hiro Tanaka for helpful remarks and comments on earlier versions of this draft.
The author would also like to thank the Mathematical Sciences Research Institute for providing an extremely pleasant working environment while much of this article was being written, as well as the National Science Foundation for their generous support.

\section{Category theory}\label{sec:ct}

\subsection{Presheaves and colimits}

\begin{definition}
A {\em presheaf} on an $\i$-category $\C$ is a functor $\C^{\op}\to\S$.
The ${\infty}$-category of {\em presheaves} on $\C$ is the functor $\i$-category
\index{presheaf $\infty$-category}
\index{$\P$}
\[
\P(\mathcal{C}) = \mathrm{Fun}(\mathcal{C}^{\op}, \mathcal{S}),
\]
where  ${\mathcal{S}}$ denotes the ${\infty}$-category of (small) spaces.
\end{definition}
\begin{remark}
Recall \cite[Proposition 5.1.3.1]{HTT} that $\P(\C)$ is equipped with a fully faithful Yoneda embedding $j:\C\to\P(\C)$,
\index{Yoneda embedding}
given by the formula $j(A)=\Map_{\C}(-,A):\C^{\op}\to\S$.
Colimits in $\P(\C)$ are computed objectwise as colimits in $\S$.
Since $\S$ is cocomplete, in the sense that it admits all small colimits, we see that $\P(\C)$ is as well.
In fact, the Yoneda embedding exhibits $\P(\C)$ as the free cocompletion of $\C$, in the sense made precise by the following statement.
\end{remark}
\begin{theorem}{\em \cite[Theorem 5.1.5.6]{HTT}}
Given a small $\infty$-category $\C$ and a cocomplete ${\infty}$-category $\D$, precomposition with the Yonda embedding
\[
\Fun^\mathrm{colim}(\P(\C),\D)\subset\Fun(\P(\C),\D)\too\Fun(\C,\D)
\]
induces an equivalence between the $\infty$-categories of functors $\C\to\D$ and colimit preserving functors $\P(\C)\to\D$.
\end{theorem}
\begin{remark}
The actual statement of \cite[Theorem 5.1.5.6]{HTT} is in terms of left adjoint functors $\P(\C)\to\D$.
Left adjoint functors preserve colimits, and in this case any colimit preserving functor is a left adjoint.
\end{remark}

We will require variants of the above construction in which we freely adjoin only certain types of colimits.

\begin{definition}
A map of simplicial sets $J\to I$ is {\em cofinal} if, for every $\i$-category $\C$ and every $I$-indexed colimit diagram $I^\triangleright\to\C$, the induced cone $J^\triangleright\to\C$ is a colimit diagram. See \cite[Proposition 4.1.1.8]{HTT} for details.
\end{definition}

\begin{definition}
Let $\kappa$ be an infinite regular cardinal.
A simplicial set $I$ is {\em $\kappa$-small}\index{$\kappa$-small simplicial set} if $I$ has fewer than $\kappa$ nondegenerate simplices.
A simplicial set $J$ is {\em $\kappa$-filtered} if every map $f:I\to J$ from a $\kappa$-small simplicial set $I$\index{$\kappa$-filtered simplicial set} extends to a cone $f^\triangleright:I^\triangleright\to\I$.
A simplicial set $K$ is {\em sifted}\index{sifted!simplicial set} if $K\to K\times K$ is cofinal.
\end{definition}

\begin{example}
Filtered simplicial sets are sifted, and the simplicial indexing category $\Delta^{\op}$ is sifted.
A functor of $\i$-categories $f:\C\to\D$ preserves sifted colimits if and only if $f$ preserves filtered colimits and geometric realizations.
\end{example}

 \begin{definition}
 An object $A$ of an $\i$-category $\C$ which admits $\kappa$-filtered colimits is {\em $\kappa$-compact}\index{$\kappa$-compact object} if $\Map_\C(A,-):\C\to\S$ commmutes with $\kappa$-filtered colimits.
An object $A$ of an $\i$-category $\C$ which admits geometric realizations is {\em projective}\index{projective object} if $\Map_\C(A,-):\C\to\S$ commutes with geometric realizations.
 \end{definition}

\begin{definition}
An {\em indexing class} is a collection of simplicial sets which we regard as indexing allowed types of diagrams; this is not standard terminology.
Given an indexing class $\I$, we write $\P_\I(\C)\subset\P(\C)$ for the full subcategory consisting of those presheaves $f:\C^{\op}\to\S$ such that, for all $I\in\I$, $f$ transforms $I$-indexed colimits in $\C$ to $I^{\op}$-indexed limits in $\S$.
\end{definition}

\begin{remark}
The indexing classes which most commonly arise when considering algebraic structures are finite discrete, finite, filtered, and sifted.
These classes often come in pairs.
A functor $f:\C\to\D$ of preserves all small colimits if and only if it preserves finite colimits and filtered colimits, or finite discrete colimits (that is, finite coproducts) and sifted colimits.
\end{remark}

\begin{proposition}{\em \cite[Proposition 5.3.6.2]{HTT}}\label{prop:fc}
Let ${\mathcal{C}}$ and $\D$ be ${\infty}$-categories.
If $\C$ admits finite coproducts and ${\mathcal{D}}$ admits sifted colimits, precomposition with the Yoneda embedding $\C\subset\P_\Sigma(\C)$ determines an equivalence of $\i$-categories
\[
\mathrm{Fun}(\mathcal{C}, \mathcal{D})\overset{\simeq}{\from}\mathrm{Fun^{\mathrm{sift}}}(\P_\Sigma(\C), \mathcal{D}) \subset\Fun(\P_\Sigma(\C),\D)
\]
of functors $\C\to\D$ and sifted colimit preserving functors ${\P_\Sigma(\C) \rightarrow \mathcal{D}}$.
If $\C$ admits $\kappa$-small colimits and $\D$ admits $\kappa$-filtered colimits, precomposition with the Yoneda embedding $\C\subset\P_{\kappa\textrm{-}\mathrm{sm}}(\C)$ determines an equivalence of $\i$-categories
\[
\mathrm{Fun}(\mathcal{C}, \mathcal{D})\overset{\simeq}{\from}\mathrm{Fun^{\kappa\textrm{-}\mathrm{filt}}}(\P_{\kappa\textrm{-}\mathrm{sm}}(\C), \mathcal{D}) \subset\Fun(\P_{\kappa\textrm{-}\mathrm{sm}}(\C),\D)
\]
of functors $\C\to\D$ and $\kappa$-filtered colimit preserving functors ${\P_{\kappa\textrm{-}\mathrm{sm}}(\C) \rightarrow \mathcal{D}}$.
\end{proposition}

\subsection{The Grothendieck construction}\label{sec:GC}

Given an ordinary category $\C$, the theory of categories fibered over $\C$ developed in \cite{SGA1} characterizes the essential image of the {\em Grothendieck construction}
\[
\Fun(\C^{\op},\CAT)\too\CAT_{/\C}
\]
as the (not full) subcategory of $\CAT_{/\C}$ consisting of the fibered categories and their morphisms.
Because of the fundamental role the Grothendieck construction plays in higher algebra, we begin with a brief overview of the basic notions of fibered $\i$-category theory \cite{HTT}.
Whenever possible we choose to phrase these notions internally inside of $\CAT_\i$ itself, with a few exceptions: it is sometimes quite useful to represent an $\i$-category by a (marked) simplicial set \cite{HTT}, in particular when specifying diagrams $\D\to\C$ in an $\i$-category $\C$.

Consider the slice $\i$-category $\CAT_{\i/\C}$: objects of $\CAT_{\i/\C}$ are functors $p:\D\to\C$ with target $\C$, and morphisms of $\CAT_{\i/\C}$ from $q:\E\to\C$ to $p:\D\to\C$ are commuting triangles of $\CAT_{\i}$ of the form
\[
\xymatrix{\E\ar[rr]^f\ar[rd]_p & & \D\ar[ld]^q\\
& \C &}.
\]
We write $\Fun_\C(\E,\D)\simeq\Fun(\E,\D)\times_{\Fun(\E,\C)}\{p\}$\index{$\Fun_\C$} for the $\i$-category of functors $f:\E\to\D$ such that $q\circ f = p$.

\begin{definition}
Let $p:\D\to\C$ be a functor\footnote{If we are working externally in the ordinary category of quasicategories, we would additionally require $p$ to be a {\em categorical fibration} \cite{HTT}; internally, however, this is a meaningless assumption, as the notion of categorical fibration is not stable under categorical equivalence.} of $\infty$-categories.\index{fibration!cartesian}
An arrow $\alpha:C\to D$ of $\D$ is {\em $p$-cartesian}
\index{cartesian!arrow}
if the induced functor
\[
\D_{/\alpha}\too\D_{/D}\times_{\C_{/p(D)}} \C_{/p(\alpha)}
\]
is an equivalence of $\i$-categories.
A {\em $p$-cartesian lift} of an arrow $f:A\to B$ in $\C$ is a $p$-cartesian arrow $\alpha:C\to D$ in $\D$ such that $p(\alpha)=f$.
A functor $p:\D\to\C$ is a {\em cartesian fibration}
if, for every arrow $f:A\to B$ of $\C$ and every object $D\in p^{-1}(B)$, there is a $p$-cartesian lift $\alpha:C\to D$ of $f$ with target $D$.
\end{definition}

\begin{remark}
When $p:\D\to\C$ is clear from context, we will refer to $p$-cartesian arrows simply as cartesian arrows.
\end{remark} 
A cartesian fibration $p:\D\to\C$ determines a functor $\mathrm{St}_\C(p):\C^{\op}\to\CAT_\i$, called the {\em straightening}\index{straightening} of $p$ in \cite{HTT}, basically by taking fibers, which happen to be contravariantly functorial in precisely this case.
This operation is inverse to the {\em unstraightening}\index{unstraightening} functor, the $\i$-categorical analogue of the Grothendieck construction,\index{Grothendieck construction} which ``integrates'' the functor $f:\C^{\op}\to\CAT_\i$ to a cartesian fibration $p:\mathrm{Un}_\C(f)\to\C$.
There are analogous versions for cocartesian fibrations,
\index{fibration!cocartesian}
which correspond to covariant functors $\C\to\CAT_\i$.

\begin{theorem}{\em \cite[Theorem 3.2.0.1]{HTT}}
For any small $\infty$-category $\C$, the unstraightening functor
\[
\mathrm{Un}_\C\colon\Fun(\C^{\op},\CAT_\i)\overset{\simeq}{\too}\CAT_{\i/\C}^{\cart}\subset\CAT_{\i/\C}
\]
induces an equivalence between the $\infty$-category $\Fun(\C^{\op},\CAT_\i)$ and the (not full) subcategory $\CAT_{\i/\C}^{\cart}\subset\CAT_{\i/\C}$ consisting of the cartesian fibrations over $\C$ and those functors over $\C$ which preserve cartesian arrows.
\end{theorem}

The straightening and unstraightening equivalence will be used throughout this article, especially in our definitions of (symmetric) monoidal $\infty$-category and (commutative) algebra object therein.
Even more fundamentally, this equivalence is used in the definition of adjunction.

\begin{definition}
Let $\C$ and $\D$ be $\i$-categories.
An {\em adjunction} between $\C$ and $\D$ is cartesian and cocartesian fibration $p:\M\to\Delta^1$ equipped with equivalences $i:\C\to\M_{\{0\}}$ and $j:\D\to\M_{\{1\}}$.
\end{definition}
\begin{remark}
Given an adjunction, the left adoint $f:\C\to\D$ is determined by cocartesian lifts of $0\to 1$, and the right adjoint $g:\D\to\C$ is determined by cartesian lifts of $0\to 1$.
It is possible to construct a unit or counit transformation \cite[Proposition 5.2.2.8]{HTT}, and the fact that both $\C$ and $\D$ fully faithfully embed into $\M$ implies that we have equivalences
\[
\Map(A,gB)\simeq\Map(iA,igB)\simeq\Map(iA,jB)\simeq\Map(jfA,jB)\simeq\Map(fA,B)
\]
for any pair of objects $A$ of $\C$ and $B$ of $\D$.
A functor $f:\Delta^1\to\CAT_\i$ is a {\em left adjoint} if and only if the associated cocartesian fibration $\M\to\Delta^1$ is also cartesian, and a functor $g:\Delta^{1}\to\CAT_\i^{\op}$ is a {\em right adjoint} if and only if the associated cartesian fibration $\M\to\Delta^1$ is also cocartesian.
\end{remark}
\subsection{Monoidal and symmetric monoidal $\infty$-categories}

Informally, a symmetric monoidal $\infty$-category is an $\infty$-category $\C$ equipped with unit $\eta:\Delta^0\to\C$ and multiplication $\mu:\C\times\C\to\C$ maps which are appropriately coherently associative, commutative, and unital.
Making this precise requires a means to organize for us the infinite hierarchy of coherence data necessary to assert that the intermediate multiplications $\C^{\times m}\to\C^{\times n}$, which should be indexed by something like functions from the set with $m$ elements to the set with $n$ elements, are suitably compatible.

\begin{remark}
Actually, as we may first project away some of the factors in the product before we multiply, these operations are indexed by functions of  {\em pointed} sets.
For convenience we allow ourselves to specify a finite pointed set by the cardinality $n$ of its elements, writing $\langle n\rangle=\{0,1,\ldots,n\}$ for the pointed set with basepoint $0$ and $n$ elements in the complement $\langle n\rangle^\circ=\{1,\ldots,n\}$.
\end{remark}
\begin{definition}
A morphism $\alpha:\langle m\rangle\to \langle n\rangle$ in $\Fin_*$ is {\em inert}
\index{inert map!in $\Fin_*$}
if, for each element $i\in\langle n\rangle^\circ$, the preimage $\alpha^{-1}(i)\subset\langle m\rangle^\circ$ consists of a single element.
A morphism $\alpha:\langle m\rangle\to \langle n\rangle$ in $\Fin_*$ is {\em active}
\index{active map}
if the preimage $\alpha^{-1}(\pt)\subset\langle m\rangle$ consists of a single element (necessarily the basepoint).
\end{definition}

\begin{example}
There are exactly $n$ inert maps $\delta_i:\langle n\rangle\to\langle 1\rangle$, namely
\[
\delta_i(j)=\begin{cases} 
      1 & j=i \\
      0 & j\neq i
   \end{cases}.
\]
\end{example}

\begin{definition}
A symmetric monoidal $\infty$-category
\index{symmetric monoidal $\infty$-category}
is a cocartesian fibration $p:\C^\otimes\to\Fin_*$ which satisfies the Segal condition: for each natural number $n$, the map $\C^\otimes_{\n}\to\prod_{i=1}^n\C^\otimes_{\langle 1\rangle}$ induced by the inert morphisms $\delta_i:\n\to\langle 1\rangle$, $1\leq i\leq n$, is an equivalence.
A morphism of symmetric monoidal $\i$-categories, or a {\em symmetric monoidal functor}\index{symmetric monoidal!functor}, from $p:\C^\otimes\to\Fin_*$ to $q:\D^\otimes\to\Fin_*$, is a functor $f:\C^\otimes\to\D^\otimes$ over $\Fin_*$ such that $f$ preserves cocartesian edges.
The $\i$-category $\CMon(\CAT_\i)$\index{$\CMon(\CAT_\i)$} of symmetric monoidal $\i$-categories is the full subcategory of $\CAT_{\i/\Fin_*}^{\cocart}$ spanned by the symmetric monoidal $\i$-categories.
\end{definition}
\begin{remark}
The straightening of a symmetric monoidal $\i$-category $\C^\otimes\to\Fin_*$ is a functor $\Fin_*\to\CAT_\i$ satisfying the Segal condition.
This is the data of a commutative monoid object in $\Cat_\i$, hence the notation $\CMon(\CAT)_\i$.
We could take this as the definition of a symmetric monoidal $\infty$-category, except that in practice these functors are often difficult to write down,
an issue which already arises in ordinary category theory: for instance, the tensor product of modules is more naturally defined by a universal property than a specific choice of representative.
It is usually easier to avoid making explicit choices altogether by constructing the cocartesian fibration $p:\C^\otimes\to\Fin_*$ instead, where one may as well take the fiber over $\langle n\rangle$ to be the $\infty$-category of {\em all} possible choices of the symmetric monoidal product of $n$ objects of $\C$.
\end{remark}
\begin{example}
Let $\C$ be an $\infty$-category with finite products.
Then there exists a cocartesian fibration $p:\C^\times\to\Fin_*$ whose fiber over $\n$ is the $\infty$-category of commutative diagrams of the form $f:(\Sub(n),\leq)^{\op}\to\C$ such that $f(\emptyset)\simeq\pt$ and $f$ carries pushout squares in $\Sub(n)$ to pullback squares in $\C$.
In particular, the map $f(I)\to\prod_{i\in I} f(\{i\})$ is an equivalence in $\C$, from which it follows that $\C^\times_{\n}\to\prod_{i=1}^n\C$ is an equivalence and therefore that $p:\C^\times\to\Fin_*$ is a symmetric monoidal $\infty$-category.
We refer to this as the {\em cartesian} symmetric monoidal structure
\index{cartesian symmetric monoidal structure}
on $\C$, which exists if and only if $\C$ has finite products.
\end{example}
\begin{definition}
Let $p:\C^\otimes\to\Fin_*$ be a symmetric monoidal $\infty$-category.
A {\em commutative algebra object}\index{commutative algebra object}
of $\C^\otimes$ is a section $s:\Fin_*\to\C^\otimes$ of $p$ which sends inert morphisms in $\Fin_*$ to cocartesian morphisms in $\C^\otimes$.
The $\infty$-category $\CAlg(\C^\otimes)$ of commutative algebra objects in $p:\C^\otimes\to\Fin_*$ is the full subcategory of the $\infty$-category $\Fun_{\Fin_*}(\Fin_*,\C^\otimes)$ of sections $s:\Fin_*\to\C^\otimes$ of $p$ consisting of the commutative algebra objects.
\end{definition}

\begin{remark}
There is an equivalence $\CAlg(\CAT_\i^\times)\simeq\CMon(\CAT_\i)$.
That is, the $\i$-category of commutative algebra objects in $\CAT_\i^\times\to\Fin_*$ is equivalent to the $\i$-category of commutative monoid objects in $\CAT_\i$.
\end{remark}

We will also be interested in monoidal $\infty$-categories.
To set up the theory we need a nonsymmetric analogue on the category $\Fin_*$ of finite pointed sets.
\begin{remark}
To keep the prerequisites to a minimum, we purposely avoid using the language of $\i$-operads.
However, from that perspective, a natural choice of indexing category for monoidal structures is the ``desymmetrization'' $q:\Fin_*^\mathrm{ord}\to\Fin_*$ of $\Fin_*$, the functor whose fiber over $\n$ is the set of total orderings of $\n^\circ$.
However it will be both more convenient and elementary to simply use simplicial objects instead.
\end{remark}
\begin{definition}
The functor $\mathrm{Cut}:\Delta^{\op}\to\Fin_*$\index{$\mathrm{Cut}$} is defined by identifying the nonempty finite ordinal $[n]=\{0\to 1\to\cdots\to n-1\to n\}$ with the set of ``cuts'' of the pointed cardinal $\n=\{0,1,\ldots,n-1,n\}$ as follows: we can cut the string $\{0\to 1\to\cdots\to n-1\to n\}$ before or after any $i\in [n]$, for a total of $n+2$ possibilities.
However, cutting before $0$ or after $n$ have the same effect (nothing), so we identify the two trivial cuts with the basepoint of $\n$.
\end{definition}

\begin{definition}
An order preserving function $f:[m]\to[n]$ is {\em inert}
\index{inert map!in $\Delta$ and $\Delta^{\op}$}
if it is it injective with convex image; that is, $f(m)=f(0)+m$.
\end{definition}
\begin{remark}
The $n$ inert morphisms $[1]\to[n]$ in $\Delta$ are the order preserving functions of the form $\sigma_i(j)=i+j$, $0\leq i<n$, $j\in [1]$.
\end{remark}
\begin{definition}
A {\em monoidal $\infty$-category}\index{monoidal $\i$-category} is a cocartesian fibration $p:\C^\otimes\to\Delta^{\op}$ which satisfies the Segal condition: for each $n\in\NN$, the map
\[
\C^\otimes_{[n]}\too\prod_{i=1}^n\C^\otimes_{[1]}
\]
induced by the inert morphisms $\sigma_i:[1]\to[n]$, $0\leq i< n$, is an equivalence.
A morphism of monoidal $\i$-categories, or a {\em monoidal functor}\index{monoidal functor}, from $p:\C^\otimes\to\Fin_*$ to $q:\D^\otimes\to\Fin_*$ is a functor $f:\C^\otimes\to\D^\otimes$ over $\Delta^{\op}$ such that $f$ preserves cocartesian edges.
The $\i$-category $\Mon(\CAT_\i)$\index{$\Mon(\CAT_\i)$} of monoidal $\i$-categories is the full subcategory of $\CAT_{\i/\Delta^{\op}}^{\cocart}$ spanned by the monoidal $\i$-categories.
\end{definition}
\begin{definition}
Let $p:\C^\otimes\to\Delta^{\op}$ be a monoidal $\infty$-category.
An {\em algebra object}\index{algebra object} of $\C^\otimes$ is a section $s:\Delta^{\op}\to\C^\otimes$ of $p$ which sends inert morphisms in $\Delta^{\op}$ to cocartesian morphisms in $\C^\otimes$.
The $\infty$-category $\Alg(\C^\otimes)$ of algebra objects in $\C^\otimes$ is the full subcategory of the $\infty$-category $\Fun_{\Delta^{\op}}(\Delta^{\op},\C^\otimes)$ of sections $s:\Delta^{\op}\to\C^\otimes$ of $p$ consisting of the algebra objects.
\end{definition}
\begin{remark}
As in the symmetric case, there is an equivalence $\Alg(\CAT_\i^\times)\simeq\Mon(\CAT_\i)$: the $\i$-category of associative algebra objects in $\CAT_\i^\times\to\Delta^{\op}$ is equivalent to the $\i$-category of monoid objects in $\CAT_\i$.
\end{remark}

\begin{remark}
Using \cite[Construction 4.1.2.9]{HA} he theory is set up so that a symmetric monoidal $\i$-category $p:\C^\otimes\to\Fin_*$ restricts to a monoidal $\i$-category $q:\C^\otimes|_{\Delta^{\op}}\to\Delta^{\op}$.
It follows that we have a forgetful functor $\CAlg(\C^\otimes)\to\Alg(\C^\otimes):=\Alg(\C^\otimes|_{\Delta^{\op}})$.
\end{remark}

\subsection{Presentable $\infty$-categories}

\begin{definition}
Given a small $\infty$-category $\C$, we write
$
\Ind_\kappa(\C)\subset\P(\C)\index{$\Ind_\kappa$}
$
for the full subcategory consisting of those functors $f:\C^{\op}\to\S$ such that the source $\D$ of the unstraightening $p:\D\to\C$ of $f$ is $\kappa$-filtered.
\end{definition}
\begin{proposition}{\em \cite[Proposition 5.3.5.10]{HTT}}
Let $\C$ be a small $\i$-category and $\D$ an $\i$-category with $\kappa$-filtered colimits. There is an equivalence of $\i$-categories
\[
\Fun^{\kappa\textrm{-}\mathrm{filt}}(\Ind_\kappa(\C),\D)\simeq\Fun(\C,\D),
\]
and hence by \autoref{prop:fc} an equivalence $\Ind_\kappa(\C)\simeq\P_{\kappa\textrm{-}\mathrm{fin}}(\C)$.
\end{proposition}

\begin{definition}
An object $A$ of an $\infty$-category $\C$ is {\em $\kappa$-compact}
\index{$\kappa$-compact object}
if the corepresentable functor $\Map_\C(A,-):\C\to\S$ commutes with $\kappa$-filtered colimits.
We write $\C^\kappa\subset\C$
\index{$\C^\kappa$}
for the full subcategory of $\C$ consisting of the $\kappa$-compact objects.
\end{definition}

\begin{definition}
An $\infty$-category $\C$ is {\em $\kappa$-compactly generated}
\index{$\kappa$-compactly generated}
if $\C$ admits all small colimits and, writing $\C^\kappa\subset\C$ for the full subcategory consisting of the $\kappa$-compact objects, the canonical map $\Ind_\kappa(\C^\kappa)\to\C$ is an equivalence.
A {\em presentable $\infty$-category}
\index{presentable $\infty$-category}
is an $\infty$-category which is $\kappa$-compactly generated for some infinite regular cardinal $\kappa$.
\end{definition}

There are dual notions of morphism of presentable $\infty$-category: namely, those functors which are left (respectively, right) adjoints.
We write $\PrL$
\index{$\PrL$}
\index{$\PrR$}
(respectively, $\PrR$) for these $\infty$-categories.
A key point \cite[Proposition 5.5.3.13 and Theorem 5.5.3.18]{HTT} is that the inclusion of subcategories $\PrL\subset\mathrm{CAT}_\i$ and $\PrR\subset\mathrm{CAT}_\i$ preserves limits: in fact, given a functor $\D\to\PrL$, a cone $\D^\triangleleft\to\PrL$ is limiting in $\PrL$ (respectively, $\PrR$) if and only if the induced cone is limiting in $\mathrm{CAT}_\i$.

\begin{remark}
For each infinite regular cardinal $\kappa$ we have subcategories (again not full) $\PrL_\kappa\subset\PrL$
\index{$\PrL_\kappa$}
\index{$\PrR_\kappa$}
(respectively, $\PrR_\kappa\subset\PrR$) consisting of the $\kappa$-compactly generated $\infty$-categories and those left (respectively, right) adjoint functors which preserve $\kappa$-compact objects (respectively, $\kappa$-filtered colimits).
\end{remark}

\begin{definition}
A {\em presentable fibration}
\index{presentable fibration}
is a cartesian fibration $p:\D\to\C$ such that the straightening $\mathrm{St}_\C(p):\C^{\op}\to\CAT_\i$ factors through the subcategory $\PrR\subset\CAT_\i$ of presentable $\i$-categories and right adjoint functors.
\end{definition}

\begin{remark}\label{rem:aft}
The {\em adjoint functor theorem}
\index{adjoint functor theorem}
\cite[Corollary 5.5.2.9]{HTT} states that any colimit preserving functor $L:\A\to\B$ of presentable $\i$-categories admits a right adjoint $R:\B\to\A$ (this is even true more generally, for instance if $\B$ is only assumed to be cocomplete).
Hence $\PrL$ can be equivalently described as the subcategory of presentable $\i$-categories and left adjoint functors.
Dually, writing $\PrR\subset\CAT_\i$ for the subcategory consisting of the presentable $\i$-categories and right adjoint functors, uniqueness of adjoints allows us to construct a canonical equivalence
$
\PrL^{\op}\simeq\PrR.
$
\end{remark}

\begin{definition}
A presentable symmetric monoidal $\infty$-category is a symmetric monoidal $\infty$-category $p:\C^\otimes\to\Fin_*$
\index{presentable symmetric monoidal $\infty$-category}
such that the underlying $\infty$-category $\C\simeq\C^\otimes_{\langle 1\rangle}$ is presentable and any choice of tensor product bifunctor $\otimes:\C\times\C\to\C$ commutes with colimits separately in each variable.
\end{definition}
\begin{theorem}{\em \cite[Proposition 4.8.1.17]{HA}}
Let $\C$ and $\D$ be presentable $\infty$-categories.
The subfunctor
\[
\Fun'(\C\times\D,-)\subset\Fun(\C\times\D,-)\colon\Prl\too\CAT_\i,
\]
whose value at $\E\in\Prl$ consists of those functors $f:\C\times\D\to\E$ which preserve colimits separately in each variable, is corepresented by an object\index{tensor product!of presentable $\infty$-categories} $\C\otimes\D\in\Prl$.
\end{theorem}
\begin{remark}
It is straightforward to check, with the definition of $\C\otimes\D$ as above, that for a presentable $\infty$-category $\E$, there is a canonical equivalence
\[
\Funl(\C\otimes\D,\E)\simeq\Funl(\C,\Funl(\D,\E)),
\]
where $\Funl(\D,\E)$
\index{$\Funl$}
denotes the $\infty$-category of left adjoint functors $\D\to\E$, which is presentable by \cite{HA}.
Dually, we write $\Funr(\E,\D)$
\index{$\Funr$}
for the $\infty$-category of right adjoint functors $\E\to\D$, and note that $\Funl(\D,\E)\simeq\Funr(\E,\D)^{\op}$.
\end{remark}
\begin{proposition}{\em \cite[Lemma 4.8.1.16]{HA}}
Let $\C,\D$ be presentable $\infty$-categories. Then $\Funr(\D^{\op},\C)$ is a presentable $\i$-category, and 
\[
\C\otimes\D\simeq\Funr(\D^{\op},\C).
\]
\end{proposition}

\begin{theorem}{\em \cite[Corollary 4.8.1.12]{HA}}
The functor $\P\colon\Cat_\i\to\PrL$ extends to a symmetric monoidal functor $\P^\otimes:\Cat_\i^{\times}\to\PrL^\otimes$.
\end{theorem}
\begin{remark}
The symmetric monoidality of the presheaves functor is an $\infty$-categorical generalization of the Day convolution product.
\index{Day convolution product}
If $\C^\otimes$ is a small symmetric monoidal $\infty$-category, then $\P^\otimes(\C)$ inherits a {\em convolution} symmetric monoidal structure \cite[Remark 4.8.1.13]{HA}, given by the formula
\[
(X_1\otimes\cdots\otimes X_n)(A)\simeq\colim_{(B_1,\ldots,B_n)\in\C^{\times n}_{/A}} X_1(B_1)\times\cdots\times X_n(B_n).
\]
Here the colimit is taken over the $\infty$-category $\C^{\times n}_{/A}\simeq\C^{\times n}\times_\C \C_{/A}$.
\end{remark}
\begin{remark}
The fact that $\P^\otimes$ is symmetric monoidal implies that the canonical map $\P(\C)\otimes\P(\D)\to\P(\C\times\D)$ is an equivalence.
\end{remark}

\begin{definition}
A space $X\in\S$ is said to be {\em $n$-truncated} \index{$n$-truncated!space} if 
$
\pi_m(X,x)\cong 0
$
for all  $m>n$ and  $x\in X$.
By convention, a space is said to be $(-2)$-truncated if it is contractible and $(-1)$-truncated if it is empty or contractible.
\end{definition}

\begin{definition}\label{def:ntrunc}
An object $A$ of an $\i$-category $\C$ is said to be $n$-truncated, $n\in\NN$,
\index{$n$-truncated!object}
if the associated representable functor $\Map_\C(-,A):\C^{\op}\to\S$
factors through the full subcategory $\tau_{\leq n}\S\subset\S$ spanned by the $n$-truncated spaces.
\end{definition}
\begin{proposition}{\em \cite[Proposition 5.5.6.18]{HTT}}
Let $\C$ be a presentable $\i$-category.
The inclusion of the full subcategory of $n$-truncated objects $\tau_{\leq n}\C\subset\C$ is stable under limits and admits a left adjoint $\tau_{\leq n}:\C\to\tau_{\leq n}\C$.
\index{$\tau_{\leq n}$}
\index{truncation!functor}
\end{proposition}

\begin{proposition}{\em \cite[Example 4.8.1.22]{HA}}
As an endofunctor of $\PrL$, $\tau_{\leq n}:\PrL\to\PrL$ is idempotent.
It therefore determines a localization of $\PrL$ with essential image the presentable $n$-categories.
In particular, for any presentable $\infty$-category $\C$, we have a canonical equivalence 
$\C\otimes\tau_{\leq n}\S\simeq\tau_{\leq n}\C$.
\end{proposition}

For any presentable $\infty$-category $\C$, the left adjoints of the inclusions $\tau_{\leq m}\C\subset\tau_{\leq n}\C$ for $m<n$ result in a tower
$
\cdots\to\tau_{\leq n}\C\to\cdots\to\tau_{\leq 0}\C
$
of truncations of $\C$ in $\PrL_{\C/}$, and consequently a comparison map
\[
\C\too\lim\{\cdots\to\tau_{\leq n}\C\to\cdots\to\tau_{\leq 0}\C\}.
\]
\begin{definition}
Let $A$ be an object of a presentable $\i$-category $\C$.
The {\em Postnikov tower} of $A$ is the tower of truncations $\cdots\to\tau_{\leq n}A\to\cdots\to\tau_{\leq 0}A$ of $A$, regarded as a diagram in $\C_{A/}$.
\end{definition}
\begin{definition}
Let $\C$ be a presentable $\i$-category.
We say that {\em Postnikov towers converge} in $\C$ if the map $\C\to\lim\{\cdots\to\tau_{\leq n}\C\to\cdots\to\tau_{\leq 0}\C\}$ is an equivalence of $\infty$-categories.
\end{definition}
\index{Postnikov tower}
More concretely, Postnikov towers converge in $\C$ if, for each object $A$, the map $A\to\lim\tau_{\leq n}A$ from $A$ to the limit of its Postnikov tower is an equivalence.

\begin{example}
Postnikov towers converge in the  $\infty$-categories $\S$ and $\S_*$.
\end{example}

\subsection{Stable $\infty$-categories}

\begin{definition}
A {\em zero object}\index{zero object} of an $\infty$-category $\C$ is an object which is both initial and final.
An $\infty$-category $\C$ is {\em pointed}\index{pointed $\i$-category} if $\C$ admits a zero object.
\end{definition} 
\begin{definition}
Let $\C$ be a pointed $\i$-category.
A {\em triangle}\index{triangle} in $\C$ is a commutative square $\Delta^1\times\Delta^1\to\C$ of the form
\[
\xymatrix{
A\ar[r]^f\ar[d] & B\ar[d]^g\\
0\ar[r] & C}
\]
where $0$ is a zero object of $\C$.
We typically write $A\overset{f}{\to}B\overset{g}{\to}C$ for a triangle in $\C$, though it is important to remember that a choice of composite $g\circ f$ and nullhomotopy $g\circ f\simeq 0\in\Map(A,C)$ are also part of the data.
We say that a triangle in $\C$ is {\em left exact}\index{triangle!left exact} if it is cartesian (a pullback), {\em right exact}\index{triangle!right exact} if it is cocartesian (a pushout), and {\em exact}\index{triangle!exact} if it is cartesian and cocartesian.
\end{definition}
\begin{remark}
Let $\C$ be a pointed $\i$-category with finite limits and colimits, and consider a triangle in $\C$ of the form
\[
\xymatrix{
A\ar[r]\ar[d] & 0\ar[d]\\
0\ar[r] & B}.
\]
Then $A\simeq\Omega B$ if the triangle is left exact and $\Sigma A\simeq B$ is the triangle is right exact.
If the triangle is exact, we have equivalences $A\simeq\Omega\Sigma A$ and $\Sigma\Omega B\simeq B$.
\end{remark}
\begin{definition}
A {\em stable $\infty$-category}\index{stable $\infty$-category} is an $\infty$-category $\C$ with a zero object, finite limits and colimits, and which satisfies the following axiom: a commutative square $\Delta^1\times\Delta^1\to\C$ in $\C$ is cartesian if and only if it is cocartesian (in other words, a pullback if and only if it is a pushout).
\end{definition}

\begin{definition}
Let $\C$ and $\D$ be a functor of $\i$-categories which admits finite limits and colimits.
A functor $f:\C\to\D$ is {\em left exact}\index{left exact!functor} if $f$ preserves finite limits, {\em right exact}\index{right exact!functor} if $f$ preserves finite colimits, and {\em exact}\index{exact!functor} if $f$ preserves finite limits and finite colimits.
We write $\CAT_\i^{\st}\subset\CAT_\i$ for the (not full) subcategory consisting of the stable $\i$-categories and the exact functors.
\end{definition}
\begin{proposition}{\em \cite[Corollary 1.4.2.11]{HA}}
Let $\C$ be a pointed $\i$-category.
The following conditions are equivalent:
\begin{enumerate}\itemsep.1em
\item[\emph{(1)}]
$\C$ is stable.
\item[\emph{(2)}]
$\C$ admits finite colimits and $\Sigma:\C\to\C$ is an equivalence.
\item[\emph{(3)}]
$\C$ admits finite limits and $\Omega:\C\to\C$ is an equivalence.
\end{enumerate}
\end{proposition}

\begin{theorem}{\em \cite[Lemma 1.1.2.13]{HA}}
The homotopy category of a stable $\infty$-category $\C$ admits the structure of a triangulated category such that:
\begin{enumerate}\itemsep.1em
\item[\emph{(1)}]
The shift functor is induced from the suspension functor $\Sigma:\C\to\C$.
\item[\emph{(2)}]
A triangle $A\overset{f}{\to} B\overset{g}{\to} C\overset{h}{\to}\Sigma A$ in the homotopy category is exact if and only if there exist exact triangles
\[
\xymatrix{
A\ar[r]^{f'}\ar[d] & B\ar[r]\ar[d]_{g'} & 0\ar[d]\\
0\ar[r] & C\ar[r]^{h'} & D}
\]
in $\C$ such that $f'$, $g'$, and $h'$ are representatives of $f$, $g$ and $h$ composed with the equivalence $\Sigma A\to D$, respectively.
\end{enumerate}\itemsep.1em
\end{theorem}
\begin{definition}
Let $\C$ be an $\infty$-category with finite limits.
A functor
\[
F:\S_*^{\f}\too\C
\]
is said to be {\em excisive}\index{excisive functor} (respectively, {\em reduced}\index{reduced functor}) if $F$ sends cocartesian squares to cartesian squares (respectively, sends initial objects to final objects).
We write
\[
\Exc_*(\S_*^{\f},\C)\subset\Fun(\S_*^{\f},\C)\index{$\mathrm{\Exc_*}$}
\]
for the full subcategory of reduced excisive functors $\S_*^{\f}\to\C$.
\end{definition}
\begin{definition}
Let $\C$ be an $\infty$-category with finite limits.
The $\infty$-category $\Sp(\C)$ of {\em spectrum objects in $\C$} is the $\infty$-category
\[
\Sp(\C)=\Exc_*(\S_*^{\fin},\C)\index{$\Sp(\C)$}\index{spectrum object}
\]
of reduced excisive functors from finite spaces to $\C$.
\index{spectrum object}
\end{definition}
The reason for the terminology ``spectrum object'' is that a reduced excisive functor determines, and is determined by, its value on any infinite sequence $\{S^{n_0}, S^{n_1}, S^{n_2},\ldots\}$ of spheres of strictly increasing dimension.
The most canonical choice is the sequence of all spheres $\{S^0,S^1,S^2,\ldots\}$, so that evaluation on this family of spheres induces a map
\[
\Exc_*(\S_*^{\fin},\C)\too\lim\{\cdots\overset{\Omega}{\to}\C_*\overset{\Omega}{\to}\C_*\overset{\Omega}{\to}\C_*\}
\]
which sends the excisive functor $F$ to the sequence of pointed objects $\{F(S^n)\}_{n\in\NN}$ and equivalences 
$F(S^n)\overset{\simeq}{\too}\Omega F(\Sigma S^n)\overset{\simeq}{\too}\Omega F(S^{n+1})$.
\begin{proposition}{\em \cite[Remark 1.4.2.25]{HA}}
Let $\C$ be an $\i$-category with finite limits. The functor
\[
\Sp(\C)=\Exc_*(\S^{\fin}_*,\C)\too\lim\{\cdots\overset{\Omega}{\to}\C_*\overset{\Omega}{\to}\C_*\overset{\Omega}{\to}\C_*\}
\]
induced by evaluation on spheres is an equivalence of $\i$-categories.
\end{proposition}

\begin{remark}
Since $\C$ has finite limits, it has a final object $\pt$, and so it makes sense to consider the $\infty$-category $\C_*\simeq\C_{\pt/}$ of pointed objects in $\C$.
There is a canonical equivalence $\Sp(\C_*)\simeq\Sp(\C)$.
\end{remark}

The main examples of stable $\infty$-categories are either presentable or embed fully faithfully inside a stable presentable $\infty$-category as the subcategory of $\kappa$-compact objects for some infinite regular cardinal $\kappa$.
Let $\Prl_{\st}\subset\Prl$\index{$\Prl_{\st}$} denote the full subcategory spanned by the stable presentable $\infty$-categories.
\begin{proposition}{\em \cite[Proposition 4.8.2.18]{HA}}
The inclusion of the full subcategory $\Prl_{\st}\subset\Prl$ admits a left adjoint $\Prl\to\Prl_{\st}$.
\end{proposition}
\begin{remark}
The left adjoint ``stabilization'' functor is given by tensoring with the $\i$-category of spectra.
By virtue of the equivalence
$\C\otimes\Sp\simeq\Sp(\C)$, the $\i$-category of spectrum objects in $\C$ is a description of its stabilization.
\end{remark}

\begin{corollary}
Let $\A$ and $\B$ be presentable $\infty$-categories such that $\A$ is stable.
Then $\A\otimes\B$ is stable.
\end{corollary}

\begin{corollary}
The symmetric monoidal structure $\Prl^\otimes$ on $\Prl$ induces a symmetric monoidal structure $\Prl^\otimes_{\st}$\index{$\Prl^\otimes_{\st}$} on the full subcategory $\Prl_\st\subset\Prl$ consisting of the stable presentable $\infty$-categories.
\end{corollary}
We write $\Prl_{\st}^\otimes\subset\Prl^\otimes$ for this symmetric monoidal subcategory.
Note that this inclusion is lax symmetric monoidal and right adjoint to the symmetric monoidal stabilization functor $\Sp(-):\Prl^\otimes\to\Prl^\otimes_{\st}$.
\begin{corollary}
The forgetful functor $\Prl_{\st}\subset\CAT_\i$ extends to a lax symmetric monoidal functor $\Prl_{\st}^\otimes\to\CAT_\i^\times$.
In particular, it carries (commutative) algebra objects of $\Prl_\st^\otimes$ to (symmetric) monoidal $\infty$-categories.
\end{corollary}
\begin{example}
The $\infty$-category $\Sp$ of spectra is a unit of the symmetric monoidal $\infty$-category $\Prl_{\st}^\otimes$ of stable presentable $\infty$-categories.
It therefore admits a presentable symmetric monoidal structure, called the {\em smash product}.
\index{smash product!of spectra}
We will write ${\otimes}^n:\Sp^{\times n}\to\Sp$\index{$\otimes^n$} for any choice of smash product multifunctor.
\end{example}
\begin{remark}
Given spectra $A$ and $B$, their smash product is usually written $A\land B$ is the literature.
We follow the convention of \cite{HA} and write $A\otimes B$ instead, emphasizing the analogy with the tensor product of abelian groups (or more precisely the derived tensor product of chain complexes).
\end{remark} 
\begin{remark}
By construction, the symmetric monoidal structure on the $\infty$-category of spectra is compatible with the cartesian symmetric monoidal structure on the $\infty$-category of spaces via the suspension spectrum functor\index{$\Sigma^\infty_+$}
\[
\Sigma^\infty_+:\S\too\Sp.
\]
This means that, for spaces $X_1,\ldots,X_n$, there is a canonical equivalence
\[
(\Sigma^\infty_+ X_1)\otimes\cdots\otimes(\Sigma^\infty_+ X_n)\simeq\Sigma^\infty_+(X_1\times\cdots\times X_n).
\]
\end{remark}

\begin{proposition}{\em \cite[Corollary 1.4.4.2]{HA}}
A stable $\infty$-category $\C$ is $\kappa$-compactly generated if and only if $\C$ admits small coproducts, a $\kappa$-compact generator for some regular cardinal $\kappa$, and (the homotopy category of) $\C$ is locally small.
\end{proposition}

\begin{definition}
A {\em Verdier sequence}
\index{Verdier sequence}
is a cocartesian square
\[
\xymatrix{
\A\ar[r]\ar[d] & \B\ar[d]\\
0\ar[r] & \C}
\]
in $\CAT_\i^\mathrm{st}$ such that $0$ is a zero object and the top horizontal map $\A\to\B$ is fully faithful.
A {\em semi-orthogonal decomposition}
\index{semi-orthogonal decomposition}
is a Verdier sequence such that the functors $\A\to\B$ and $\B\to\C$ admit right adjoints.
It is common to simply write $\A\to\B\to\C$ for a Verdier sequence, leaving implicit the requirement that $\A\to\B$ be fully faithful with cofiber $\C$.
\end{definition}

\begin{remark} 
Given a Verdier sequence $\A\overset{i}\to\B\overset{j}{\to}\C$, the fact that $i$ is fully faithful implies that the left adjoint functor $i_!:\Ind(\A)\to\Ind(\B)$ is fully faithful with right adjoint $i^*$, and the right adjoint $j^*:\Ind(\C)\to\Ind(\B)$ of the functor $j_!:\Ind(\B)\to\Ind(\C)$ is also fully faithful.
Thus the Verdier sequence $\Ind(\A)\to\Ind(\B)\to\Ind(\C)$ is actually a semi-orthogonal decomposition.
\end{remark}

\subsection{Homotopy groups and t-structures}

Unfortunately, the notion of truncation (see \autoref{def:ntrunc} above) considered earlier is poorly behaved in stable $\i$-categories $\C$, even if $\C$ happens to be presentable.
This is because it is a direct consequence of stability that an object $A$ of $\C$ is $n$-truncated if and only if it is $(n-1)$-truncated, so that the only finitely truncated objects at all are the zero objects, all of which are canonically equivalent to one another.
Thus any attempt to study $\C$ itself via the standard obstruction theoretic truncation type methods is doomed to fail.

The notion of a t-structure on a stable $\i$-category remedies this situation, providing all kinds of other interesting and potentially effective ways of (co)filtering objects of $\C$ --- ideally even $\C$ itself --- as limits of ``Postnikov type'' towers of {\em truncation} functors $\tau_{\leq n}:\C\to\C$, or dually as colimits of ``Whitehead type'' telescopes of {\em connective cover} functors $\tau_{\geq n}:\C\to\C$.
We emphasize that while a t-structure is {\em structure}, it is encoded as innocuously as possible, via a choice of ``orthogonal'' full subcategories of $\C$ in the sense made precise below.

First we need some notation.
If $\C$ is a stable $\infty$-category and $A$ is an object of $\C$, we also write
$A[n]$ for a choice of $n$-fold suspension $\Sigma^n A$ of $A$, $n\in\ZZ$. If $\D\subset\C$ is a
full subcategory of $\C$, we write $\D[n]\subset\C$ for the full subcategory
consisting the objects of the form $B[n]$, where $B$ is an object of the full subcategory $\D$.
Finally, for any pair of objects $A$ and $B$ of $\C$, we write $\Ext^n(A,B)$ for the abelian group $\pi_0\Map_{\C}(A,B[n])$.

\begin{definition}\label{def:t}
A {\em t-structure}\index{t-structure} on a stable $\infty$-category $\C$ consists of a pair of full subcategories $\C_{\geq 0}\subset \C$ and $\C_{\leq 0}\subset \C$ satisfying the following conditions:
\begin{enumerate}\itemsep.1em
\item[(1)]
$\C_{\geq 0}[1]\subset \C_{\geq 0}$ and $\C_{\leq 0}\subset \C_{\leq 0}[1]$;
\item[(2)]
If $A\in \C_{\geq 0}$ and $B\in \C_{\leq 0}$, then $\Map_\C(A,B[-1])=0$;
\item[(3)]   Every $A\in \C$ fits into an exact triangle of the form
$
\tau_{\geq 0}A\to A\to\tau_{\leq -1}A
$
with $\tau_{\geq 0}A\in\C_{\geq 0}$ and $\tau_{\leq -1}A\in \C_{\leq 0}[-1]$.
\end{enumerate}
\end{definition}
\begin{definition}
An exact functor $\C\rightarrow\D$ between stable $\infty$-categories equipped with t-structures is {\em left $t$-exact}\index{left $t$-exact functor} (respectively, {\em right $t$-exact}\index{right $t$-exact functor}) if it sends $\C_{\leq 0}$ to $\D_{\leq 0}$ (respectively, $\C_{\geq 0}$ to $\D_{\geq 0}$). An exact functor is {\em $t$-exact}\index{$t$-exact functor} if is both left and right $t$-exact. We set $\C_{\geq n}=\C_{\geq 0}[n]$ and $\C_{\leq n}=\C_{\leq 0}[n]$.
\end{definition}
\begin{remark}
The data of a t-structure on a stable $\i$-category $\C$ is equivalent to the data of a t-structure on its triangulated homotopy category $\Ho(\C)$ in the sense originally defined and studied by Beilinson-Bernstein-Deligne~\cite{bbd}.
\end{remark}
\begin{example}\label{ex:t}
If $\A$ is a small abelian category, then the bounded derived $\infty$-category $\D^b(\A)$\index{$\D^b(\A)$}  admits a canonical t-structure, where $\D^b(\A)_{\geq n}$ consists of the complexes $A$ such that $\mathrm{H}_i(A)=0$ for $i<n$, and similarly for $\D^b(\A)_{\leq n}$.
If $\A$ is a Grothendieck abelian category, then the {\em unbounded}\index{derived $\infty$-category!unbounded} derived $\infty$-category $\D(\A)$ admits a t-structure with the same description as the previous example.
This stable $\infty$-category and its t-structure are studied in~\cite[Section~1.3.5]{HA}.
\end{example}
\begin{example}\label{ex:t'}
If $R$ is a connective associative ring spectrum, then the stable presentable $\infty$-category $\LMod_R$ of left $R$-module spectra admits a t-structure with $\LMod_{R\,\geq 0}\simeq\LMod_R^{\cn}$, the $\infty$-category of connective left $R$-module spectra.
We call this the {\em Postnikov t-structure}.\index{Postnikov t-structure}
\end{example}

\begin{remark}
The mapping space $\Map_\C(A,B[-1])$ is actually
contractible for $A\in \C_{\geq 0}$ and $B\in \C_{\leq 0}$. This is not
the case for the mapping spectra, as
$\pi_{0}\Map_{\D(\A)}(A[-n],B[-1])\cong\Ext^{n-1}_\A(A,B)$ for any pair of objects $A$ and $B$ in a Grothendieck abelian category $\A$
(see~\cite[Proposition~1.3.5.6]{HA}).
\end{remark}
\begin{proposition}{\em \cite[Remark 1.2.1.12 and Warning 1.2.1.9]{HA}}
The intersection $\C_{\geq 0}\cap \C_{\leq 0}$ is an abelian category equivalent to the full subcategory of $\C_{\geq 0}$ consisting of the discrete objects (that is, $0$-truncated in the sense of \autoref{def:ntrunc}).
\end{proposition}

\begin{definition}
The abelian category $\C_{\geq 0}\cap \C_{\leq 0}$ is called the {\em heart}\index{heart of a t-structure}\index{${}^\heartsuit$} of the t-structure $(\C_{\geq 0},\C_{\leq 0})$ on $\C$, and is denoted $\C^{\heartsuit}$.
\end{definition}

\begin{example}
The hearts of the t-structures in \autoref{ex:t} are both equivalence to the abelian category $\A$ itself.
The heart of the t-structure in \autoref{ex:t'} is $\LMod_{\pi_0R}^{\heartsuit}$, the abelian category of left $\pi_0R$-modules.
\end{example}
\begin{remark}
It turns out that the truncations $\tau_{\geq n}$ and $\tau_{\leq n}$\index{truncation!associated to a t-structure} are functorial in the
sense that the inclusions $\C_{\geq n}\rightarrow \C$ and $\C_{\leq n}\rightarrow
\C$ admit right and left adjoints, respectively, by~\cite[Corollary~1.2.1.6]{HA}.
\end{remark}
\begin{definition}\label{tshg}
Let $\C$ be a stable $\i$-category equipped with a t-structure $(\C_{\geq 0},\C_{\leq 0})$ and let $A$ be an object of $\C$.
There is an equivalence of functors $\tau_{\geq 0}\tau_{\leq 0}\simeq\tau_{\leq 0}\tau_{\geq 0}:\C\to\C^\heartsuit$ \cite[Proposition 1.2.1.10]{HA}.
The {\em homotopy groups}\index{homotopy groups} $\pi_n A$, $n\in\ZZ$, are defined via the formula $\pi_nA=\tau_{\geq 0}\tau_{\leq 0}A[-n]\in \C^{\heartsuit}$.
\end{definition}
\begin{remark}
In the stable setting, passing to homotopy groups is a homological functor\index{homological functor} in the sense that there are long exact sequences
$$\cdots\rightarrow\pi_{n+1}C\rightarrow\pi_nA\rightarrow\pi_nB\rightarrow\pi_nC\rightarrow\pi_{n-1}A\rightarrow\cdots$$
in $\C^\heartsuit$ whenever $A\rightarrow B\rightarrow C$ is an exact triangle in $\C$.
\end{remark}

\begin{definition}
A t-structure $(\C_{\geq 0},\C_{\leq 0})$ on a stable $\infty$-category $\C$ is {\em left complete}\index{left complete t-structure} if the full subcategory of infinitely connected objects $\C_{\geq\infty}=\bigcap_{n\in\ZZ}\C_{\geq n}\subset\C$
consists only of zero objects.
{Right complete} t-structures\index{right complete t-structure} are defined similarly.
\end{definition}

\begin{definition}
A t-structure $(\C_{\geq 0},\C_{\leq 0})$ on a stable $\infty$-category $\C$ is  {\em bounded}\index{bounded t-structure} if the inclusion of the full subcategory of bounded objects $\C^b=\bigcup_{n\in\NN}\C_{\geq -n}\cap \C_{\leq n}\subset \C$ is an equivalence.
There are analogous notions of right and left bounded.
\end{definition}

Many of the most important examples of stable $\infty$-categories may already be familiar from homological algebra.
If ${\mathcal{A}}$ is an abelian category with enough projectives, the derived ${\infty}$-category ${\D^-(\mathcal{A})}$\index{${\D^-(\mathcal{A})}$} of bounded below chain complexes in $\A$ admits a universal property which characterizes the stable $\i$-categories which arise in this way.
We will suppose that $\A$ has enough projective and write $\A^{\proj}\subset\A$ for the full subcategory consisting of the projective objects.

\begin{definition}
Let $\A$ be an abelian category with enough projectives.
The bounded below {\em derived $\infty$-category}\index{derived $\infty$-category}\index{$\D(\A)$} of $\A$ is the homotopy coherent  nerve
$
\D^-(\A):=\mathrm{N}(\mathrm{Ch}^-(\A))
$
of the differential graded category $\mathrm{Ch}^-(\A)$ (viewed as a simplically enriched category as in \autoref{DK}) of bounded below chain complexes in $\A$ (those complexes $A\in\mathrm{Ch}(\A)$ such that $\mathrm{H}_n(A)=0$ for $n\ll 0$.
\end{definition} 

\begin{proposition}{\em \cite[Corollary 1.3.2.18 and Proposition 1.3.2.19]{HA}}
Let $\A$ be an abelian category with enough projectives.
Then $\D^-(\A)$ is a stable $\i$-category, and the full subcategories $\D^-_{\geq 0}(\A)\subset\D^-(\A)$ and $\D^-_{\leq 0}(\A)\subset\D^-(\A)$ consisting of those complexes $A$ such that $\mathrm{H}_n(A)=0$ for $n<0$ and $\mathrm{H}_n(A)=0$ for $n>0$, respectively, comprise a right bounded and left complete t-structure on $\D^-(\A)$.
\end{proposition}

\begin{theorem}{\em \cite[Theorem 1.3.4.4]{HA}}
Let $\A$ be an abelian category with enough projectives and $W\subset\Fun(\Delta^1,\mathrm{Ch}^-(\A))$ the class of quasi-isomorphisms.
There is a canonical equivalence of $\i$-categories $\A[W^{-1}]\simeq\D^-(\A)$.
\end{theorem}

\begin{remark} ${\mathcal{A}}$ sits inside ${\D^-(\mathcal{A})}$ as the full subcategory of complexes concentrated in degree zero.
Moreover, straightforward homological algebra arguments show that $\pi_0: \D^-(\mathcal{A}) \rightarrow \mathcal{A}$ restricts to an equivalence $\D^-(\A)^{\heartsuit}\simeq\A$.
\end{remark}

\begin{remark}
Let $\A$ be a small abelian category with enough projectives.
The category $\Ind(\A)$ of inductive objects of $\A$ is again an abelian category with enough projective objects and is equivalent to the large abelian category $\Fun^{\Pi}((\A^{\proj})^{\op},\Set)$ of product preserving functors from $(\A^{\proj})^{\op}$ to $\Set$.
\end{remark}

\begin{theorem}{\em \cite[Theorem 1.3.3.8]{HA}}
Let $\A$ be an abelian category with enough projectives, let $\C$ be an $\i$-category which admits geometric realizations, and let $\Fun^\mathrm{geom}(\D^-_{\geq 0}(\A),\C)\subset\Fun(\D^-_{\geq 0}(\A),\C)$ denote the full subcategory consisting of those functors which preserve geometric realizations.
Then restriction along the embedding $\A^{\proj}\subset\D_{\geq 0}^-(\A)$ induces an equivalence of $\i$-categories
\[
\Fun(\A^{\proj},\C)\overset{\simeq}{\from}\Fun'(\D^-_{\geq 0}(\A),\C)\subset\Fun(\D^-_{\geq 0}(\A),\C)
\]
between realization preserving functors $\D^-_{\geq 0}(\A)\to\C$ and functors $\A^{\proj}\to\C$.
\end{theorem}
\begin{theorem}{\em \cite[Theorem HA.1.3.3.2]{HA}}
Let $\A$ be an abelian category with enough projectives, let $\C$ be a stable $\infty$-category with a left complete t-structure, and let $\Fun'(\D^-(\A),\C)\subset\Fun(\D^-(\A),\C)$ denote the full subcategory of functors $\D^-(\A)\to\C$ which are right t-exact and send projective objects of $\A$ to $\C^\heartsuit$.
Then restriction along $i:\A\to\D^-(\A)$ induces an equivalence of $\i$-categories
\[
\Fun^\mathrm{rex}(\A,\C^\heartsuit)\overset{\simeq}{\from}\Fun'(\D^-(\A),\C)\subset\Fun(\D^-(\A),\C)
\]
between right t-exact functors $\D^-(\A)\to\C$ which send $\A^{\proj}\subset\A$ to $\C^\heartsuit\subset\C$ and right exact functors $\A\to\C^\heartsuit$.
\end{theorem}

\section{Ring theory}\label{sec:rt}

\subsection{Spectra}

So far we have considered the $\i$-category of spectra from two seemingly different but equivalent perspectives: on the one hand, as reduced excisive functors $\Sp=\Sp(\S)=\Exc_*(\S^{\fin}_*,\S)$, and on the other, as a unit object $\Sp\in\Prl^\otimes_{\st}$ of the symmetric monoidal $\i$-category of presentable stable $\i$-categories.
The former  is more explicit and yields lots of examples, while the latter is more abstract and suggests a universal property.

\begin{remark}
The reader may be wondering what any of this has to do with the classical notion of spectrum in algebraic topology.
By definition, $\Sp=\Exc_*(\S^{\fin}_*,\S)$, but evaluation on the family of spheres $\{S^n\}_{n\in\NN}$ induces an equivalence of $\i$-categories
\[
\Sp\simeq\lim\left\{\cdots\overset{\Omega}{\too}\S_*\overset{\Omega}{\too}\S_*\overset{\Omega}{\too}\S_*\right\}.
\]
\end{remark}

\begin{remark}
Our convention is that we work in the $\infty$-category $\S$ of spaces and its pointed variant $\S_*$, so that all (pointed) spaces have the homotopy type of cell complexes.
If we wish to model a spectrum, however, it is often convenient to work in a category $\C$ equipped with a class of weak equivalences $W$ and functor $\C\to\S_*$ which sends the maps in $W$ to equivalences in $\S_*$, such as the category $\mathrm{Top}_*$ of pointed topological spaces.
In such models we need only ask that the maps $\eta_n:A_n\to\Omega A_{n+1}$ are weak equivalences, a structure which is often referred to as an $\Omega$-spectrum\index{$\Omega$-spectrum} in the literature.
In the case $\C=\mathrm{Top}_*$ we can even require the $\eta_n:A_n\to\Omega A_{n+1}$ to be homeomorphisms.
\end{remark}

\begin{remark}
The forgetful functor $\Prr\to\CAT_\i$ preserves limits, meaning we may also form this limit in $\Prr$. Equivalently, $\Sp$ is given as the colimit
\[
\Sp\simeq\colim\left\{\S_*\overset{\Sigma}{\too}\S_*\overset{\Sigma}{\too}\S_*\overset{\Sigma}{\too}\cdots\right\}
\]
in $\Prl$, by virtue of the antiequivalence $\Prl\simeq\Prr^{\op}$ which takes the adjoint.
\end{remark}

\begin{remark}
The forgetful functor $\Prl\to\CAT_\i$ does not commute with filtered colimits; nevertheless, the $\infty$-category of finite --- or equivalently, in this case, compact --- spectra $\Sp^\mathrm{fin}\simeq\Sp^\omega$ is computed as the filtered colimit of the suspension functor on finite pointed spaces $\S_*^\mathrm{fin}$.
This equivalence
\[
\Sp^\mathrm{fin}\simeq\colim\left\{\S_*^\mathrm{fin}\overset{\Sigma}{\too}\S_*^\mathrm{fin}\overset{\Sigma}{\too}\S_*^\mathrm{fin}\overset{\Sigma}{\too}\cdots\right\}
\]
in $\Cat_\i$ is the $\i$-categorical analogue of the {\em Spanier-Whitehead} category.\index{Spanier-Whitehead $\infty$-category}
\end{remark}

\begin{remark}
We have taken the ``coordinate free'' convention that spectra are reduced excisive functors.
Nevertheless, it is often the case in practice that a spectrum $A$ is given in terms of an infinite delooping of its infinite loop space $\Omega^\infty A$, in which case the associated excisive functor is given by the formula
\[
T_A(-):=\Omega^\infty(-\otimes A)\colon\S_*^\mathrm{fin}\too\S.
\]
Here, $\otimes$ refers to the fact that $\Sp$, as a commutative $\S_*$-algebra in $\Prl$, is canonically left tensored over $\S_*$, and hence over the symmetric monoidal subcategory $\S_*^\mathrm{fin}\subset\S_*$ as well. Alternatively, this tensor is computed by first applying $\Sigma^\infty:\S_*\to\Sp$ and then tensoring with $A$ in $\Sp^\otimes$.
\end{remark}

\begin{remark}
The $\infty$-category $\Sp$ of spectra admits a left and right complete t-structure with heart
$\Sp^\heartsuit\simeq\Ab$ the category of abelian groups.
The identity functor $\Ab\to\Ab$ is right exact and so determines a right t-exact functor $\D^-(\ZZ)\simeq\D^-(\Ab)\to\Sp$ with image the {\em Eilenberg-MacLane spectra}.\index{Eilenberg-MacLane spectrum}
Strictly speaking, these are the {\em generalized} Eilenberg-MacLane spectra:\index{Eilenberg-MacLane spectrum!generalized} an Eilenberg-MacLane spectrum is a generalized Eilenberg-MacLane spectrum which has homotopy concentracted in a single degree, or equivalently is the image of a chain complex with homology concentrated in a single degree.
\end{remark}

\begin{definition}
A {\em cohomology theory}
\index{cohomology theory}
$\{F^n\}_{n\in\ZZ}$ is a $\ZZ$-graded family of functors $F^n:\Ho(\S_*^{\op})\to\Ab$ and natural isomorphisms $\sigma^n:F^n\to F^{n+1}\circ \Sigma$ satisfying the following exactness conditions:
\begin{itemize}\itemsep.1em
\item[\rm{(1)}]
For any cofiber sequence $X\to Y\to Z$ integer $n$ the sequence of abelian groups $F^n(Z)\to F^n(Y)\to F^n(X)$ is exact.
\item[\rm{(2)}]
For any (possibly infinite) wedge decomposition $X\simeq\bigvee_{i\in I} X_i$ and integer $n$, the homomorphism $F^n(X)\to\prod_{i\in I} F^n(X_i)$ is an isomorphism.
\end{itemize}
\end{definition}

\begin{remark}
These are the {\em Eilenberg-Steenrod axioms} \cite{ES52} for a (generalized) cohomology theory.
\index{Eilenberg-Steenrod axioms}
The original formulation included the {\em dimension axiom},
\index{dimension axiom}
which required that $F^n(S^0)=0$ for all $n\neq 0$.
One can show without much difficulty that the cohomology theories $F$ which satisfy the dimension axiom are necessarily of the form $F^n(X)\cong H^n(X;F^0(S^0))$, which is to say cohomology with coefficients in the abelian group $F^0(S^0)$.
This axiom was disregarded as interesting ``generalized'' cohomology theories were discovered.
\end{remark}

\begin{remark}
The {\em Brown representability theorem}
\index{Brown representability theorem}
\cite{Br62}
asserts that any suitably left exact functor $F:\S^{\op}\to\mathrm{Set}_*$ is representable.
In particular,
\[
F^n\cong\pi_0\Map_{\S_*}(-,A_n)\colon\S^{\op}\too\Ab
\]
for some pointed space $A_n$, and the resulting sequence of pointed spaces $\{A_n\}$ can be chosen so that the suspension isomorphisms
\index{suspension isomorphism}
$\sigma^n:F^n\to F^{n+1}\circ\Sigma$ are induced by equivalences $\eta_n:A_n\to\Omega A_{n+1}$ in $\S_*$ (in fact the loop space structure induces the group structure on the represented functor, and any such group structure arises in this way). This is how spectra arose in practice.
\end{remark}

\begin{remark}
By choosing a representing spectrum, any cohomology theory $F=\{F^n\}$ in the classical sense as above gives rise to a {\em cohomological functor}\index{cohomological functor} in the $\i$-categorical sense, by which we mean a limit preserving functor $F:\S_*^{\op}\to\Sp$.
One can show that any such functor is necessarily a right adjoint, so that the $\i$-categories of spectra and cohomological functors are canonically equivalent: 
\[
\Sp\simeq\S_*\otimes\Sp\simeq\Funr(\S_*^{\op},\Sp).
\]
Similarly, $\Sp\simeq\Funr(\Sp^{\op},\Sp)$, so that any cohomological functor of pointed spaces extends uniquely to a cohomological functor of spectra.
\end{remark}
\begin{remark}
A cohomological functor $F:\Sp^{\op}\to\Sp$ induces a functor $\pi_0 F:\Sp^{\op}\to\Ab$ which necessarily factors through the triangulated homotopy category of spectra.
There is a version of Brown representability for triangulated categories which are compactly generated in the appropriate sense, due to Neeman \cite{Nee96}.
This is not a corollary of the corresponding formal result for compactly generated stable $\i$-categories, namely that $\Sp(\C)\subset\Fun(\C^{\op},\Sp)$ is the full subcategory consisting of the cohomological functors.
Indeed, triangulated categories aren't always homotopy categories of stable $\i$-categories and don't necessarily even admit finite limits or colimits; rather, the requisite exactness properties are encoded by the triangulated structure. 
\end{remark}

\begin{remark}\label{exc}
Much of the classical algebraic topology literature models spectra as a full subcategory of local objects inside of a larger category.
Our definition of spectra $\Sp\subset\Fun_*(\S_*^{\fin},\S_*)$ is as the full subcategory of (reduced) excisive functors. The {\em excisive approximation}\index{excisive approximation} \cite[Example 6.1.1.28]{HTT} $\partial F:\S_*^{\fin}\to\S_*$ of a reduced functor $F:\S_*^{\fin}\to\S_*$ is
given by the formula
\[
(\partial F)(T)\simeq\colim_{n\to\infty}\Omega^n F(\Sigma^n T).
\]
\end{remark}
As any finite space admits a cell decomposition, a reduced excisive functor is determined by its values on spheres; in fact, any sequence of spheres $\{S^{n_0},S^{n_1},\ldots\}$ with $n_i>n_j$ whenever $i>j$.
Moreover, for any reduced $F:\S_*^{\fin}\to\S_*$, we have maps
\[
S^1\too\Map_*(S^n,S^{n+1})\too\Map_*(F(S^n),F(S^{n+1}))
\]
and hence maps $S^1\land F(S^n)\to F(S^{n+1})$. This motivates the next definition.

\begin{definition}
A {\em prespectrum}\index{prespectrum} consists of an $\NN$-indexed collection of pointed spaces $\{Z_n\}_{n\in\NN}$ and structure maps $\{\epsilon_{n}:\Sigma Z_{n-1}\to Z_{n}\}_{n\in\NN}$.
\end{definition}

To organize prespectra into an $\infty$-category $\mathrm{PSp}$,\index{$\mathrm{PSp}$} it will be useful to write
\[
\Sigma[1]:\S_*^\NN\too\S_*^\NN\qquad\textrm{and}\qquad\Omega[-1]:\S_*^\NN\too\S_*^\NN
\]
for the ``shifted'' suspension and loops functors $\Sigma[1](X)_n=\Sigma(X_{n-1})$ and $\Omega[-1](X)_n=\Omega(X_{n+1})$, where we take the convention that $X_{-1}\simeq\pt$ is contractible.
Then $\Sigma[1]$ is left adjoint to $\Omega[-1]$, and the $\i$-category of prespectra is defined by forming either of following the pullbacks $\mathrm{CAT}_\i$:
\[
\xymatrix{
\Fun(\Delta^1,\S_*^\NN)\ar[d] & \mathrm{PSp}\ar[r]\ar[d]\ar[l] & \Fun(\Delta^1,\S_*^\NN)\ar[d]\\
\Fun(\partial\Delta^1,\S_*^\NN) & \S_*^\NN\ar[r]^-{(\Sigma[1],\id)}\ar[l]_-{(\id,\Omega[-1])} & \Fun(\partial\Delta^1,\S_*^\NN)}.
\]

\begin{remark}
A prespectrum $Z=\{Z_n\}$ gives rise to a sequence of representable functors
$
\Map(-,Z_n):\S_*^{\op}\to\S_*
$
which collectively represent a cohomology theory if and only if the map $Z\to\Omega[-1] Z$ is an equivalence.
The ``diagonal'' map $\S_*^\NN\to\Fun(\Delta^1,\S_*^\NN)$ identifies $\S_*^\NN$ with the full subcategory of $\Fun(\Delta^1,\S_*^\NN)$ consisting of the equivalences,
and the iterated pullback square 
\[
\xymatrix{
\Sp\ar[r]\ar[d] & \S_*^\NN\ar[d]\\
\mathrm{PSp}\ar[r]\ar[d] & \Fun(\Delta^1,\S_*^\NN)\ar[d]\\
\S_*^\NN\ar[r]^-{(\id,\Omega[-1])} & \Fun(\partial\Delta^1,\S_*^\NN)}
\]
exhibits $\Sp\subset\PSp$ as the equalizer of the endofunctors $\id,\Omega[-1]:\S_*^\NN\to\S_*^\NN$.
\end{remark}
\begin{remark}
This inclusion admits a left adjoint {\em spectrification} functor $\PSp\to\Sp$.
Indeed, if $Z=\{Z_n\}$ is a prespectrum, the $n^\mathrm{th}$-space of the associated spectrum $A$ is given by the formula
\[
A_n\simeq\underset{m\to\i}{\colim}\,\Omega^m Z_{m+n}.
\]
The equivalence $A_n\overset{\!\!\sim}{\to}\Omega A_{n+1}$ is induced by the maps $Z_{n+m}\to\Omega Z_{n+m+1}$, which becomes the equivalence $\colim\Omega^m Z_{n+m}\simeq\Omega^{m+1} Z_{n+m+1}$ after passing to the colimit.
This is essentially the same formula as in \autoref{exc}.
\end{remark}

\begin{remark}
The $\i$-category of prespectra is quite useful in practice.
As a source of examples, a pointed space $X$ evidently determines a suspension presprectrum $\{\Sigma^n X\}_{n\in\NN}$ whose associated spectrum is $\Sigma^\infty X$.
Since the spectrification functor preserves colimits, the $\i$-category of prespectra can be used as a tool for computing colimits and smash products of spectra.
\end{remark}

\subsection{The smash product}

The construction of a symmetric monoidal model category of spectra was a major foundational problem in the subject for quite some time.
One issue is that there isn't an obvious candidate for the smash product of prespectra; instead, given prespectra $A=\{A_m\}$ and $B=\{B_n\}$, their smash product $A\otimes B$ is most naturally indexed on the poset $\NN\times\NN$; that is, $(A\otimes B)_{m,n}=A_m\land B_n$.
Adams' theory of ``handicrafted smash products'' \cite{Adams} shows that any cofinal choice of poset $\NN\subset\NN\times\NN$ results in a prespectrum representing the smash product, and verifies that this procedure descends to a symmetric monoidal structure on the homotopy category of spectra.
However this is insufficient for many purposes, especially as the homotopy category doesn't admit even the most basic sorts of limits and colimits like pullbacks and pushouts.

Another more significant issue is already apparent in homological algebra: in the derived $\i$-category of chain complexes of modules over a commutative ring, the derived tensor product is really only defined up to quasi-isomorphism, so there's no philosophical reason to expect this to lift to a symmetric monoidal structure on the ordinary category of chain complexes.
The first resolutions of this problem in homotopy theory were the {\em symmetric spectra} of Hovey-Shipley-Smith \cite{HSS} and the {\em $S$-modules} of Elmendorf-Kriz-Mandell-May \cite{EKMM}.

There are morphisms of $\otimes$-idempotent objects of $\Prl$
\[
\S\too\S_*\too\CMon(\S)\too\CMon^\mathrm{gp}(\S)\too\Sp,
\]
necessarily symmetric monoidal, which enable us to calculate the tensor product in these presentable $\infty$-categories.
Specifically, writing $\Sigma^\infty:\S_*\to\Sp$\index{$\Sigma^\infty$} for unique symmetric monoidal left adjoint functor, we deduce that
\[
(\Sigma^\infty X_1)\otimes\cdots\otimes(\Sigma^\infty X_n)\simeq \Sigma^\infty(X_1\land\cdots\land X_n)
\]
for any finite collection of pointed spaces $X_1,\ldots,X_n$.
The analogous result remains true for unpointed spaces by addition of a disjoint basepoint.
\begin{remark}
Using the description of spectra as the limit of tower associated to the endofunctor $\Omega:\S_*\to\S_*$, we obtain maps
\[
\Omega^{\infty-n}:\Sp\too\S_*
\]
by projection to the $n^\mathrm{th}$ factor.
The reason for this notation is that we have equivalences $\Omega^\infty\simeq\Omega^n\Omega^{\infty-n}$; indeed, if $A$ is a spectrum, $\Omega^\infty A\simeq A_0\simeq\Omega^n A_n\simeq\Omega^n\Omega^{\infty-n}A$.
The collection of functors $\{\Omega^{\infty-n}\}_{n\in\NN}$, or any infinite subset thereof, form a conservative family of functors $\Sp\to\S_*$ in $\PrR$.
They admit left adjoints $\Sigma^{\infty-n}:\S_*\to\Sp$ which factor through $\PSp\to\Sp$.
\end{remark}
\begin{remark}
By \cite[Proposition 6.3.3.6]{HTT}, any spectrum $A$ admits a canonical presentation by desuspended suspension spectra
\[
A\simeq\colim_{n\to\infty}\Sigma^{\infty-n}\Omega^{\infty-n}A\simeq\colim_{n\to\infty}\Sigma^{-n}\Sigma^\i A_n
\]
in which the maps $\Sigma^{\infty-n} A_n\to\Sigma^{\infty-n-1} A_{n+1}$ correspond to the composites
\[
A_n\to\Omega^{\infty}\Sigma^n\Omega^n\Sigma^{\infty} A_{n}
\to\Omega^{\infty}\Sigma^n\Omega^{n+1}\Sigma^{\infty}A_{n+1}\simeq\Omega^{\infty-n}\Sigma^{\infty-n-1}A_{n+1}.
\]
Using the fact that $\Sigma^{\infty-m}X\otimes\Sigma^{\infty-n}Y\simeq\Sigma^{\infty-m-n}X\land Y$, we can write down explicit formulas for the spaces in the smash product of any finite sequence of spectra.
\end{remark}

\begin{remark}
Another characterization of the $n$-fold smash product functor $\otimes^n:\Sp^{\times n}\to\Sp$ is as the derivative of the $n$-fold cartesian product functor $\times^n:\S^{\times n}\to\S$ or its {\em coreduction}, the $n$-fold smash product functor $\land^n:\S_*^{\times n}\to\S_*$.
See \cite[Example 6.2.3.28]{HA} for further details.
\end{remark}

\begin{remark}
If $B$ is a fixed spectrum, the functor
$
\Sp\overset{(\id,B)}\too\Sp\times\Sp\overset{\otimes^2}\too\Sp
$
which sends $A$ to $A\otimes B$ preserves colimits and therefore, according to the adjoint functor theorem, admits a right adjoint.
This right adjoint is the {\em mapping spectrum} functor $\F(B,-):\Sp\to\Sp$, which admits the following description: if $A$ is a spectrum, then $\F(B,A)$ is the spectrum given by the formula
\[
\F(B,A)_n\simeq\Map_{\Sp}(B,\Sigma^n A),
\]
with structure maps
\[
\Map_{\Sp}(B,\Sigma^n A)\simeq\Map_{\Sp}(B,\Omega\Sigma^{n+1}A)\simeq\Omega\Map_{\Sp}(B,\Sigma^{n+1}A),
\]
where the first equivalence uses the fact that $\Sp$ is a stable $\i$-category and therefore that the composite functor $\Omega\Sigma$ is equivalent to the identity.
\end{remark}

In order to be able to work effectively with spectra, we need ordinary algebraic invariants such as homotopy groups.
One way to obtain the homotopy groups of a spectrum is via the Postnikov t-structure, defined as follows.
\begin{definition}
Let $\Sp_{\leq -1}\subset\Sp$ denote the full subcategory of spectra consisting of those objects $A$ such that $\Omega^\infty A$ is contractible.
\end{definition}

\begin{remark}
The functor $\Omega^\infty\colon\Sp\to\S$ is corepresented by the unit object $\SS$ of $\Sp$, which is compact.
It therefore preserves limits and filtered colimits; furthermore, limits and filtered colimits of contractible spaces are contractible.
Hence the inclusion $\Sp_{\leq -1}\subset\Sp$ preserves limits and filtered colimits and therefore, by the adjoint functor theorem, it admits a left adjoint.
We write $\tau_{\leq -1}\colon\Sp\to\Sp$ for the left adjoint followed by the right adjoint, so that any spectrum $A$ admits a natural unit map $A\to\tau_{\leq -1} A$.
The fiber of the unit map then determines an endofunctor $\tau_{\geq 0}:\Sp\to\Sp$, the {\em connective cover}.\index{connective cover}
We write $\Sp_{\geq 0}$ for the essential image of $\tau_{\geq 0}$ and $\Sp_{\leq 0}=\Sp_{\leq -1}[1]$.
\end{remark}

\begin{proposition}{\em \cite[Proposition 1.4.3.6]{HA}}
The pair of full subcategories $\Sp_{\geq 0}\subset\Sp$ and $\Sp_{\leq 0}\subset\Sp$ are the connective and coconnective parts of a t-structure on $\Sp$.
Moreover, this t-structure is left and right complete, and its heart $\Sp^\heartsuit\simeq\Ab$ is canonically equivalent to the category of abelian groups.
\end{proposition}

\begin{remark}
The spectra that lie in the full subcategory $\Sp_{\geq 0}\subset\Sp$ are called the connective\index{connective!spectrum} spectra, and we will often write $\Sp^{\cn}$ in place of $\Sp_{\geq 0}$.
\end{remark}

\begin{remark}
The homotopy groups of a spectrum $A$, defined via the t-structure on $\Sp$, collectively form a $\ZZ$-graded abelian group
\[
\pi_* A=\bigoplus_{m\in\ZZ}\pi_m A .
\]
Viewing $A$ as a sequence of pointed spaces $\{A_n\}_{n\in\NN}$ equipped with equivalences $A_n\simeq\Omega A_{n+1}$, we can also obtain the homotopy groups of $A$ via the homotopy groups of the pointed spaces $\{A_n\}_{n\in\NN}$.
Specifically, the nonnegative homotopy groups of $A$ are the homotopy groups of the underlying infinite loop space; that is, for $m\geq 0$,
\[
\pi_m A\cong\pi_m A_0\cong\pi_{m+1} A_1\cong\pi_{m+2}A_2\cong\cdots.
\]
The negative homotopy groups of $A$, on the other hand, are the homotopy groups of a sufficiently high delooping; that is, for $m<0$,
\[
\pi_m A\cong\pi_0 A_{-m}\cong\pi_1 A_{-m+1}\cong\pi_2 A_{-m+2}\cong\cdots.
\]
More generally, $\pi_m A\cong\colim\pi_{n} A_{n-m}$, where for $m>n$ we take $A_{n-m}$ to mean the homotopy type of the space $\Omega^{m-n}A_0\simeq\Omega^{m-n+1}A_1\simeq\Omega^{m-n+2}A_2\simeq\cdots$.
\end{remark}

\begin{definition}
Let $A$ and $B$ be spectra.
The $A$-homology of $B$, $A_*(B)$, is the graded abelian group $\pi_*(A\otimes B)$.
Dually, the $A$-cohomology of $B$, $A^*(B)$, is the graded abelian group $\pi_*\F(B,A)$.
\end{definition}

\begin{remark}
If $X$ is a pointed space, we write $A_*(X)=A_*(\Sigma^\infty X)$ and $A^*(X)=A^*(\Sigma^\infty X)$.
If $A$ is a spectrum representing a generalized cohomology theory, this recovers the $A$-cohomology groups of the pointed space $X$.
\end{remark}

\subsection{Associative and commutative algebras}

\begin{definition}
An {\em associative ring spectrum}\index{associative ring spectrum} is an algebra object of the monoidal $\i$-category of spectra.
A {\em commutative ring spectrum}\index{commutative ring spectrum} is a commutative algebra object of the symmetric monoidal $\i$-category of spectra.
\end{definition}

We write $\Alg=\Alg(\Sp^\otimes)$ and $\CAlg=\CAlg(\Sp^\otimes)$ for the $\infty$-categories of associative and commutative ring spectra, respectively.

\begin{remark}
As the notation suggests, we often refer to associative or commutative ring spectra as associative or commutative algebra spectra, or just associative or commutative algebras when the symmetric monoidal $\i$-category of spectra is clear from context.
The latter terminology is especially convenient in relative contexts, such as when we wish to work over\footnote{Over $R$ means over $\Spec R$, reversing the direction of the arrows. Algebraically this means working {\em under} $R$, but the terminology is influenced by geometric intuition.} a fixed base commutative ring spectrum $R$, where the corresponding notion is that of an associative or commutative $R$-algebra (spectrum here is implicit).
\end{remark}

\begin{remark}
Recall that a commutative algebra object of the symmetric monoidal $\i$-category of spectra is a section
\[
A\langle -\rangle:\Fin_*\too\Sp^\otimes
\]
of the cocartesian fibration $p:\Sp^\otimes\to\Fin_*$ such that $A\m\to A\n\in\Sp^\otimes$ is cocartesian whenever $\m\to\n\in\Fin_*$ is inert.
We will typically write $A$ in place of $A\langle 1\rangle\in\Sp^\otimes_{\langle 1\rangle}\simeq\Sp$ and refer to  $A$ as the algebra object.
The value of $A\langle-\rangle$ on the active map $\langle 2\rangle\to\langle 1\rangle$, pushed forward to the fiber over $\langle 1\rangle$ via this map, is the multiplication $\mu:A^{\otimes 2}\to A$.

Similarly, an associative algebra object $A\in\Alg(\Sp^\otimes)$ is a section
\[
A[-]:\Delta^{\op}\too\Sp^\otimes\times_{\Fin_*}{\Delta^{\op}}
\]
of the restricted (along the cut map $\Delta^{\op}\to\Fin_*$) fibration which carries inert arrows to cocartesian arrows.
We will typically also write $A=A[1]$ for the underling spectrum of $A[-]$ and refer to $A$ as the associative ring spectrum.
\end{remark}
\begin{remark}
A section $A[-]:\Delta^\op\to\Sp^\otimes\times_{\Fin_*}\Delta^{\op}$ 
amounts to a diagram of the form
\[
\xymatrix{
A[0]\ar[r] & A[1]\ar@<-.75ex>[l]\ar@<.75ex>[l]\ar@<-.75ex>[r]\ar@<.75ex>[r] & A[2]\ar@<-1.5ex>[l]\ar[l]\ar@<1.5ex>[l]\ar@<-1.5ex>[r]\ar[r]\ar@<1.5ex>[r] & \cdots\ar@<-2.25ex>[l]\ar@<-.75ex>[l]\ar@<.75ex>[l]\ar@<2.25ex>[l]}
\]
in $\Sp^\otimes$.
If it is an algebra object, the $n$ inert maps $[1]\to [n]$ in $\Delta$ force $A[n]$ to decompose as the $n$-fold product $(A[1],\ldots, A[1])$ of $A[1]$ under the equivalence $\Sp^\otimes_{\n}\simeq\Sp^{\times n}$, where we identify $\Sp$ with $\Sp^\otimes_{\langle 1\rangle}$ as usual.
Restricting to the active maps and pushing everything forward to the the fiber over $[1]\in\Delta$ via the inert maps, and writing $\SS$, $A$, $A\otimes A$ for the images of $A[0]$, $A[1]$, $A[2]$,
we obtain a diagram
\[
\xymatrix{
\SS\ar[r] & A\ar@<-.75ex>[r]\ar@<.75ex>[r] & A\otimes A\ar@<-1.5ex>[r]\ar[r]\ar@<1.5ex>[r]
\ar[l] & \cdots\ar@<-.75ex>[l]\ar@<.75ex>[l]}
\]
in $\Sp$, encoding exactly the maps one would expect from an algebra object.
The commutative case is similar but more complex due to the permutations.
\end{remark}
In certain situations, some of which we will encounter later, it is convenient to work with only the connective spectra.
Recall that $\Sp^{\cn}\subset\Sp$ denotes the full subcategory consisting of the spectra $A$ whose negative homotopy groups $\pi_n A$ vanish for all $n<0$.
In other words, $\Sp^{\cn}=\Sp_{\geq 0}$ as full subcategory of spectra, where $\Sp_{\geq 0}$ is defined via the standard Postnikov t-structure.
Since the tensor product of connective spectra is again connective, we actually obtain $\Sp_{\geq 0}^\otimes\subset\Sp^\otimes$ as a symmetric monoidal subcategory.
Thus we have $\i$-categories $\Alg^{\cn}=\Alg(\Sp^{\cn})$ and $\CAlg^{\cn}=\CAlg(\Sp^{\cn})$ of connective associative and commutative algebra spectra.

\begin{proposition}{\em \cite[Proposition 7.1.3.19]{HA}}
Postnikov towers converge in the presentable $\infty$-categories $\Sp^{\cn}$, $\Alg^{\cn}$, and $\CAlg^{\cn}$.
\end{proposition}

\begin{definition}
An {\em $\infty$-monoid}\index{$\i$-monoid} is an algebra object of the symmetric monoidal $\i$-category $\S^\times$ of spaces.
An {\em abelian $\infty$-monoid}\index{abelian $\i$-monoid} is a commutative algebra object of $\S^\times$.
An $\i$-monoid, or abelian $\i$-monoid, $G$ is said to an $\i$-group,\index{$\i$-group} or abelian $\i$-group,\index{abelian $\i$-group} if the ordinary discrete monoid $\pi_0 G$ is a group.
\end{definition}

\begin{remark}
We have $\i$-categories $\Mon_\i=\Mon(\S^\times)=\Alg(\S^\times)$ and $\mathrm{Gp}_\i\subset\Mon_\i$ of $\i$-monoids and $\i$-groups, as well as their abelian variants $\mathrm{AbMon}_\i=\CMon(\S^\times)=\CAlg(\S^\times)$ and $\mathrm{AbGp}_\i\subset\mathrm{AbMon}_\i$.  
The {\em group completion}\index{group completion} of an $\i$-monoid $M$ is the left adjoint of the fully faithful inclusion $\mathrm{Gp}_\i\subset\Mon_\i$, and is given by the formula $G\simeq\Omega\mathrm{B} M$.
\end{remark}

\begin{remark}
The reader might be wondering why we did not define associative and commutative ring spectra as homotopy coherently associative and commutative ring objects in the $\i$-category of spaces.
The primary reason is that this only yields the connective ring spectra, and nonconnective spectra, even ring spectra, such as topological K-theory, are among our most important examples.
Nevertheless it is true, and a good sanity check, that we have an equivalence $\Sp^{\cn}\simeq\mathrm{AbGp}_\i$, which induces equivalences $\Alg^{\cn}\simeq\mathrm{Ring}_\i$ and $\CAlg^{\cn}\simeq\mathrm{CRing}_\i$, where the latter notions are defined as in \cite{GGN} (among other algebraic theories, such as semirings and their $\EE_n$ variants).
\end{remark}

\begin{definition}
A {\em homotopy associative ring spectrum} is an associative algebra object in the monoidal homotopy category of spectra.
A {\em  homotopy commutative ring spectrum} is a commutative algebra object in the symmetric monoidal homotopy category of spectra.
\end{definition}
\begin{remark}
Homotopy categories are ordinary categories.
Hence there are no coherences to specify, and a homotopy associative ring spectrum is the data of a ring spectrum $R$ equipped with a unit map $\eta:\SS\to R$ and a multiplication map $\mu:R\otimes R\to R$ which are associative and unital  in the sense that the diagrams
\[
\xymatrix{
R\otimes R\otimes R\ar[r]^{\mu\otimes\id}\ar[d]_{\id\otimes\mu} & R\otimes R\ar[d]^\mu & & R\ar[r]^{\eta\otimes\id}\ar[d]_{\id\otimes\eta}\ar[rd]^{\id} & R\otimes R\ar[d]^{\mu}\\
R\otimes R\ar[r]^{\mu} & R & & R\otimes R\ar[r]^{\mu} & R}
\]
commute up to homotopy (a choice of homotopy is not part of the data).
As in ordinary algebra, commutativity in this case is a {\em property} of a homotopy associative ring spectrum $R$.
While the theory of homotopy associative or commutative ring spectra plays an important role in algebraic topology, the categories of left (or right) modules over homotopy associative ring spectra are poorly behaved, lacking basic structure such as finite limits and colimits.
\end{remark}

\begin{proposition}{\em \cite[Proposition 7.2.4.27]{HA}}
The $\i$-categories $\Alg^{\cn}$ and $\CAlg^{\cn}$ are compactly generated. Moreover, they are generated under sifted colimits by compact projective\index{projective!object} objects.
\end{proposition}

\begin{definition}
A connective associative ring spectrum $A$ is said to be {\em locally of finite presentation}\index{locally of finite presentation} if $A$ is compact as an object of $\Alg^{\cn}$.
A connective commutative ring spectrum $A$ is said to be {\em locally of finite presentation} if $A$ is a compact as an object of $\CAlg^{\cn}$.
\end{definition}

\subsection{Left and right modules}

An associative (respectively, commutative) algebra $A$ is something that exists in a monoidal (respectively, symmetric monoidal) $\i$-category $\A^\otimes$.
Classically it was common to take $\A^\otimes$ to be the symmetric monoidal category of abelian groups, or (left) $R$-modules for a commutative ring $R$.
Recall that a monoidal $\i$-category $\A^\otimes$ can be regarded as a cocartesian fibration over $\Delta^{\op}$ (equivalently a functor $\Delta^{\op}\to\CAT_\i$)  satisfying the Segal condition, the underlying category $\A$ of $\A^\otimes$ is the value $\A^\otimes_{[1]}$ at the ordinal $[1]$, and the simplicial object $\A^\otimes$ can be regarded as a categorical bar construction on $\A$, using its monoidal structure.

The ordinary bar construction generalizes as follows:  if $\M$ admits a left $\A$-action and $\N$ admits a right $\A$-action, then we may form a simplicial object which in degree $n$ is equivalent to $\N\times\A^{\times n}\times\M$.
The $\infty$-categories $\M$ and $\N$ which arise in this way are said to be left and right tensored over $\A$, respectively; this is the structure which, when given an algebra object $A$ of $\A$, allows us to define the notions of left and right $A$-module object of $\M$ and $\N$.
Such $\i$-categories are themselves cocartesian fibrations over $\A^\otimes$, and hence over $\Delta^{\op}$, but they satisfy a slight variant of the Segal condition.

\begin{definition}\label{def:rsc}
Let $p:\A^\otimes\to\Delta^{\op}$ be a monoidal $\infty$-category.
An $\i$-category {\em left tensored over $\A^\otimes$}\index{left tensored $\i$-category} is a cocartesian fibration $q:\M^\otimes\to\Delta^{\op}$ together with a morphism of cocartesian fibrations $f:\M^\otimes\to\A^\otimes$ over $\Delta^{\op}$ which satisfies the following relative version of the Segal condition: for each natural number $n$, the map
\[
\M^\otimes_{[n]}{\too}\A^{\otimes}_{[n]}\times\M^\otimes_{\{n\}}
\]
induced by $f$ and the inclusion of the final vertex $\{n\}\subset[n]$ is an equivalence.
\end{definition}
\begin{remark}
The morphism of cocartesian fibrations $f:\M^\otimes\to\A^\otimes$ over $\Delta^{\op}$ is an example of a {\em left action object} of $\CAT_\i$.
More precisely, we write $\mathrm{LMon}(\CAT_\i)\subset\Fun(\Delta^1,\CAT_\i{/\Delta^{\op}}^{\cocart})$ for the full subcategory consisting of those morphisms of cocartesian fibrations $f:\M^\otimes\to\A^\otimes$ which satisfy the relative Segal condition of \autoref{def:rsc}.
As in \cite[Notation 4.2.2.5]{HA}, there is a functor $\mathrm{LMon}(\CAT_\i)\to\Mon(\CAT_\i)$ which sends the left action object $f:\M^\otimes\to\A^\otimes$ to the monoidal $\i$-category $\A^\otimes$.
\end{remark}
\begin{remark}
The {\em underlying $\i$-category} of the left tensored $\i$-category $q:\M^\otimes\to\Delta^{\op}$ is the fiber $\M^\otimes_{[0]}$ over $[0]$, which we will typically denote $\M$.
There are equivalences $\M_{\{n\}}\simeq\M$ for all $n$.
The Segal condition at $n=1$ gives an equivalence $\M^\otimes_{[1]}\simeq\A\times\M$ and the inclusion $[0]\to[1]$ of the initial vertex induces a {\em left action} morphism $\A\times\M\to\M$.
\end{remark}
\begin{remark}
Dually, there is an entirely analogous notion of {\em right tensored} $\i$-category\index{right tensored $\i$-category} whose definition instead involves the inclusions of the initial vertices, with the action coming from the inclusion of the final vertices.
\end{remark}

\begin{remark}
A monoidal $\infty$-category $p:\A^\otimes\to\Delta^{\op}$ is canonically left (respectively, right) tensored over itself.
Roughly, restricting the cocartesian fibration $p$ along the right (respectively, left) cone functor $^\triangleright:\Delta\to\Delta$ which sends $[n]$ to $[n]\star[0]$ (respectively, $[n]$ to $[0]\star[n]$) yields a cocartesian fibration $q:\A^{\triangleright\otimes}\to\Delta^{\op}$ whose fiber over $[n]$ is equivalent to $\A^{\times n+1}$. The morphism of cocartesian fibrations $f:\A^{\triangleright\otimes}\to\A^\otimes$ is obtained as in \cite[Example 4.2.2.4]{HA}.
See \cite[Variant 4.2.2.11]{HA} for details and a comparison to the operadic approach.
\end{remark}

\begin{definition}
Let $p:\A^\otimes\to\Delta^{\op}$ be a monoidal $\i$-category and let $f:\M^\otimes\to\A^\otimes$ be an $\i$-category left tensored over $\A^\otimes$. A {\em left module object}\index{left module!object} of $\M^\otimes$ is a map $s:\Delta^{\op}\to\M^\otimes$ such that the composite $f\circ s$ is an algebra object of $\A^\otimes$ and, if $i : [m]\to [n]$ is an inert map in such that $i(m) = n$, $f(i)$ is a cocartesian morphism of $\M^\otimes$.
We write
\[
\LMod(\M^\otimes)\subset\Fun_{\Delta^{\op}}(\Delta^{\op},\M^\otimes)
\]
for the full subcategory consisting of the left module objects of $\M^\otimes$.
\end{definition}
\begin{remark}
In the special case in which $\M^\otimes\simeq\A^{\triangleright\otimes}$ is equivalent to $\A$, regarded as being left tensored over itself via $f:\A^{\triangleright\otimes}\to\A^\otimes$, we simply write $\LMod(\A^\otimes)$ in place of $\LMod(\A^{\triangleright\otimes})$.
\end{remark}
\begin{remark}
There is an evident notion of {\em right module object}\index{right module!object} of an $\infty$-category right tensored over a monoidal $\i$-category.
\end{remark}

\begin{definition}
A {\em left module spectrum}\index{left module!spectrum} is a left module object of the monoidal $\i$-category $\Sp^\otimes\times_{\Fin_*}\Delta^{\op}$ of spectra.
\end{definition}

We write $\LMod=\LMod(\Sp^\otimes)$ for the $\i$-category of left module spectra.

\begin{remark}
The $\i$-category $\LMod$ of left module spectra comes equipped with a forgetful functor $\LMod\to\Alg$ which returns the algebra object in the definition of a left module object.
Given an associative algebra spectrum $A$, we write $\LMod_A$ for the fiber over $A$ of the forgetful functor.
We also have a forgetful functor $\LMod\times_{\Alg}\CAlg\to\CAlg$ obtained by pulling back along the forgetful functor $\CAlg\to\Alg$.
We will sometime simply write $\LMod\to\CAlg$ for this forgetful functor, when it is clear from context that we are only considering left modules over commutative algebra spectra.
\end{remark}

\subsection{Localization}

In this section we briefly review the theory of localization of ring spectra.
Recall that a multiplicatively closed subset $S$, containing the unit, of an associative ring $R$ is said to
satisfy the {\em left Ore condition}\index{left Ore condition} if, for or every pair of elements $r\in R$ and $s\in S$ there exist elements $r'\in R$ and $s'\in S$ such that
$s'r = r's$, and if $r\in R$ is an element such that $rs = 0$ for some element $s\in S$ then there exists an
element $s'\in S$ such that $s'r=0$.

\begin{remark}
If $A$ is an associative algebra spectrum, any homogeneous element $s\in\pi_n A$ can be represented by a ``degree $n$'' left $A$-module map $\Sigma^n A\to A$ which is unique up to homotopy.
\end{remark}
\begin{definition}
Let $A$ be an associative algebra spectrum and $S\subset\pi_* A$ a set of homogeneous elements satisfying the left Ore condition.
A left $A$-module spectrum $M$ is {\em $S$-local}\index{local object} if, for every element $s\in S$, left multiplication by $s$ induces an isomorphism $\pi_* M\to\pi_* M$.
We write $S^{-1}\LMod_A\subset\LMod_A$ for the full subcategory consisting of the $S$-local left $A$-module spectra.
\end{definition}

\begin{proposition}{\em \cite[Remark 7.2.3.18]{HA}}
The inclusion of the full subcategory $S^{-1}\LMod_A\subset\LMod_A$ of $S$-local objects admits a left adjoint $S^{-1}:\LMod_A\to S^{-1}\LMod_A$, the $S$-localization functor.\index{localization!functor}\index{$S^{-1}$}
\end{proposition}

\begin{remark}
The theory of Bousfield localization generalizes that of Ore localization.
In the stable setting, this is the data of a left adjoint functor $L:\B\to\C$ of stable presentable $\i$-categories such that the right adjoint $R:\C\to\B$ is fully faithful.
The kernel of $L:\B\to\C$ is the stable presentable subcategory $\A\subset\B$ consisting of those objects $A\in\B$ such that $L(A)\simeq 0$ in $\C$.
An instance of this construction is the localization of the $\i$-category $\B\simeq\Sp$ of spectra at a fixed spectrum $E$: in this case, a spectrum $N$ is $E$-acyclic if $E\otimes N\simeq 0$ and a spectrum $M$ is $E$-local if every map $N\to M$ from an $E$-acyclic object $M$ is null.
While it turns out that Bousfield localization preserves (commutative) algebra structures \cite[Proposition 2.2.1.9]{HA}, it is quite difficult to control the localization in this generality, which is why we focus on the Ore localization instead.
\end{remark}

\begin{remark}
It is possible to use the left Ore condition to give a reasonable explicit description of the homotopy groups of the left Ore localization $S^{-1}M$, for $M\in\LMod_A$ and $S\subset\pi_* A$ a set of homogeneous elements satisfying the left Ore condition.
See \cite[Construction 7.2.3.19 and Proposition 7.2.3.20]{HA} for details.
\end{remark}

\begin{definition}
Let $A$ be an associative algebra spectrum and $S\subset\pi_* A$ a set of homogeneous elements.
A map of associative algebra spectra $\eta:A\to A'$ exhibits $A'$ as the {\em left Ore localization}\index{localization!left Ore} of $A$ at $S\subset\pi_* A$ if, for every associative algebra spectrum $B$, $\eta^*:\Map_{\Alg}(A',B)\to\Map_{\Alg}(A,B)$ is fully faithful with image those $f:A\to B$ such that $f(s)$ is invertible in $\pi_* B$ for all $s\in S$.
\end{definition}

\begin{theorem}{\em \cite[Proposition 7.2.3.27]{HA}}
Let $A$ be an associative algebra spectrum and $S\subset\pi_* A$ a set of homogeneous elements satisfying the left Ore condition.
The localization $S^{-1} A$ admits the structure of an associative algebra $A[S^{-1}]$ equipped with an algebra map $\eta:A\to A[S^{-1}]$ such that, for any associative algebra spectrum $B$, precomposition with $\eta$ is fully faithful
\[
\eta^*:\Map_{\Alg}(A[S^{-1}],B)\subset\Map_{\Alg}(A,B)
\]
with image those $f:A\to B$ such that $f(s)\in\pi_* B$ is
invertible for all $s\in S$.
\end{theorem}

\begin{remark}
If $S\subset\pi_* A$ is a set of homogeneous elements satisfying the left Ore condition, the canonical map $S^{-1}\pi_*(A)\to\pi_*(S^{-1}A)$ is an isomorphism of graded rings, where $S^{-1}\pi_*A$ denotes the localization as graded rings.
\end{remark}

\begin{remark}
If $A$ is an associative algebra spectrum such that the graded ring $\pi_* A$ is graded commutative, then the left Ore condition on a multiplicative subset $S\subset\pi_* A$ is automatically satisfied.
\end{remark}

There is an analogous statement for commutative algebra spectra.

\begin{definition}
Let $A$ be a commutative algebra spectrum and $S\subset\pi_* A$ a set of homogeneous elements.
A map of commutative algebra spectra $\eta:A\to A'$ exhibits $A'$ as the {\em localization}\index{localization} of $A$ at $S\subset\pi_* A$ if, for every commutative algebra spectrum $B$, $\eta^*:\Map_{\CAlg}(A',B)\to\Map_{\CAlg}(A,B)$ is fully faithful with image those $f:A\to B$ such that $f(s)$ is invertible in $\pi_* B$ for all $s\in S$.
\end{definition}
\begin{theorem}{\em \cite[Example 7.5.0.7]{HA}}
Let $A$ be a commutative algebra spectrum and $S\subset\pi_* A$ a multiplicative set of homogeneous elements.
The localization $S^{-1} A$ admits the structure of a commutative algebra $A[S^{-1}]$ equipped with a commutative algebra map $A\to A[S^{-1}]$ such that, for any commutative algebra specrum $B$, precomposition with $\eta$ is a fully faithful functor
\[
\eta^*:\Map_{\CAlg}(A[S^{-1}],B)\subset\Map_{\CAlg}(A,B)
\]
with image those $f:A\to B$ such that $f(s)\in\pi_* B$ is invertible for all $s\in S$.
\end{theorem}

\begin{remark}
Given a commutative algebra spectrum $A$ and an arbitrary subset $T\subset\pi_* A$ of homogeneous elements of $A$, we often write $A[T^{-1}]$ in place of $A[S^{-1}]$, where $S$ denotes the multiplicative closure of $T$ in $\pi_*A$.
Indexing the elements of $T$ by some ordinal $I$, so that $T=\{t_i\}_{i\in I}$, we have equivalences of commutative algebra spectra $\bigotimes_{i\in I} A[t_i^{-1}]\simeq A[T^{-1}]$, where the infinite tensor product is defined to be the filtered colimit of the finite tensor products.
\end{remark}

A map of ordinary commutative rings $f:A\to B$ is a {\em Zariski localization}\index{localization!Zariski} if there exists a finite collection of elements $x_i\in\pi_0 A$, defining basic Zariski open sets $A\to A[x_i^{-1}]$, such that $f$ is isomorphic to the product map $A\to\Pi_i A[x_i^{-1}]$ as objects of $\CAlg_A$.
There is a similar notion for ring spectra.

\begin{example}\label{eg:zar}
Let $R$ be a commutative ring spectrum and suppose given an element $f\in\pi_0 R$, which we regard as an $R$-linear map $f:R\to R$ via the equivalence $R\simeq\End_R(R)$.
Then the filtered colimit
\[
R[f^{-1}]\simeq\colim\{ R\overset{f}{\too} R\overset{f}{\too} R\overset{f}{\too}\cdots\}
\]
is a $\otimes$-idempotent left $R$-module: that is, the relative tensor product (see \autoref{sec:rtp}) is an equivalence $R[f^{-1}]\otimes_R R[f^{-1}]\simeq R[f^{-1}]$.
If follows that the associated Bousfield localization of $\LMod_R$ is given the formula $M\mapsto M[f^{-1}]\simeq M\otimes_R R[f^{-1}]$, and that $R[f^{-1}]\in\CAlg_R$ has the following universal property: $\Map_{\CAlg_R}(R[f^{-1}],A)$ is contractible if $f\in\pi_0(A)^\times$ and empty otherwise.
\end{example}
\begin{example}\label{eg:zars}
Let $R$ be a commutative ring spectrum and consider the affine scheme $X=\Spec\pi_0 R$. The structure sheaf $\O_{X}$ is determined by its values $\O_{X}(U_f)=\pi_0 R[f^{-1}]$ on the basic open sets $U_f=\Spec\pi_0 R[f^{-1}]$, $f\in\pi_0 R$.
Localizing $R$ at elements of $\pi_0 R$ allows us to enhance $\O_{X}$ to a sheaf of commutative ring spectra $\O_{\Spec R}$ on (the topological space of) $X$ by the formula $\O_{\Spec R}(U_f)=R[f^{-1}]$.
This is the {\em affine spectral scheme}\index{affine spectral scheme} $\Spec R$.
\end{example}

\begin{definition}
Let $R$ be a commutative ring spectrum.
Then $R$ is {\em local}\index{local commutative ring spectrum} if $\pi_0 R$ is local, in the sense that there exists a unique maximal ideal $\mathfrak{m}\subset\pi_0 R$.
Equivalently, $R$ is local if, for any $f\in\pi_0 R$, either $f$ or $1-f$ is a unit.
\end{definition}

\section{Module theory}\label{sec:mt}

\subsection{Monads}

Given an $\infty$-category $\C$, the $\infty$-category $\Fun(\C,\C)$ of endofunctors of $\C$ is canonically a monoidal $\infty$-category with respect to composition.
Viewing $\C$ as a simplicial set and $\Fun(\C,\C)$ as simplicial sets, where composition is already a functor on the nose, we obtain a simplicial model for this monoidal $\infty$-category, which moreover comes equipped with a strict left action on $\C$ via the evaluation pairing $\Fun(\C,\C)\times\C\to\C$.
\begin{definition}
Let $\C$ be an $\infty$-category.
A {\em monad}\index{monad} $T$ on $\C$ is an algebra object of the monoidal $\i$-category $\Fun(\C,\C)$ of endofunctors of $\C$.
\end{definition}
\begin{remark}
We write $T:\C\to\C$ for the underlying endofunctor of the monad.
The unit and multiplication maps are usually denoted $\eta:\id_\C\to T$ and $\mu:T\circ T\to T$. There are also homotopy coherent higher multiplications.
\end{remark}
\begin{remark}
Let $T$ be a monad on an $\infty$-category $\C$.
As $\C$ is left tensored over $\Fun(\C,\C)$ via the evaluation map $\Fun(\C,\C)\times\C\to\C$, it makes sense to consider left $M$-module objects in $\C$ (these are often referred to instead as $M$-algebras, especially in ordinary category theory).
We write $\LMod_T(\C)$ for the $\infty$-category of left $T$-modules.
\end{remark}

\begin{example}
Let $g:\D\to\C$ be a functor of $\infty$-categories which admits a left adjoint $f:\C\to\D$. Then the composite functor $g\circ f:\C\to\C$ admits a canonical structure of a monad on $\C$.
Here, the unit $\eta:\id\to g\circ f$ is the unit of the adjunction, and the multiplication $\mu:g\circ f\circ g\circ f\to g\circ f$ is induced by the counit $\epsilon:f\circ g\to\id$.
The formalism of adjunctions allows one to fill in all the coherences in a unique (up to a contractible space of choices) way.
\end{example}

\begin{remark}
Let $g:\D\to\C$ be a functor of $\infty$-categories which admits a left adjoint $f:\C\to\D$ and let $T=g\circ f$ denote the resulting monad on $\C$, as in the example above.
Then $g:\D\to\C$ factors as the composite $g=p\circ g'$, where $g':\D\to\LMod_{T}(\C)$ is the functor which sends the object $B$ to the left $T$-module $\epsilon:TB\to B$,
and $p:\LMod_{M}(\C)\to\C$ is the forgetful functor.
\end{remark}

\begin{definition}
Let $g:\D\to\C$ be a functor of $\infty$-categories which admits a left adjoint $f:\C\to\D$ and let $T=g\circ f$ be the resulting monad on $\C$.
Then $g$ is {\em monadic over $\C$}\index{monadic functor} if the induced map $g':\D\to\LMod_T(\C)$ is an equivalence.
\end{definition}

Knowing when a functor is monadic is quite important: for instance, monadic functors are conservative, and they exhibit the source as a kind of $\i$-category of modules in the target.
Fortunately there is a recognition principle for monadic functors. To state it we'll need the following notions.
\begin{definition}
The category $\Delta_{-\i}$ has an object $[n]$ for each integer $n\geq -1$ and an arrow $\alpha:[m]\to[n]$ for each nondecreasing function $[m]\cup\{-\i\}\to[n]\cup\{-\i\}$ with $\alpha(-\i)=\i$ which composes in the obvious way.
The subcategory $\Delta_+\subset\Delta_{-\i}$ has the same objects but only those arrows $\alpha:[m]\to[n]$ such that $\alpha^{-1}(-\i)=\{-\i\}$.
Note that we may identify $\Delta$ with the full subcategory of $\Delta_+$ consisting of those objects $[n]$ with $n\geq 0$; in fact, the category $\Delta_+$ parametrizes augmented simplicial objects.
\end{definition} 

\begin{definition}
An augmented simplicial object $A_\bullet:\Delta_+^{\op}\to\C$ is {\em split} if $A_\bullet$ extends to a functor $\Delta_{-\i}^{\op}\to\C$.
A simplicial object $A_\bullet:\Delta^{\op}\to\C$ is {\em split} if it extends to a split augmented simplicial object.\index{split augmented simplicial object}
Finally, given a functor $g:\D\to\C$, an (augmented) simplicial object $A_\bullet$ of $\D$ is $g$-split if $g\circ A_\bullet$ is split as an (augmented) simplicial object of $\C$.
\end{definition}
The monadicity theorem is a higher categorical analogue of the Barr-Beck Theorem.
The result plays a considerably more important role higher categorically due to the difficulty of producing explicit constructions in this context.
\begin{theorem}{\em \cite[Theorem 4.7.3.5]{HA}}\label{thm:bbl}
Let $g:\D\to\C$ be a functor of $\infty$-categories which admits a left adjoint $f:\C\to\D$.
Then $g$ is monadic over $\C$ if and only if $g$ is conservative, $\D$ admits colimits of $g$-split simplicial objects, and $g$ preserves colimits of $g$-split simplicial objects.
\end{theorem}

\subsection{Relative tensor products}\label{sec:rtp}

We now consider the $\infty$-category of left modules over a base commutative ring spectrum $R$ (which could be the sphere itself, in the absolute case).
To generalize ordinary algebra, we'd like to have a notion of (commutative) $R$-algebra spectrum.
In order to make this notion precise, we need a (symmetric) monoidal structure on the $\infty$-category $\LMod_R$ of left $R$-modules.
\begin{remark}
Using the language of $\infty$-operads, these results can be refined to produce $\EE_n$-monoidal $\infty$-categories of left $R$-modules when $R$ is only an $\EE_{n+1}$-algebra spectrum.
See \cite{HA} for the details of this approach.
\end{remark}

As in ordinary algebra, the $\infty$-category of (either left or right) modules for an associative ring spectrum $A$ will not carry a symmetric monoidal structure, which has $A$ as the unit and commutes with colimits in each variable, unless the algebra structure on $A$ extends to commutative algebra structure.
Nevertheless, given a left $A$-module $M$ and a right $A$-module $N$, we may form the iterated tensor products
\[
N\otimes A\otimes\cdots\otimes A\otimes M.
\]
The algebra structure on $A$ and the left and right module structures on $M$ and $N$ organize these into a simplicial spectrum $\mathrm{Bar}_A(M,N)$ with
\[
\mathrm{Bar}_A(M,N)_n\simeq N\otimes A^{\otimes n}\otimes M.
\]
\begin{definition}
The {\em relative tensor product}\index{relative tensor product} $N\otimes_A M$ is a spectrum equivalent to the geometric realization of the simplicial spectrum $\mathrm{Bar}_A(M,N)$:
\[
N\otimes_A M\simeq |\mathrm{Bar}_A(M,N)|.
\]
\end{definition}
\begin{remark}
No noncanonical choices were involved in this construction, and indeed the relative tensor product can be shown to extend to a functor
\[
-\otimes_A -\colon\RMod_A\times\LMod_A\too\Sp
\]
which preserves colimits separately in each variable.
We therefore obtain a morphism $\LMod_{A^{\op}\otimes A}\simeq\RMod_A\otimes\LMod_A\to\Sp$ in $\Prl$.
By Morita theory, such a map is determined by a left $A\otimes A^{\op}$-module, which in this case is $A$ itself.
\end{remark}

Left adjoint functors between $\i$-categories of modules determine, and are determined by, bimodules.
More precisely, if $A$ and $B$ are associative ring spectra, an $(A,B)$-bimodule $M$ determines a functor $\LMod_A\to\LMod_B$ via the relative tensor product, and conversely.
We will avoid the theory of $(A,B)$-bimodules by simply using the equivalent $\infty$-category $\LMod_{A^{\op}\otimes B}$.
\begin{theorem}{\em \cite[Theorem 7.1.2.4]{HA}}
Let $A, B$ be associative algebra spectra.
The relative tensor product induces an equivalence of $\i$-categories
\[
\LMod_{A^{\op}\otimes B}\simeq\Funl(\LMod_A,\LMod_B).
\]
\end{theorem}
\begin{remark}
The above equivalence sends the left $A^{\op}\otimes B$-module $N$ to the functor $M\mapsto N\otimes_A M$ and a left adjoint functor $f:\LMod_A\to\LMod_B$ to $f(A)$, viewed as a left $A^{\op}\otimes B$-module via its right $A$-action.
\end{remark}

\begin{remark}
The right adjoint of the left adjoint functor $N\otimes_A -:\LMod_A\to\LMod_B$ is $\Map_B(N,-):\LMod_B\to\LMod_A$.
\end{remark}

The following version of Morita theory,\index{Morita theory} originally due to Schwede-Shipley \cite{SS}, is a convenient recognition principle for $\infty$-categories of modules.

\begin{theorem}{\em \cite[Theorem 7.1.2.1]{HA}}
Let $\C$ be a stable presentable $\infty$-category and let $P$ be an object of $\C$. Then $\C$ is compactly generated by $P$ if and only if the functor $\Map_{\C}(P,-)\colon  \C\to\RMod_{\End_{\C}(P)}$ is an equivalence.
\end{theorem}

\begin{remark}
If $f:A\to B$ is a map of associative ring spectra, then we may view $B$ as a left $A^{\op}\otimes B$-module, in which case the resulting left adjoint functor $B\otimes_A -:\LMod_A\to\LMod_B$ is the basechange functor, with right adjoint $\Map_B(B,-):\LMod_B\to\LMod_A$ the forgetful functor.
Note that the forgetful functor preserves colimits, as they are detected on underlying spectra, so that this is the same as tensoring with the left $B^{\op}\otimes A$-module $B$, i.e. $\Map_B(B,-)\simeq B\otimes_B -$.
In particular, there is a further right adjoint $\Map_A(B,-):\Mod_A\to\Mod_B$.
\end{remark}

\begin{example}
Suppose that $B\simeq A[S^{-1}]$ is a localization of $A$.
Then the forgetful functor $\LMod_B\to\LMod_A$ is fully faithful with essential image the $S$-local left $A$-module spectra, namely those $M$ such that $M\simeq S^{-1}M$.
\end{example}

\subsection{Projective, perfect, and flat modules}

We now  study the $\infty$-categories $\LMod_A$ for $A$ an associative ring spectrum.

\begin{remark}
The Morita theory of the previous section can be used to identify stable $\i$-categories $\C$ of the form $\LMod_A$ for an associative algebra spectrum $A$, where now $A\simeq\End_\C(P)^{\op}$ for some compact generator $P$ of $\C$.
It follows that $\LMod_A$ is compactly generated.
Moreover, a stable $\infty$-category $\C$ is of the form $\LMod_A$ if and only if $\C$ admits a compact generator, and $\LMod_A\simeq\Ind(\LMod_A^\omega)$ is the $\Ind$-completion of its full subcategory $\LMod_A^\omega$ of compact objects.
\end{remark}
\begin{definition}
Let $A$ be an associative ring spectrum.
A left $A$-module $M$ is {\em perfect}\index{perfect module} if $M$ is a compact object of $\LMod_A$.
\end{definition}

\begin{definition}\label{def:freemod}
Let $A$ be an associative ring spectrum.
A left $A$-module $M$ is {\em free}\index{free!module} if $M$ is a (possibly infinite) coproduct of (unshifted) copies of $A$, viewed as a left module over itself.
\end{definition}

\begin{definition}
A free left $A$-module $M$ is {\em finitely generated}\index{finitely generated!free module} if it is equivalent to a finite coproduct of copies of $A$.
\end{definition}

Just as in ordinary algebra, there are other notions of ``freeness'' corresponding to various forgetful-free adjunctions.
For instance, if $M$ is a left $A$-module and $f:A\to B$ to a ring map, then $B\otimes_A M$ is often referred to as the ``free'' left $B$-module associated to the left $A$-module $M$, even though $B\otimes_A M$ is rarely a free $B$-module (although it will be of course if $M$ was actually a free left $A$-module).

Using the long exact sequence on homotopy groups, it is straightforward to check that a map $f:M\to N$ of connective left $A$-modules (over a connective associative ring $A$) is surjective if and only if $\fib(f)$ is connective.

\begin{definition}
Let $A$ be a connective associative ring spectrum. A left $A$-module $P$ is {\em projective}\index{projective module} if it is connective and projective as object of the $\infty$-category $\LMod^{\cn}_A$ of connective left $A$-modules, in the sense that the corepresented functor $\Map(P,-)\LMod^{\cn}_A\to\S$ preserves geometric realizations (that is, colimits of simplicial diagrams). See \cite[Definition 5.5.8.18]{HTT} for details.
\end{definition}

\begin{remark}
It is unreasonable to ask for a left $A$-module to be projective as an object of the $\infty$-category $\LMod_A$ itself, as the only projective objects of this $\infty$-category are the zero objects.
Indeed, suppose that $M$ is a projective object of $\LMod_A$. For any $A$-module N, we can write the suspension of $N$ as the geometric realization of the simplicial $A$-module
\[
\Sigma N\simeq \left| 0\llarrow N\lllarrow N\oplus N\llllarrow\cdots\right|.
\]
But $
\Map_{\LMod_A}(\Sigma^{-1}M,N)\simeq\Map_{\LMod_A}(M,\Sigma N)\simeq
B\Map_{\LMod_A}(M,N)$,
so $\pi_0 \Map_{\LMod_A}(\Sigma^{-1} M, N)\simeq 0$ for all left $A$-modules $N$ and therefore $M\simeq 0$.
\end{remark}

\begin{proposition}{\em \cite[Proposition 7.2.2.6]{HA}}
Let $\C$ be a stable $\i$-category with a left complete t-structure $(\C_{\geq 0},\C_{\leq 0})$ and let $P\in\C_{\geq 0}$ be an object.
The following conditions are equivalent:
\begin{itemize}\itemsep.1em
\item[\em{(1)}]
$P$ is projective (as an object of $\C_{\geq 0}$).
\item[\em{(2)}]
For every $M\in\C_{\geq 0}$, $\Ext^1_{\C}(P,M)\cong 0$.
\item[\em{(3)}]
For every $M\in\C_{\geq 0}$ and every integer $n > 0$, $\Ext^n_{\C}(P,M)\cong 0$.
\item[\em{(4)}]
For every $M\in\C^\heartsuit$ and every integer $n>0$,  $\Ext^n_{\C}(P,M)\cong 0$.
\item[\em{(5)}]
For every exact triangle $L\to M\to N$ in $\C_{\geq 0}$, the map $\Ext^0_\C(P,M)\to\Ext^0_\C(P,N)$ is surjective.
\end{itemize}
\end{proposition}
\begin{proposition}{\em \cite[Corollary 7.2.2.19]{HA}}\label{prop:projectives}
Let $f:A\to B$ be a map of connective associative algebra spectra such that $\pi_0 f:\pi_0 A\to\pi_0 B$ is an isomorphism.
The basechange functor $f^*=(-)\otimes_A B:\Mod_A\to\Mod_B$ restricts to an equivalence
\[
f^*:\Ho(\LMod_A^{\proj})\overset{\simeq}{\too}\Ho(\LMod_B^{\proj})
\]
on homotopy categories of projective objects.
In particular, the $0$-truncation map $f:A\to\pi_0 A$ induces an equivalence $\Ho(\LMod_A^{\proj})\simeq\LMod_{\pi_0 A}^{\heartsuit\proj}$.
\end{proposition}

\begin{proposition}{\em \cite[Proposition 7.2.2.7]{HA}}
Let $A$ be a connective associative algebra spectrum and $P$ a  projective left $A$-module.
Then there exists a free  left $A$-module $M$ such that $P$ is a retract of $M$.
If additionally $P$ is finitely generated, we may take $M$ to be finitely generated as well.
\end{proposition}
While the notion of projectivity really only makes sense over connective ring spectra, the notion of flatness makes sense over arbitrary ring spectra.

\begin{definition}
Let $A$ be an associative ring spectrum.
A left $A$-module spectrum $M$ is said to be {\em flat}\index{flat module} over $A$ if $\pi_0 M$ is a flat left $\pi_0 A$-module, in the sense of ordinary algebra, and for each integer $n$, the canonical map
\[
\pi_n A\otimes_{\pi_0 A}\pi_0 M\too\pi_n M
\]
is an isomorphism.
\end{definition}

Nevertheless, over a connective associative ring spectrum $A$, the notion of flatness behaves in a manner more similar to that of ordinary algebra. For example we have the following important generalization of Lazard's theorem.

\begin{theorem}{\em \cite[Theorem 7.2.2.15]{HA}}
Let $A$ be a connective associative algebra spectrum and $M$ a connective left $A$-module. The following conditions are equivalent:
\begin{enumerate}
\itemsep.1em
\item[\em{(1)}]
$M$ is flat.
\item[\em{(2)}]
$M$ is a filtered colimit of finitely generated free left $A$-modules.
\item[\em{(3)}]
$M$ is a filtered colimit of finitely generated projective left $A$-modules.
\item[\em{(4)}]
The functor $\RMod_A\to\Sp$ given by $N\mapsto N\otimes_A M$ is left $t$-exact.
\item[\em{(5)}]
If $N$ is a discrete right $A$-module then $N\otimes_A M$ is discrete.
\end{enumerate}
\end{theorem}
\begin{remark}
A free (respectively, projective) left $A$-module $M$ is a filtered colimit of finitely generated free (respectively, projective) left $A$-modules.
Thus the statement of Lazard's theorem remains true (albeit less precise) if we disregard finite generation.
\end{remark}

\begin{remark}
The {\em Tor spectral sequence}\index{spectral sequence} has $E_2$-page
\[
E_2^{p,q} =\mathrm{Tor}^{\pi_*A}_p(\pi_* M,\pi_* N)_q
\]
and converges to the homotopy groups $\pi_{p+q}(M\otimes_A N)$ of the relative tensor product.
If $M$ or $N$ is flat over $A$, $E_2^{p,q}$ vanishes for $p>0$, the spectral sequence degenerates at the $E_2$-page, and
$
\pi_*(M\otimes_A N)\cong\pi_* M\otimes_{\pi_* A}\pi_* N
$
is calculated as graded tensor product of $\pi_* M$ and $\pi_* N$ over $\pi_* A$.
\end{remark}
\begin{remark}
As a consequence we observe that if $M$ and $N$ are both flat over $A$, then their tensor product $M\otimes_A N$ is again
flat over $A$. Since the unit object $A$ of $\LMod_A$ is
flat, we see that the full subcategory $\LMod_A^\flat\subset\LMod_A$ inherits the structure
of a symmetric monoidal $\infty$-category.
\end{remark}

In order to calculate in the $\i$-category $\LMod_A$ of left $A$-modules, it is useful to have something analogous to a projective resolution.
If $A$ is connective, the theory of {\em tor-amplitude} plays this role, giving a means of construct any perfect $A$-module inductively, in a finite number of steps, as an iterated cofiber of maps from shifted finitely projective left $A$-modules.

\begin{definition}
Let $R$ be a connective commutative ring spectrum.
A left  $R$-module   $P$   has   {\em tor-amplitude   contained   in   the   interval   $[a,b]$}\index{tor-amplitude} if  for  any  discrete left $\pi_0  R$-module  $M$, $\pi_i(P\otimes_R M)=0$ for $i\notin [a,b]$. If such integers $a$ and $b$ exist,  $P$ is said to have {\em finite tor-amplitude}.
\end{definition}

If $P$ is an $R$-module, then $P$ has tor-amplitude contained in $[a,b]$ if and only if $P\otimes_R \pi_0
R$  is  a  complex  of  left $\pi_0  R$-modules  with  tor-amplitude  contained  in  $[a,b]$ in
the ordinary sense. Note, however, that this definition  differs  from  that
in~\cite[I  5.2]{SGA6}  as  we  work homologically as opposed to cohomologically.

\begin{proposition}{\em \cite[Proposition 2.13]{AG14}}\label{prop:tor}
Let $A$ be a connective associative algebra spectrum and let $M$ and $N$ be left $A$-module spectra.
\begin{enumerate}
\itemsep.1em
\item[\em{(1)}]  If  $M$  is   perfect  then $M$   has   finite   tor-amplitude.
\item[\em{(2)}]   If $B\in\CAlg_A^{\cn}$ and $M$ has tor-amplitude contained in $[a,b]$ then the left $B$-module $B\otimes_A M$ has tor-amplitude contained in $[a,b]$.
\item[\em{(3)}]  If $M$ and $N$ have tor-amplitude contained in $[a,b]$ then the fiber and cofiber of a map $M\to N$ have tor-amplitude contained in $[a-1,b]$ and $[a,b+1]$.
\item[\em{(4)}]   If $M$ is perfect  with tor-amplitude contained in $[0,b]$ then $M$ is connective and $\pi_0   M\cong\pi_0(\pi_0 A\otimes_A M)$.
\item[\em{(5)}] If $M$ is  perfect  with  tor-amplitude  contained  in  $[a,a]$ then $M$ is equivalent to $\Sigma^{a}P$ for a finitely generated projective left $A$-module $P$.
\item[\em{(6)}]   If $M$ is perfect with tor-amplitude contained in $[a,b]$ then there exists an exact triangle $\Sigma^a P\rightarrow M\to Q$ with $P$ finitely generated projective and $Q$ perfect with tor-amplitude contained in $[a+1,b]$.
\end{enumerate}
\end{proposition}

\begin{remark}
If additionally $A$ is a connective commutative algebra spectrum and $M$ and $N$ are left $A$-modules such that $M$ has tor-amplitude contained in $[a,b]$ and $N$ has tor-amplitude contained in $[c,d]$, then   $M\otimes_A N$ has tor-amplitude contained in $[a+c,b+d]$.
\end{remark}

\subsection{Tensor powers, symmetric powers, and free objects}

Given a map of associative algebra spectra $f:A\to B$, the basechange functor $f^*:\LMod_A\to\LMod_B$ (corresponding to the left $A^{\op}\otimes B$-module $B$) is left adjoint to the forgetful functor $f_*:\LMod_B\to\LMod_A$ (corresponding to the left $B^{\op}\otimes A$-module $B$).
Indeed, the counit of the adjunction $f^*f_*\to\id_{\LMod_B}$ is induced from the left $B$-module action map $B\otimes_A N\to N$, and the unit of the adjunction $\id_{\LMod_A}\to f_*f^*$ is induced from the left $A$-module unit map $M\simeq A\otimes_A M\to B\otimes_A M$.
In particular, $f^*M\simeq B\otimes_A M$ is the {\em free left $B$-module}\index{free!left module} on the left $A$-module $M$, by virtue of the equivalence
\[
\Mod_B(f^*M,N)\simeq\Mod_A(M,f_*N).
\]
\begin{remark}
Observe that $f^*M\simeq B\otimes_A M$ need not be free as a left $B$-module in the sense of \autoref{def:freemod} unless $M$ is free as a left $A$-module.
\end{remark}
\begin{remark}
The forgetful functor $f_*:\LMod_B\to\LMod_A$ is conservative and preserves all small limits and colimits.
Hence it is monadic, and exhibits $\LMod_B$ as the $\i$-category of left modules for the monad $T=f_* f^*$.
\end{remark}

The forgetful functors $\CAlg\to\Alg\to\Sp$ preserve limits, so it is natural to ask whether or not they admit left adjoints.
Using presentability and the adjoint functor theorem (see \autoref{rem:aft}), this is the case if and only if they preserve $\kappa$-filtered colimits for some regular cardinal $\kappa$.
As we might expect from algebraic kinds of categories, they preserve all filtered colimits (as well as geometric realizations), so the left adjoints exist, and are instances of {\em free algebra}\index{free!algebra} functors \cite{HA}.
As is ordinary algebra, it is useful to consider the relative situation, so we work over a base commutative ring spectrum $R$.

\begin{proposition}{\em \cite[Corollaries 3.2.2.4 and 3.2.3.2]{HA}}
The forgetful functors $\CAlg_R\to\Alg_R\to\Sp$ preserve small limits and sifted colimits.
\end{proposition}

Basically by construction, these forgetful functors are also conservative, hence monadic by the Barr-Beck-Lurie \autoref{thm:bbl}.
This means that the free-forgetful adjunctions exhibit $\Alg$ and $\CAlg$ as $\i$-categories of left modules for their respective monads.
The underlying endofunctors of these monads are just in ordinary algebra, the tensor and symmetric algebra constructions.
\begin{remark}
Straightening the symmetric monoidal $\infty$-category $\LMod_R^\otimes$ to a functor $\Fin_*\to\CAT_\i$ and restricting to the active maps $\n\to\langle 1\rangle$ for each $n\in\NN$, we obtain $n$-fold tensor power functors $\mathrm{Ten}_R^n:\LMod_R\to\LMod_R$.
That is, $\mathrm{Ten}_R^n(M)\simeq M^{\otimes n}$,\index{$\mathrm{Ten}_R^n$} where the tensor product is taken over $R$.
\end{remark}
\begin{proposition}{\em \cite[Proposition 4.1.1.18]{HA}}
Let $M$ be a left $R$-module.
The free associative $R$-algebra\index{free!associative $R$-algebra} on $M$ is the tensor\index{tensor algebra}\index{$\mathrm{Ten}_R(M)$} algebra
\[
\mathrm{Ten}_R(M)\simeq\bigoplus_{k\in\NN}\mathrm{Ten}^k_R(M)\simeq\bigoplus_{k\in\NN} M^{\otimes k}.
\]
\end{proposition}
\begin{remark}
This only describes the underlying left $R$-module of the free associative $R$-algebra on $M$.
The multiplication
\[
\bigoplus_{i\in\NN} M^{\otimes i}\underset{R}{\otimes}\,\,\bigoplus_{j\in\NN} M^{\otimes j}\simeq\bigoplus_{i,j\in\NN} M^{\otimes i+j}\too\bigoplus_{k\in\NN} M^{\otimes k}
\]
is given by concatenation of tensor powers.
This is still only a small, but important, piece of the homotopy coherently associative algebra structure. 
\end{remark}

The free commutative algebra functor, also known as the symmetric algebra, uses the symmetric power functors $\Sym^{n}_R:\LMod_R\to\LMod_R$,\index{$\Sym^n_R$} given by the formula
\[
\Sym^n_R(M)\simeq\Ten^n_R(M)_{h\Sigma_n}\simeq M^{\otimes n}_{h\Sigma_n},
\]
where the tensor product is taken in the symmetric monoidal $\infty$-category $\LMod_R^\otimes$ of left $R$-module spectra.
Here the subscript $h\Sigma_n$ refers to the {\em homotopy quotient} of $M^{\otimes n}$ by the permutation action of the symmetric group $\Sigma_n$, which is to say that quotient in the $\infty$-categorical sense.
This is to distinguish from more strict versions of quotients by group actions which arise through either ordinary categorical models of the $\infty$-category $\LMod_R$, where it can make sense to ask for an ordinary categorical quotient, or in the context of equivariant homotopy theory, where there are several notions of fixed points.

\begin{proposition}{\em \cite[Example 3.1.3.14]{HA}}
Let $M$ be a left $R$-module.
The free commutative $R$-algebra\index{free!commutative $R$-algebra} on $M$ is the symmetric\index{symmetric algebra} algebra
\[
\mathrm{Sym}_R(M)\simeq\bigoplus_{k\in\NN}\Sym_R^n(M)\simeq\bigoplus_{k\in\NN} M^{\otimes k}_{h\Sigma_k}.
\]
\end{proposition}
\begin{remark}
While it isn't terribly difficult to describe the multiplication on $\Sym_R(M)$ explicitly, organizing all of the higher multiplications in a coherent manner seems difficult without abstract machinery.
The theory of {\em operadic left Kan extensions} \cite{HA} is one way to make this precise.
\end{remark}

\begin{example}
In the case where $M=R$, we have $\Ten_R(R)\simeq\bigoplus_{n\in\NN} R$.
This is often denoted $R[t]$, as we have that $\pi_*(R[t])\cong(\pi_* R)[t]$ on homotopy groups.
As the notation suggests, the free tensor algebra on one generator in degree zero, $R[t]$, happens to be a commutative $R$-algebra spectrum, though it is {\em not} the free commutative $R$-algebra on one generator in degree zero.
Instead we have that $\Sym_R(R)\simeq\bigoplus_{n\in\NN} \Sym^n_R(R)\simeq \bigoplus_{n\in\NN} R^{\otimes n}_{h\Sigma_n}$,
which is sometimes denoted $R\{t\}$ in order to distinguish it from $R[t]$.
In this case, $\pi_*(\Sym_R(R))\cong\bigoplus_{k\in\NN} R_*(B\Sigma_n)$
is the coproduct of the $R$-homologies of the symmetric groups $\Sigma_n$.
\end{example}

\begin{remark}
If $R$ is a $\QQ$-algebra then, for all $n\in\NN$, the map $R_*(\pt)\to R_*(B\Sigma_n)$ is an isomorphism.
This follows from the Atiyah-Hirzebruch spectral sequence\index{spectral sequence!Atiyah-Hirzebruch} $E^2_{p,q}\cong H_p(B\Sigma_n,\pi_q R)\Rightarrow R_{p+q}(B\Sigma_n)$ and the vanishing of group cohomology in characteristic zero by Maschke's theorem.
Hence the map
\[
R[t]\simeq\Ten_R(R)\too\Sym_R(R)\simeq R\{ t\}
\]
obtained by taking the homotopy quotient is an equivalence.
\end{remark}

\begin{remark}
Away from characteristic zero the map $R_*(\pt)\to R_*(B\Sigma_n)$ is typically not an isomorphism.
In particular, the map $\Ten_R(R)\to\Sym_R(R)$ is rarely an equivalence.
Nevertheless, the fact that $\Ten_\SS(\SS)\simeq\SS[t]$ admits a commutative algebra structure (although it is not free as a commutative algebra) implies, by basechange along the commutative algebra map $\SS\to R$, that there is always a commutative $R$-algebra map $R\{t\}\to R[t]$, which is an equivalence when $R$ is a $\QQ$-algebra and not typically otherwise.
\end{remark}

There are other sorts of free algebra functors as well.
The forgetful functor $\Sp\to\S$ fails to be monadic because it isn't conservative; it is, however, monadic on the full subcategory $\Sp^{\cn}\subset\Sp$ of connective spectra, the obvious subcategory on which it is conservative.
The resulting monad $Q\simeq\Omega^\infty\Sigma^\infty_+$ is a higher categorical analogue of the free abelian group monad.

\begin{definition}
Let $G$ be an $\infty$-group (respectively, $\infty$-monoid).
The {\em group ring}\index{group ring} (respectively, {\em monoid ring}\index{monoid ring}) $R[G]$ is the associative ring spectrum $R[G]\simeq R\otimes\Sigma^\infty_+ G$, with $R$-algebra structure induced from that of $G$ via the symmetric monoidal functor $\Sigma^\infty_+:\S\to\Sp$.
\end{definition}

\begin{remark}
Group and monoid rings are a rich source of examples of ring spectra.
For instance, toric varieties, which are locally modeled on the group and monoid rings like $R[\ZZ^{\times k}]$ and $R[\NN^{\times k}]$, are combinatorial enough that they descend from the integers to the sphere, giving basic examples of spectral schemes such as the {\em projective space} $\PP^n_R$~\cite[Construction 5.4.1.3]{SAG}.
\end{remark}

\begin{definition}
The {\em general linear group}\index{general linear group} $\GL_n(R)$ of a ring spectrum $R$ is the $\i$-group $\Aut_R(R^{\oplus n})$ of left $R$-module automorphisms of $R^{\oplus n}$.
\end{definition}

\begin{example}
Let $M=\coprod_{n\in\NN} B\Sigma_n$, the free abelian $\i$-monoid on one generator.
Then $\SS[M]\simeq\bigoplus\Sigma^\infty_+ B\Sigma_n$ is the free commutative algebra spectrum on one generator in the sense that, via the equivalences
\[
\Map_{\CAlg}(\SS[M],A)\simeq\Map_{\mathrm{AbMon}_\infty}(M,\Omega^\infty A)\simeq\Map_\S(\pt,\Omega^\infty A)\simeq\Omega^\infty A,
\]
specifying a commutative algebra map $\SS[M]\to A$ is the same as specifying a point of $\Omega^\infty A$, the image of the generator $\pt\simeq B\Sigma_1\to M$ of $M$.
\end{example}

\begin{example}
Let $G\simeq\Omega^\infty\SS$ be the abelian $\infty$-group given by the infinite loop space of the sphere $\SS$.
Under the equivalence between the $\infty$-categories of connective spectra and abelian $\i$-groups, it follows that $G$ is the free abelian $\i$-group on one generator, or equivalently the group completion of the free abelian $\i$-monoid $M\simeq\coprod_{n\in\NN} B\Sigma_n$ on one generator of the previous example.

By adjunction, for any commutative algebra spectrum $A$, we have equivalences
\[
\Map(\SS[G],A)\simeq\Map(G,\Omega^\infty A)\simeq\GL_1(A)\subset\Omega^\infty A,
\]
where the $\infty$-group $\GL_1(A)\simeq\Aut_{\Mod_A}(A)$ is equivalently the subspace of the $\infty$-monoid $\Omega^\infty A\simeq\End_{\Mod_A}(A)$ consisting of the invertible components; that is, it fits into the pullback square
\[
\xymatrix{
\GL_1(A)\ar[r]\ar[d] & \Omega^\infty A\ar[d]\\
(\pi_0 A)^\times\ar[r] & \pi_0 A.
}
\]
Since $\SS[G]$ corepresents the functor which sends the commutative algebra $A$ to its space of units, the ``derived scheme'' $\Spec\SS[G]$ can be regarded as derived version of the multiplicative group scheme.
\end{example}

\begin{example}
Any abelian $\i$-group $G$ with $\pi_0 G\cong\ZZ$ determines a ``connective spectral abelian group scheme'' with underlying ordinary scheme $\GG_m$, the multiplicative group.
While the previous example represents the functor of units, it is sometimes necessary to consider more ``strictly commutative'' derived analogues of $\GG_m$. At the other extreme, one can consider $\ZZ$ as an $\i$-group, in which case the spherical group ring $\SS[\ZZ]$ is a flat commutative $\SS$-algebra such that, for any commutative ring spectrum $A$,
\[
\Map_{\CAlg}(\SS[\ZZ],A)\simeq\Map_{\mathrm{AbMon}_\infty}(\ZZ,\Omega^\infty A)\simeq\Map_{\mathrm{AbGp}_\i}(\ZZ,\GL_1(A)).
\]
This space of ``strict units'' of $A$ plays a central role in elliptic cohomology \cite{EC1}.
\end{example}

Thom spectra\index{Thom spectra} (\cite{Thom}, \cite{MQRT}, \cite{ABGHR}) are spectra which occur as quotients of the sphere  by the action of a group.
Again, it is somewhat more useful to have the version relative to a fixed  commutative ring spectrum $R$.
By construction the $\i$-group $\GL_1(R)$ is the universal group which acts on $R$ by $R$-linear maps, so any $\i$-group $G$ over $GL_1(R)$ acts as well.
The group homomorphism $f:G\to\GL_1(R)$ deloops to a map of pointed spaces $f:\mathrm{B}G\to \mathrm{BGL}_1(R)$.

\begin{example}
The equivalence $\GL_1(\SS)\simeq\colim_{n\to\infty}\Aut_*(S^n)$ and compatible families of $\i$-group maps $O(n)\to\Aut_*(S^n)$ and $U(n)\to\Aut_*(S^{2n})$ give $\i$-group maps $O\to\GL_1(\SS)$ and $U\to\GL_1(\SS)$.
The Thom spectra $MO$ and $MU$ are the homotopy quotients $MO=\SS_{hO}$ and $MU=\SS_{hU}$.
While defined in an apparently abstract and formal way, there's a surprising connection to geometry: $\pi_* MO$ is  the bordism ring of unoriented manifolds, $\pi_*MU$ is the bordism ring of stably almost complex manifolds.
The analogues for other tangential structures hold; for instance, $\pi_*\SS$ is the stably framed bordism ring.
\end{example}

\subsection{Smooth, proper, and dualizable objects}

\begin{definition}
A symmetric monoidal $\i$-category $p:\C^\otimes\to\Fin_*$ is {\em closed}\index{closed!symmetric monoidal $\i$-category} if, for any object $A$ of $\C=\C^{\otimes}_{\langle 1\rangle}$, the right multiplication by $A$ functor $(-)\otimes A:\C\to\C$ admits a right adjoint $\F_\C(A,-):\C\to\C$.
\end{definition}
\begin{remark}
This allows for the construction of {\em function objects}\index{function object} $\F_\C(B,C)$ of $\C$, naturally as a functor $\F_\C:\C^{\op}\times\C\to\C$.
In particular, there are natural equivalences
$
\Map_\C(A\otimes B,C)\simeq\Map_\C(A,\F_\C(B,C)).
$
\end{remark}
\begin{definition}
Let $\C^\otimes$ be a closed symmetric monoidal $\i$-category.
A {\em dual}\index{dual} of an $A$ is an object of the form $\F_\C(A,\mathbf{1})$, where $\mathbf{1}$ denotes a unit object.
\end{definition}
As duals are uniquely determined, we will write $\Dual_\C A$, or $\Dual A$, for a dual of $A$.
\begin{remark}
There is a canonical {\em evaluation}\index{evaluation map} map
\[
\epsilon:A\otimes\Dual_\C A\simeq A\otimes\F_\C(A,\mathbf{1})\too\mathbf{1},
\]
any map corresponding to the identity of $\F_\C(A,\mathbf{1})$ under the equivalence $\Map(A\otimes\F_\C(A,\mathbf{1}),\mathbf{1})\simeq\Map(\F_\C(A,\mathbf{1}),\F_\C(A,\mathbf{1}))$.
\end{remark}
\begin{definition}
An object $A$ of a closed symmetric monoidal $\i$-category $\C^\otimes$ is {\em dualizable}\index{dualizable object} if there exists
a {\em coevaluation}\index{coevaluation map} map $\eta: \mathbf{1}\too\Dual A\otimes A$
such that the compositions
\begin{gather*}
A\xrightarrow{A\otimes\eta} A\otimes\Dual A\otimes A\xrightarrow{\epsilon\otimes A}A\\
\Dual A\xrightarrow{\eta\otimes\Dual A}\Dual A\otimes A\otimes\Dual A\xrightarrow{\Dual A\otimes\epsilon}\Dual A
\end{gather*}
are equivalent to the identity.
\end{definition}

We write $\C^{\mathrm{dual}}\subset\C$ for the full subcategory consisting of the dualizable objects of $\C^\otimes$.
It always contains at least one object, any unit object $\mathbf{1}$.

\begin{example}
Let $R$ be a commutative ring spectrum and consider the closed symmetric monoidal $\i$-category $\LMod_R^\otimes$ of left $R$-module spectra.
Stability forces the full subcategory $\LMod_R^{\mathrm{dual}}$ of dualizable objects to be closed under finite limits, colimits, and retracts \cite[Theorem 2.1.3]{HPS}.
Since $\LMod_R$ is compactly generated, it follows that the full subcategories of compact and dualizable objects coincide.
\end{example}

Given a commutative algebra object $\T$ of $\Prl$, we may form the $\i$-category $\LMod_{\T}(\Prl)$ of left $\T$-module objects of $\Prl$. For instance, $\Prl_\st\simeq\LMod_{\Sp}(\Prl)$.
Of course, $\Sp\simeq\LMod_\SS$, so $\Prl_\st\simeq\LMod_{\LMod_\SS}(\Prl)$.
It is useful to consider the relative version of this: given a commutative ring spectrum $R$, the relative tensor product equips the presentable $\i$-category $\LMod_R$ with the structure of a commutative algebra object.
\begin{definition}
Let $R$ be a commutative ring spectrum.
An {\em $R$-linear $\i$-category} $\C$ is a left $\LMod_R$-module in $\Prl^\otimes$.
\end{definition}
We write $\Cat_R=\LMod_{\LMod_R}(\Prl)$
for the $\i$-category of $R$-linear categories and $R$-linear functors (left $\LMod_R$-module maps).
\begin{remark}
Since $\LMod_R$ is stable, we have
$\Cat_R\simeq\LMod_{\LMod_R}(\Prl_{\st})$.
\end{remark}

\begin{remark}
$\Cat_R$ is in fact {\em closed} symmetric monoidal: if $\C$ and $\D$ are $R$-linear $\i$-categories, the $\i$-category of $R$-linear functors (that is, $\Mod_R$-module morphisms in $\Prl$)
$
\Funl_R(\C,\D)
$
from $\C$ to $\D$ is again an $R$-linear $\i$-category (a $\Mod_R$-module in $\Prl$).
The dual
\begin{equation*}
    \mathrm{D}_R\C=\Funl_R(\C,\LMod_R)
\end{equation*}
of $\C$ is the $\i$-category of $R$-linear functors from $\C$ to $\LMod_R$.
Hence an $R$-linear $\i$-category $\C$ is dualizable if there exists a coevaluation map $\eta: \LMod_R\xrightarrow{\eta}\Dual_R\C\otimes_R\C$ such that the compositions
\begin{gather*}
\C\xrightarrow{\C\otimes_R\eta}\C\otimes_R\Dual_R\C\otimes_R\C\xrightarrow{\epsilon\otimes_R\C}\C\\
\Dual_R\C\xrightarrow{\eta\otimes\Dual_R\C}\Dual_R\C\otimes_R\C\otimes_R\Dual_R\C\xrightarrow{\Dual_R\C\otimes\epsilon}\Dual_R\C
\end{gather*}
are equivalent to the identity.
\end{remark}

\begin{proposition}
Let $R$ be a commutative ring spectrum.
Then $\Cat_R^\otimes$ is a rigid symmetric monoidal $\i$-category; that is, all objects are dualizable.
\end{proposition}

\begin{remark}
Consider the subcategory
$
\Prl_{\st}^{\mathrm{cg}}\subset\Prl_{\st}
$
of compactly generated stable presentable $\infty$-category and left adjoint functors which preserve compact objects. Then $\Prl_{\st}^{\mathrm{cg}}$ is equivalent, via the functor which restricts to the full subcategories of compact objects, to the $\i$-category of idempotent complete (that is, all idempotents split) small stable $\i$-categories and exact functors.
The inverse equivalence is given by $\Ind:\Cat_\i\to\Prl$, which when restricted to $\Cat_\i^{\st}\subset\Cat_\i$, factors through the subcategory $\Prl_{\st}^{\mathrm{cg}}\subset\Prl$.
Since stability and idempotent completeness are properties of small $\i$-categories, and a functor of small stable $\i$-categories is exact if and only if it is right exact, we deduce that $\Prl_{\st}^\mathrm{cg}$ is equivalent to a full subcategory of $\Cat_\i^\mathrm{rex}$ (see \cite[Proposition 5.5.7.8]{HTT} for details).
It follows from \cite[Proposition 4.8.1.4]{HA}, using \cite[Remark 2.2.1.2]{HA}, that $\Prl_{\st}^\mathrm{cg}$ inherits a symmetric monoidal structure which is compatible with the symmetric monoidal structure on $\Prl_{\st}$ or $\Prl$.
The $\infty$-category $\Mod_R^\otimes$ is a commutative algebra object in
$\Prl_{\st}^{\mathrm{cg}}$, so we have a subcategory $\Cat_{R}^{\mathrm{cg}}\subset\Cat_R$ of compactly
generated $R$-linear categories and colimit and compact object preserving functors.
Similar arguments show that $\Cat_{R}^{\mathrm{cg}}$ inherits the structure of a symmetric monoidal $\i$-category from $\Cat_R^\otimes$.
\end{remark}

\begin{proposition}{\em \cite[Proposition 3.5]{AG14}}\label{prop:Acpt}
An $R$-algebra $A$ is compact in $\Alg_R$ if and only if $\LMod_A$ is compact in $\Cat_{R}^{\mathrm{cg}}$.
\end{proposition}

\begin{remark}
An object $\C$ is dualizable in $\Cat_{R}^{\mathrm{cg}}$ if and only if it is dualizable in $\Cat_R$ and the evaluation and coevaluation morphisms lie in the (not full) subcategory $\Cat_{R}^{\mathrm{cg}}\subset\Cat_R$.
\end{remark}

\begin{definition}
A compactly generated $R$-linear category $\C$ is {\em proper} if its evaluation map is in $\Cat_{R}^{\mathrm{cg}}$; it is {\em smooth} if it is dualizable in $\Cat_R$ and its  coevaluation map is in $\Cat_{R}^{\mathrm{cg}}$.
An $R$-algebra $A$ is {\em proper} if $\LMod_A$ is proper; it is {\em smooth} if $\LMod_A$ is smooth.
\end{definition}

\begin{remark}
The property of being compact, smooth, or proper in $\Alg_R$ is invariant under Morita equivalence.
\end{remark}
\begin{remark}
If $A$ is an $R$-algebra, then $\LMod_A$ is proper if and only if $A$ is a perfect
$R$-module. Indeed,  the evaluation map is the map
\begin{equation*}
\LMod_{A\otimes_R A^{\op}}\simeq\LMod_A\otimes_R\LMod_{A^{\op}}\longrightarrow\LMod_R
\end{equation*}
that sends $A\otimes_R A^{\op}$ to $A$.
Similarly,
$\LMod_A$ is smooth if and only if the coevaluation map
\[
\LMod_R\longrightarrow\LMod_{A^{\op}\otimes_R A},
\]
which sends $R$ to the
$A^{\op}\otimes_R A$-module $A$, exists and is in
$\Cat_{R}^{\mathrm{cg}}$. So we see that $\LMod_A$ is smooth if and only
if $A$ is perfect as an $A^{\op}\otimes_R A$-module.
\end{remark}
\begin{proposition}{\em \cite[Lemma 3.9]{AG14}}
If $\C$ is a smooth $R$-linear category then $\C\simeq\LMod_A$ for some $R$-algebra spectrum $A$.
\end{proposition}

\begin{proposition}{\em \cite[Proposition 7.3.5.8]{HA}}
Let $R$ be a commutative ring spectrum and $A$ an associative $R$-algebra.
If $A$ is compact in $\Alg_R$ then $A$ is smooth.
If $A$ is a smooth and proper then A is compact in $\Alg_R$.
\end{proposition}

\begin{corollary}{\em \cite[Corollary 7.3.5.9]{HA}}
Let $R$ be a commutative ring spectrum and let $A$ be an associative $R$-algebra.
Then A is smooth and proper if and only if it is compact as an object of $\Alg_R$ and $\LMod_R$.
\end{corollary}
\begin{remark}
The {\em noncommutative cotangent complex} is the fiber
\[
\Omega_{A/R}\too A\otimes_R A^{\op}\too A
\]
of the multiplication map.
Thus, in the noncommutative setting, $A$ is a smooth $R$-algebra if and only if $\Omega_{A/R}$ is a perfect left $A\otimes_R A^{\op}$-module.
\end{remark}

\subsection{Nilpotent, local, and complete objects}
In this section we fix a commutative ring spectrum $R$ and a finitely generated ideal $I=(f_1,\ldots,f_n)\subset\pi_0 R$.
We write $R[I^{-1}]\simeq\bigotimes_{i=1}^n R[f_i^{-1}]$ and, for any left $R$-module $M$, $M[I^{-1}]\simeq  R[I^{-1}]\otimes_R M$, regarded as a left $R$-module via the commutative algebra map $R\to R[I^{-1}]$.
\begin{remark}
There does not seem to be a good notion of {\em ideal} in the $\i$-category of commutative ring spectra.
Rather, it seems that for most part, which is relevant are ideals in the discrete ring $\pi_0 R$ or graded ring $\pi_* R$.
The problem is that the cofiber $R/f$ of an $R$-module map $f:R\to R$ need not admit a commutative multiplication, or even any multiplication at all.
A good example is the Moore spectrum $\SS/p$ for $p\in\ZZ\cong\pi_0\SS$ a prime, which does not carry an associative algebra structure; in fact, $\SS/2$ doesn't even support a unital binary multiplication map $\SS/2\otimes\SS/2\to\SS/2$ (see \cite[Proposition 4]{Sc10}).
\end{remark}

\begin{remark}
The {\em nilpotence theorem} of Devinatz-Hopkins-Smith \cite{DHS} states that if $R$ is a homotopy associative ring spectrum, then the kernel of the map $\pi_*(R)\to\pi_*(MU\otimes R)$ consists entirely of nilpotent elements.
This generalizes the {\em Nishida nilpotence theorem} \cite{Ni73}, which states that every element of $\pi_*\SS$ of positive degree is nilpotent.
It is used to show that  the only homotopy associative ring spectra $R$ with the property that $\pi_* R$ is a {\em graded field} are (extensions of) the Morava K-theory spectra $K(n)$ at the prime $p$ (suppressed from the notation) at height $n$.
Here $K(0)\simeq\QQ$,  $\pi_* K(n)\cong\FF_p[v^{\pm}]$ for a generator $v$ in degree $2(p^n-1)$, and $K(\infty)\simeq\FF_p$.
By \cite{An11}, there is an essentially unique associative $\SS$-algebra structure on $K(n)$, at least at odd primes.
\end{remark}

\begin{remark}
Given a left $\SS$-module $E$, the {\em Bousfield localization} $L_E\Sp$ of the $\i$-category of spectra at $E$ sits in a Verdier sequence
\[
K_E\Sp\too\Sp\too L_E\Sp,
\]
where $K_E\Sp\subset\Sp$ denotes the full subcategory of those spectra $M$ such that $M\otimes E\simeq 0$ (the kernel of the multiplication by $E$ map $\Sp\to\Sp$).
Since the inclusion of the full subcategory $K_E(\Sp)\subset\Sp$ evidently preserves colimits, it follows that this is a semi-orthogonal decomposition.
It is a recollemont precisely when the right adjoint inclusion $L_E\Sp\to\Sp$ preserves colimits, which implies that the localization is smashing: $L_E M\simeq M\otimes L_E\SS$.
\end{remark}

\begin{definition}
Let $R$ be a commutative ring spectrum, let $I\subset\pi_0 R$ be a finitely generated ideal, and let $A$ be an associative $R$-algebra.
\begin{itemize}\itemsep.1em
\item[\rm{(1)}]
A left $A$-module $M$ is $I$-nilpotent if $M[I^{-1}]\simeq 0$.
\item[\rm{(2)}]
A left $A$-module $M$ is $I$-local if $\Map(N,M)\simeq 0$ for every $I$-nilpotent $N$.
\item[\rm{(3)}]
A left $A$-module $M$ is $I$-complete if $\Map(L,M)\simeq 0$ for every $I$-local $L$.
\end{itemize}
\end{definition}
\begin{remark}
Taking homotopy groups, we see that $N$ is $I$-nilpotent if each $\pi_k N$, $k\in\ZZ$, is an $I$-nilpotent left $\pi_0 R$-module.
Equivalently, $N$ is $I$-nilpotent if every element of each $\pi_k N$ is annihilated by some power of $I$.
\end{remark}

\begin{remark}
It suffices to test $I$-locality only on those $I$-nilpotent objects which are also compact as left $A$-modules.
This is because the fully faithful inclusion $i_\lor:\LMod_A^{I\textrm{\!-\!}\nil}\subset\LMod_A$ of the $I$-nilpotent left $A$-modules is a colimit preserving functor of compactly generated stable $\i$-categories \cite[Proposition 7.1.1.12]{SAG}, which implies that it admits a right adjoint $i^\lor$.
For instance, if $I=(f)$ is principle, then $A/f$ is a compact generator of $\LMod_A^{I\textrm{\!-\!}\nil}\subset\LMod_A$ and $A[f^{-1}]$ is generator of $\LMod_A^{I\textrm{\!-\!}\loc}\subset\LMod_A$ which is compact as an object of $\LMod_A^{I\textrm{\!-\!}\loc}$ but not typically as an object of $\LMod_A$.
In particular, $\LMod_A^{I\textrm{\!-\!}\loc}\simeq\LMod_{A[I^{-1}]}$.
\end{remark}
\begin{remark}
The fully faithful inclusion  $j_*:\LMod_A^{I\textrm{\!-\!}\loc}\to\LMod_A$  of the $I$-local objects admits both a left adjoint $j^*$ and a right adjoint $j^\times$.
The existence of the left adjoint follows from the fact that $j_*$ preserves limits, essentially by definition, and filtered colimits by the previous remark.
The existence of the right adjoint follows from the fact that $j_*$ is exact and so it also preserves finite colimits, hence all colimits.
It follows that the $I$-localization functor $j^*:\LMod_A\to\LMod_A^{I\textrm{\!-\!}\loc}$ is given by tensoring with $j^*A\simeq A[I^{-1}]$.
\end{remark}
\begin{theorem}{\em \cite[Theorem 7.3.4.1]{SAG}}
A left $A$-module $M$ is $I$-complete if each homotopy group $\pi_n M$ is a derived $I$-complete left $\pi_0 A$-module, in the sense that $\Ext^0_{\pi_0 A}(\pi_0 A[f^{-1}],M)\cong 0\cong\Ext^1_{\pi_0 A}(\pi_0 A[f^{-1}],M)$.
\end{theorem}
\begin{remark}
The inclusion of the $I$-complete objects preserves limits, again by construction, and $\kappa$-filtered colimits for $\kappa\gg 0$.
Indeed, this follows from the fact that we need only test completeness on a generator $A[I^{-1}]$ of $\LMod_A^{I\textrm{\!-\!}\loc}$, and this generator is $\kappa$-compact in $\LMod_A$ for some $\kappa\gg 0$.
\end{remark}
\begin{proposition}{\em \cite[Proposition 7.2.4.4 and Proposition 7.3.1.4]{SAG}}
The fully faithful inclusions $i_\lor:\LMod_A^{I\textrm{\!-\!}\nil}\to\LMod_A$ and $j_*:\LMod_A^{I\textrm{\!-\!}\loc}\to\LMod_A$ induce Verdier sequences
\begin{gather*}
\LMod_A^{I\textrm{\!-\!}\nil}\overset{i_\lor}{\too}\LMod_A\overset{j^*}{\too}\LMod_A^{I\textrm{\!-\!}\loc}\\
\LMod_A^{I\textrm{\!-\!}\loc}\overset{j_*}{\too}\LMod_A\overset{i^\land}{\too}\LMod_A^{I\textrm{\!-}\mathrm{cpl}}
\end{gather*}
which are semi-orthogonal decompositions of $\LMod_A$.
\end{proposition}

\begin{remark}
In the stable setting, a {\em recollement}
\index{recollement}
is the data of a fully faithful inclusion of a stable subcategory which admits both a left and a right adjoint \cite{BG}.
An example is the inclusion $j_*:\LMod_A^{I\textrm{\!-\!}\loc}\to\LMod_A$ of the $I$-local objects.
The composite $i^\land i_\lor$ 
is an equivalence of categories with inverse $i^\lor i_\land$.
A stable recollement determines a {\em fracture square}, a cartesian square
\[
\xymatrix{
\id\ar[r]\ar[d] & i_\land i^\land\ar[d]\\
j_*j^*\ar[r] & j_*j^*i_\land i^\land}
\]
of endofunctors of $\LMod_A$.
For $A\simeq\SS$, the {\em arithmetic square} is the cartesian square
\[
\xymatrix{
M\ar[r]\ar[d] & \Big(\prod_{p} M_{p}^\wedge\Big)\ar[d]\\
M\otimes\QQ\ar[r] & \Big(\prod_{p} M_{p}^\wedge\Big)\otimes\QQ}
\]
obtained by completing a spectrum $M$ at all primes $p$ and rationalization.
\end{remark}
\begin{theorem}{\em \cite[Proposition 7.4.1.1]{SAG}}
Let $R$ be a commutative ring spectrum and let $I\subset\pi_0 R$ a finitely generated ideal.
Suppose given a map of associative $R$-algebras $f:A\to B$, and consider the commutative square
\[
\xymatrix{
\LMod_A\ar[r]^{f^*}\ar[d] & \LMod_B\ar[d]\\
\LMod_{A[I^{-1}]}\ar[r]^{f[I^{-1}]^*} & \LMod_{B[I^{-1}]}}
\]
in $\Prl$.
If $f^\wedge_I:A^\wedge_I\to B^\wedge_I$ is an equivalence, this square is cartesian.
\end{theorem}

\section{Deformation theory}\label{sec:dt}

\subsection{The tangent bundle and the cotangent complex}

The {\em cotangent complex formalism} is an instance of the fiberwise stabilization of a presentable  fibration.
Given a pair of presentable fibrations $p:\D\to\C$ and $q:\E\to\C$, we write $\Funr_\C(\E,\D)\subset\Fun_\C(\E,\D)$ for the full subcategory of those functors $g:\E\to\D$ over $\C$ which admit a left adjoint $f:\D\to\E$ such that $p(\eta(D))$ is an equivalence in $\C$ for every object $D\in\D$, where $\eta:\id_\D\to gf$ denotes any choice of unit transformation exhibiting the adjunction.
\begin{remark}
This is the precise condition needed to ensure that $g:\E\to\D$ restricts to a right adjoint on fibers over any object of $\C$, or even after pullback along any morphism $\C'\to\C$. 
\end{remark}
\begin{definition}
A {\em stable envelope}\index{stable envelope} of a presentable fibration $p:\D\to\C$ is a presentable fibration $q:\E\to\C$ equipped with a morphism 
\[
\xymatrix{
\E\ar[rr]^g\ar[rd]_q & & \D\ar[ld]^p\\
& \C &}
\]
of presentable fibrations over $\C$ which exhibits $\E$ as the fiberwise stabilization of $\D$ over $\C$ in the following sense: if $q':\E'\to\C$ is a stable presentable fibration, the induced map $g_*:\Funr_\C(\E',\E)\to\Funr_\C(\E',\D)$ is an equivalence.
\end{definition}

\begin{definition}
A {\em tangent bundle}\index{tangent bundle} $q:T_\C\to\C$ of a presentable $\infty$-category $\C$ is a stable envelope of the target fibration $p:\Fun(\Delta^1,\C)\to\Fun(\{1\},\C)\simeq\C$.
\end{definition}

\begin{remark}The fiber of $q:T_\C\to\C$ over the object $A\in\C$ is the stabilization $\Sp(\C_{/A})$ of the fiber $\C_{/A}$ of the target fibration $p:\Fun(\Delta^1,\C)\to\C$.
\end{remark}
\begin{remark}
The tangent bundle $T_\C$ admits an explicit construction as the $\i$-category of {\em unreduced} excisive functors\index{excisive functor!unreduced}
\[
\E\simeq\mathrm{Exc}(\S_*^\mathrm{fin},\C).
\]
In this case, the morphism of presentable fibrations $g:\E\to\D$ over $\C$ is given by evaluating  at the arrow $S^0\to\pt$; that is, $g(X)=\{X(S^0)\to X(\pt)\}$.
\end{remark}

Let $\C$ be a presentable $\i$-category, and consider the commutative triangle
\[
\xymatrix{
T_\C\ar[rr]^g\ar[rd]_q & & \Fun(\Delta^1,\C)\ar[ld]^p\\
& \C &}
\]
where $p$ is the target fibration.
A relative version of the adjoint functor theorem implies that $f$ admits a left adjoint,
allowing for the following construction.

\begin{definition}
The {\em absolute cotangent complex functor}
\index{cotangent complex!absolute}
$L:\C\to T_\C$ is  the composition $\C\to\Fun(\Delta^1,\C)\to T_\C$, where the first map is the diagonal (constant) embedding and the second is left adjoint to $g:T_\C\to\Fun(\Delta^1,\C)$.
We will write $L_A$ for the value of $L:\C\to T_\C$ at the object $A$ of $\C$, and refer to $L_A$ as the {\em (absolute) cotangent complex} of $A$.
\end{definition}
\begin{remark}
The cotangent complex has a long history.
Its first incarnation was as the sheaf of K\"ahler differentials, and a derived version of this was introduced by Berthelot and Illusie (\cite{SGA6}, \cite{Il71}, \cite{Il72}). 
Around the same time, Quillen (\cite{Qu67}, \cite{Qu70}) and Andr\'e (\cite{An74}) defined a derived version in the context of (simplicial) commutative rings, and later Basterra and Mandell (\cite{Ba99}, \cite{BM05}) developed the theory in the more general context of commutative algebra spectra, where they refer to it as {\em topological Andr\'e-Quillen homology}.
\end{remark}
\begin{remark}
Let $A$ be an object of $\C$.
The identification of the fiber of the tangent bundle $T_\C$ over $A\in\C$ with $\Sp(\C_{/A})$ is such that $L_A\in\Sp(\C_{/A})$ corresponds to the image of $\id_A\in\C_{/A}$ under  $\Sigma^\infty_+:\C_{/A}\to\Sp(\C_{/A})$.
\end{remark}
\begin{remark}
The diagonal embedding $\C\to\Fun(\Delta^1,\C)$ is left adjoint to the evaluation $\Fun(\Delta^1,\C)\to\Fun(\{0\},\C)\to\C$ at $0\in\Delta^1$. Hence the cotangent complex functor $L:\C\to T_\C$ is left adjoint to the composite functor
\[
T_\C\too\Fun(\Delta^1,\C)\too\Fun(\{0\},\C)\to\C.
\]
\end{remark}
\begin{definition}
The {\em relative cofiber over $\C$} functor
\[
\mathrm{cof}_\C\colon\Fun(\Delta^1,T_\C)\to T_\C
\]
is the functor which sends the morphism $f:X\to Y$ in $T_\C$ to the pushout
\[
\xymatrix{
X\ar[r]\ar[d] & Y\ar[d]\\
Z\ar[r] & \cof_\C(f)}
\]
where $Z$ is any zero object of the fiber of $T_\C$ over $p(X)$.
\end{definition}
\begin{definition}
The {\em relative cotangent complex functor}\index{cotangent complex!relative} is the composition
$
\Fun(\Delta^1,\C)\overset{L}{\too}\Fun(\Delta^1,T_\C)\overset{\mathrm{cof}_\C}{\too} T_\C.
$
\end{definition}

\begin{remark}
A zero object of a fiber of $q:T_\C\to\C$ need not be a zero object of $T_\C$ itself.
Given a morphism $f:X\to Y$ in $T_\C$, the relative cofiber $\cof_\C(f)$ of $f$ is an object of the fiber  of $T_\C$ over $q(Y)$.
\end{remark}
We write $L_{B/A}$ for the value of the relative cotangent complex functor on an object $f:A\to B$ of the $\i$-category $\Fun(\Delta^1,\C)$ of arrows in $\C$.
\begin{remark}
By construction, the relative cotangent complex of a morphism $f:A\to B$ fits into a {\em relative cofiber sequence}\index{cofiber sequence!relative}
\[
\xymatrix{
L_A\ar[r]\ar[d] & L_B\ar[d]\\
0\ar[r] & L_{B/A}}
\]
in  $T_\C$.
This induces an actual cofiber sequence $f_! L_A\to L_B\to L_{B/A}$ in the $\i$-category $T_\C\times_{\C}\{B\}\simeq\Sp(\C_{/B})$.
Here $f_!:\Sp(\C_{/A})\to\Sp(\C_{/B})$ is a straightening of the restriction $f^*q:T_\C\times_\C\Delta^1\to\Delta^1$ of $q:T_\C\to\C$ along $f:\Delta^1\to\C$.
\end{remark}

\begin{remark}
It follows that the commutative square in $T_\C$
\[
\xymatrix{L_{B/A}\ar[r]\ar[d] & L_{C/A}\ar[d]\\
L_{B/B}\ar[r] & L_{C/B}}
\]
associated to a pair of composible morphisms $A\to B$ and $B\to C$ in $\C$ is cocartesian, hence a relative cofiber sequence since $L_{B/B}\simeq 0$ in $\Sp(\C_{/B})$.
\end{remark}

\begin{definition}
The {\em tangent correspondence}\index{tangent correspondence} of a presentable $\i$-category $\C$ is the cocartesian fibration $\M_\C\to\Delta^1$ associated to the cotangent complex functor $L:\C\to T_\C$.
\end{definition}
\begin{remark}
The tangent correspondence $\M_\C\to\Delta^1$ is also a cartesian fibration since $L:\C\to T_\C$ is a left adjoint.
\end{remark}
\begin{remark}
The cocartesian fibration $\M\to\Delta^1$ associated to a functor $f:\C\to\D$ can be constructed as the pushout
\[
\xymatrix{
\C\ar[r]^f\ar[d] & \D\ar[d]\\
\Delta^1\times\C\ar[r] & \M}
\]
in which the left vertical map is the inclusion at $\{1\}\subset\Delta^1$.
The functor $\M\to\Delta^1$ has fiber over $0$ and $1$ the full subcategories $\C\to\M$ and $\D\to\M$, and over the unique map $\epsilon:0\to 1$ the functor $f:\C\to\D$; this is evidently a cocartesian fibration as we can push forward objects of $\C$ along $\epsilon$ via $f$.
Thus a functor $h:\M\to\C$ amounts to the data of a functor $g:\D\to\C$ and a natural transformation $\eta:\Delta^1\times\C\to\C$ from $\id_\C$ to $gf$.
If $\M$ is also a cartesian fibration, a unit transformation $\eta:\id_\C\to gf$ induces a canonical functor $h:\M\to\C$.
\end{remark}

\begin{remark}
A {\em derivation} in $\C$ is a morphism $\Delta^1\to\M_\C$ such that the composite $\Delta^1\to\M_\C\to\Delta^1$ is the identity and $\Delta^1\to\M_\C\to\C$ is constant.
More concretely, a derivation in $\C$ is a morphism in $\M_\C$ from an object $A$ in the fiber $\C$ over $0$ to an object $M$ in the fiber $\Sp(\C_{/A})$ of $q:T_\C\to\C$ over $A$ (see \cite[Remark 7.4.1.2]{HA}).
The $\i$-category $\Der(\C)$\index{$\Der(\C)$} of derivations in $\C$ is the pullback
\[
\xymatrix{
\Der(\C)\ar[r]\ar[d] & \Fun(\Delta^1,\M_\C)\ar[d]\\
\C\ar[r] & \Fun(\Delta^1,\Delta^1\times\C)}
\]
in which the right vertical map is induced by the functor $\M_\C\to\Delta^1\times\C$ and the bottom horizontal map is adjoint to the identity of $\Delta^1\times\C$.
\end{remark}

\subsection{Derivations and square-zero extensions}

We now specialize to the case in which $\C$ is the presentable $\infty$-category $\CAlg$ of commutative algebra spectra.

\begin{remark}
All algebras will be assumed commutative for the remainder of this article.
If $A$ is a commutative algebra object there is a canonical equivalence $A\simeq A^{\op}$, hence a canonical equivalence $\LMod_A\simeq\LMod_{A^{\op}}\simeq\RMod_A$, so we needn't distinguish left and right module structures in the commutative case. We thus write $\Mod_A$
\index{$\Mod_A$}
in place of $\LMod_A$ and $\RMod_A$.
\end{remark}

\begin{theorem}{\em \cite[Corollary 7.3.4.14]{HA}}
The functor $p\colon\Mod\to\CAlg$ which sends the pair $(A,M)\in\Mod$ to $A\in\CAlg$ exhibits $\Mod$ as a tangent bundle of $\CAlg$.
In particular, for any commutative algebra spectrum $A$, we have an equivalence $\Mod_A\simeq\Sp(\CAlg_{/A})$. 
\end{theorem}

We will write $\Sym^{\leq 1}:\Mod\to\Fun(\Delta^1,\CAlg)$
\index{$\Sym^{\leq 1}$}
for a functor which corresponds to $\Omega^\infty:T_{\CAlg}\to\Fun(\Delta^1,\CAlg)$.
We use this notation because $\Sym^{\leq 1}\simeq\Sym^0\oplus\Sym^1$ is the formula for the {\em split square-zero extension}\index{square-zero extension!split} in ordinary algebra, though we have obtained its augmented commutative algebra structure through abstract stabilization techniques.
\begin{remark}
Let $A$ be a commutative ring spectrum and $M$ an $A$-module.
Then the split square-zero augmented commutative $A$-algebra structure on $A\oplus M$ is square-zero in the sense that the compositions
\[
\Sym^n_A(M)\too\Sym^n_A(A\oplus M)\overset{\otimes^n}{\too} A\oplus M\too M
\]
are null whenever $n>1$.
\end{remark}
\begin{remark}
On homotopy groups, the split square-zero extension is an ordinary split square-zero extension of graded commutative rings.
In other words, for any pair of elements $(a_0,m_0)$ and $(a_1,m_1)$ of $\pi_*(A\oplus M)$, the multiplication on $\pi_*(A\oplus M)\cong\pi_*(A)\oplus\pi_*(M)$ is given by the expected formula
\[
(a_0,m_0)(a_1,m_1)=(a_0a_1,a_0m_1+(-1)^{|a_1||m_0|}a_1m_0).
\]
\end{remark}

\begin{remark}
Any split square-zero extension $A\oplus M$, viewed as an augmented commutative $A$-algebra, is canonically a spectrum object in the pointed $\infty$-category $\CAlg_{A/A}$. Indeed, in degree $n$,
$(A\oplus M)^n\simeq A\oplus\Sigma^n M$, and the map
\[
A\oplus\Sigma^n M\too\Omega(A\oplus\Sigma^{n+1} M)\simeq A\oplus\Omega\Sigma^{n+1}M\simeq A\oplus\Sigma^n M
\]
is an equivalence, as $A$ is a zero object of $\CAlg_{A/A}$ and
\[
\xymatrix{
A\oplus\Sigma^n M\ar[r]\ar[d] & A\ar[d]\\
A\ar[r] & A\oplus\Sigma^{n+1}M}
\]
is a cartesian square of $\CAlg_{A/A}$.
Moreover, as $\Mod_A$ is stable, the functor $\Sym^{\leq 1}_A:\Mod_A\to\CAlg_{A/A}$ factors through $\Omega^\infty:\Sp(\CAlg_{A/A})\to\CAlg_{A/A}$.
There is even an evident map back in the order direction given by taking the fiber, which sends the spectrum object $\{B^n\}_{n\in\NN}$ to $\fib\{B^0\to A\}\in\Mod_A$.
\end{remark}

In ordinary commutative algebra, a derivation $d:A\to M$ over $R$ is an $R$-module map satisfying the Leibniz rule
\[
d(ab)=ad(b)+bd(a).
\]
Instead of using elements, we could instead define a derivation $d:A\to M$ as a section of the projection $A\oplus M\to A$, taken in the category of commutative $R$-algebras.
Replacing $R$ with $\SS$ we obtain the following notion.
\begin{definition}
Let $A$ be a commutative algebra spectrum and $M$ an $A$-module.
A {\em derivation}\index{derivation} from $A$ to $M$ is a section $A\to A\oplus M$ in $\CAlg$ of the projection $A\oplus M\to A$.
\end{definition}
\begin{remark}
This is a special case of the notion of derivation introduced in the previous section.
Here we avoid explicit mention of the tangent correspondence by pulling back to the fiber of $\M_{\CAlg}\to\Delta^1$ via the split square zero extension functor $T_{\CAlg}\simeq\Mod\to\CAlg$.
\end{remark}
\begin{remark}
The space $\Der(A,M)$
\index{$\Der(A,M)$}
of derivations $A\to A\oplus M$ is the fiber
\[
\Der(A,M)\too\Map_{\CAlg}(A,A\oplus M)\too\Map_{\CAlg}(A,A)
\]
over the identity $\id_A\in\Map_{\CAlg}(A,A)$.
The composite of the commutative algebra map $A\to A\oplus M$ with the second projection $A\oplus M\to M$ is a map of spectra $d:A\to M$ which we will abusively refer to as the derivation.
\end{remark}
\begin{remark}
For any connective commutative algebra $A$ and left $A$-module $M$, there is an equivalence
$
\Map_{\Mod_A}(L_{A},M)\to\Der(A,M).
$
That is, the cotangent complex $L_A$ corepresents the functor $\Der(A,-):\Mod_A\to\S$.
\end{remark}

\begin{remark}
Using the tangent correspondence formalism from the previous subsection, the $\infty$-category $\Der=\Der(\CAlg)$
\index{$\Der$}
of derivations has objects derivations $d:A\to M$ and morphisms commutative squares of the form
\[
\xymatrix{
A\ar[r]\ar[d] & M\ar[d]\\
B\ar[r] & N},
\]
with $A\to B$ a commutative algebra map and $M\to N$ an $A$-module map.
\end{remark}

Let $A$ be a commutative ring spectrum, $M$ an $A$-module, and $\eta:A\to\Sigma M$ a derivation with associated section $(\id,\eta):A\to A\oplus\Sigma M$ of the projection $A\oplus\Sigma M\to A$.
By construction, $(\id,\eta)$ is a morphism of commutative algebras.

\begin{definition}
A map of  commutative algebra spectra $\epsilon:A'\to A$ is a {\em square-zero extension}
\index{square-zero extension}
of $A$ by the $A$-module $M$ if there exists a derivation $\eta:A\to\Sigma M$ and a cartesian square in $\CAlg_{/A}$ of the form
\[
\xymatrix{A'\ar[r]\ar[d] & A\ar[d]^{(\id,\eta)}\\
A\ar[r]^{(\id,0)} & A\oplus \Sigma M}.
\]
\end{definition}
\begin{remark}\label{Phi}
There is a functor $\Phi:\Der\to\Fun(\Delta^1,\CAlg)$
\index{$\Phi$}
which sends the derivation $\eta:A\to\Sigma M$ to the map $A^\eta\to A$, where $A^\eta$ is a pullback
\[
\xymatrix{A^\eta\ar[r]\ar[d] & A\ar[d]^{(\id,\eta)}\\
A\ar[r]^{(\id,0)} & A\oplus \Sigma M} 
\]
in $\CAlg_{/A}$,
and whose essential image consists of the square-zero extensions.
\end{remark}

\begin{remark}
Note that, if $f:A^\eta\to A$ is a square-zero extension of $A$ by $M$, then the fiber sequence of $A$-modules $M\to A\to A\oplus\Sigma M$ implies that we also have a fiber sequence $M\to A^\eta\to A$ and hence a canonical equivalence $\fib(f)\simeq M$.
Hence derivations $\eta:A\to\Sigma M$ give rise to square-zero extensions $A^\eta$ of $A$ by $M$.
This is why we use $\Sigma M$ instead of $M$ itself.
\end{remark}

\begin{remark}
As the name suggests, square-zero extensions $f:A'\to A$ are actually square-zero.
That is, the fiber $M\to A'\to A$ has the property that the $n$-fold symmetric power map $\Sym_A^n(M)\to M$ is null for any $n>1$.
\end{remark}

\begin{definition}
A morphism $\epsilon:A'\to A$ in $\CAlg^{\cn}$ is an $n$-{\em small extension}
\index{$n$-small extension}
if $\fib(\epsilon)$ has homotopy concentrated in degrees $[0,2n]$ and the multiplication map $\fib(\epsilon)\otimes_{A'}\fib(\epsilon)\to\fib(\epsilon)$ is nullhomotopic.
A derivation $\eta:A\to M$ is $n$-{\em small} if the associated square-zero extension $A^\eta\to A$ is $n$-small.
\end{definition}

\begin{theorem}{\em \cite[Theorem 7.4.1.26]{HA}}
The composition
\[
\Der_n\subset\Der\overset{\Phi}{\too}\Fun(\Delta^1,\CAlg^{\cn})
\]
is fully faithful with essential image the full subcategory of $\Fun(\Delta^1,\CAlg^{\cn})$ consisting of the $n$-small extensions (here $\Phi$ is as in \autoref{Phi} above).
\end{theorem}

\begin{corollary}
Any $n$-small extension is a square-zero extension.
\end{corollary}
One of the primary source of examples of square-zero extensions comes from the Postnikov tower of a connective commutative algebra $A$.
\begin{proposition}{\em \cite[Corollary 7.4.1.28]{HA}}
Let $A\in\CAlg^{\cn}$.
For each $n>0$, the map
\[
\tau_{\leq n}A\too\tau_{\leq n-1}A
\]
exhibits $\tau_{\leq n}A$ as a square-zero extension of $\tau_{\leq n-1} A$ by $\Sigma^n\pi_{n} A$.
\end{proposition}

\begin{remark}
The fact that the Postnikov tower is composed of square-zero extensions is one of the main reasons why the cotangent complex plays such an important role in spectral algebra and geometry.
The space of maps
\[
\Map_{\CAlg}(A,B)\simeq\lim\Map_{\CAlg}(A,\tau_{\leq n} B)\simeq\lim\Map_{\CAlg}(\tau_{\leq n}A,\tau_{\leq n}B)
\]
between connective commutative algebra spectra $A$ and $B$ decomposes as the limit of the space of maps between their truncations, and for any $n>0$ we have a pullback diagram
\[
\xymatrix{
\tau_{\leq n}B\ar[r]\ar[d] & \tau_{\leq n-1} B\ar[d]\\
\tau_{\leq n-1}B\ar[r] & \tau_{\leq n-1}B\oplus\Sigma^{n+1}\pi_n B}.
\]
This implies that the fibers of
$
\Map_{\CAlg}(A,B)\too\Map_{\CAlg^\heartsuit}(\pi_0 A,\pi_0 B)
$
are accessible via infinitesimal methods: arguing inductively up the Postnikov tower, we are reduced to understanding spaces of derivations from $A\to\Sigma^{n+1}\pi_nB$, an $A$-linear question concerning maps $L_A\to\Sigma^{n+1}\pi_n B$.
\end{remark}

\subsection{Deformations of commutative algebras}

Given a connective commutative algebra spectrum $A$ and a square-zero extension $A'\to A$ of $A$ by a connective $A$-module $M$, it is natural to study the space of {\em deformations} of the square-zero extension $A'\to A$ to a connective commutative $A$-algebra $f:A\to B$.
That is, we wish to understand the space of cocartesian squares in $\CAlg^{\cn}$ of the form
\[
\xymatrix{
A'\ar[r]^{f'}\ar[d] & B'\ar[d]\\
A\ar[r]^f & B}
\]
(it turns out that $B'\to B$ is automatically also a square-zero extension).
As in ordinary algebra, this reduces to a module theoretic question concerning the cotangent complex:
such a commutative $A'$-algebra $B'$ exists if and only if the map $B\otimes_A L_A\to  B\otimes_A\Sigma M$ induced by a derivation $\eta:A\to\Sigma M$ classifying $A'\to A$ factors through the absolute cotangent complex $L_B$ of $B$.

\begin{definition}
Let $A$ be a commutative algebra spectrum, $A'\to A$
a square-zero extension of $A$ by an $A$-module $M$, and $B$ a commutative $A$-algebra.
A {\em deformation of $B$ to $A' $}
\index{deformation!of commutative algebra spectra}
is a commutative $A'$-algebra $B'$ equipped with an equivalence
$B'\otimes_{A'} A\to B$ of commutative $A$-algebras.
\end{definition}
\begin{remark}
We need not assume that $B'$ is flat over $A'$, as this is the case if and only if $B$ is flat over $A$.
Indeed, if $A'\to B'$ is flat, then the basechange $A\to B\simeq B'\otimes_{A'} A$ is flat by a spectral sequence argument.
Conversely, the fact that any deformation $B'\to B$ of $A'\to A$ along a flat map $A\to B$ results in a flat map $A'\to B'$ follows from a tor-amplitude argument.
\end{remark}
\begin{remark}
In the connective case, just like in ordinary commutative algebra, there are cohomological obstructions to the existence and uniqueness of deformations.
\index{deformation!obstruction classes}
Given an $A$-linear map $\eta:L_A\to\Sigma M$ with $M$ a connective $A$-module, the associated square-zero extension $A^\eta\to A$ of $A$ by $M$ is again connective.
Given a morphism of connective commutative algebra spectra $A\to B$, we obtain a map $\eta_B:B\otimes_A L_A\to B\otimes_A\Sigma M$, and a deformation $B'\to B$ of $A'\to A$ along $A\to B$ exists if and only if $\eta_B$ factors as a composite
\[
B\otimes_A L_A\overset{\epsilon_B}{\too} L_B\overset{\eta'}{\too} B\otimes_A\Sigma M
\]
of $\epsilon_B:B\otimes_A L_A\to L_B$ and a $B$-module map $\eta':L_B\to B\otimes_A\Sigma M$.
In this case, the deformation is recovered as the square-zero extension $B'\simeq B^{\eta'}\to B$.
Notice, though, that such a factorization exists if and only if the composite
\[
\Sigma^{-1}L_{B/A}\too B\otimes_A L_A\too B\otimes_A\Sigma M,
\]
is null, so that a deformation exists if and only if the {\em obstruction class} in $\Ext^2_B(L_{B/A}, B\otimes_A M)$ corresponding to this map vanishes.
The set of equivalence classes of deformations is a torsor for the group $\Ext^1_B(L_{B/A}, B\otimes_A M)$; in particular, it might be empty, or only admit an element locally on $\Spec\pi_0 B$.
\end{remark}
\begin{definition}
A derivation $\eta:A\to\Sigma M$ is said to be {\em connective} if $A$ is a connective commutative algebra spectrum and $M$ is an connective $A$-module.
\end{definition}
Let $\D\subset\Der$ denote the subcategory consisting of the connective derivations $A\to\Sigma M$ and those morphisms of connective derivations
\[
\xymatrix{A\ar[r]\ar[d] &\Sigma M\ar[d]\\
B\ar[r] & \Sigma N}
\]
such that the induced map $B\otimes_A M\to N$ is an equivalence.
Similarly, let $\C\subset\Fun(\Delta^1,\CAlg)$ denote the subcategory consisting of those objects $A'\to A$ such that both $A$ and $A'$ are connective and those morphisms the squares
\[
\xymatrix{
A'\ar[r]\ar[d] & A\ar[d]\\
B'\ar[r] & B}
\]
which are cocartesian in $\CAlg$; in other words, $B'\otimes_{A'} A\to B$ is an equivalence.
When restricted to $\D\subset\Der$, the square zero extension functor $\Phi:\Der\to\Fun(\Delta^1,\CAlg)$ of \autoref{Phi} factors through $\C\subset\Fun(\Delta^1,\CAlg)$.
\begin{theorem}{\em \cite[Theorem 7.4.2.7]{HA}}\label{thm:lf}
The composition
\[
\D\subset\Der\overset{\Phi}{\too}\Fun(\Delta^1,\CAlg)
\]
factors through the subcategory $\C\subset\Fun(\Delta^1,\CAlg)$, and the resulting functor $\Phi':\D\to\C$ is a left fibration (a cocartesian fibration with $\i$-groupoid fibers).
\end{theorem}

\begin{proposition}{\em \cite[Proposition 7.4.2.5]{HA}}\label{prop:key}
For any connective derivation $\eta:A\to\Sigma M$,
$\Phi:\Der\to\Fun(\Delta^1,\CAlg)$ induces an equivalence
\[
\Phi_{\eta/}:\D_{\eta/}\overset{\simeq}{\too}\CAlg^{\cn}_{A^\eta},
\]
where $\D\subset\Der$ denotes the subcategory defined above.
\end{proposition}

\subsection{Connectivity results}

A {\em universal derivation}
\index{derivation!universal}
is any derivation $d:A\to L_A$ which is corresponds to an equivalence $L_A\to L_A$.
Given a map of commutative ring spectra $f:A\to B$,
\[
\xymatrix{
A\ar[r]\ar[d]^f & L_A\ar[d]\\
B\ar[r] & L_B}
\]
is a commutative square in $\Der$.
Taking vertical cofibers, we obtain map of $A$-modules $\cof(f)\to L_{B/A}$ which is adjoint to a map of $B$-modules $\cof(f)\otimes_A B\to L_{B/A}$.
We will write $\epsilon_f:\cof(f)\otimes_A B\to L_{B/A}$ for this map.

\begin{theorem}{\em \cite[Theorem 7.4.3.1]{HA}}\label{thm:2ncon}
Let $f:A\to B$ be a morphism in $\CAlg^{\cn}$ such that $\cof(f)\in\Mod_A^{\geq n}$ for some $n\geq 0$.
Then $\fib(\epsilon_f)\in\Mod_B^{\geq 2n}$.
\end{theorem}

\begin{example}
We note the elementary fact that if $M\in\Mod_A^{\geq n}$ then, for any natural number $m$, $\Ten_A^m(M)\in\Mod_A^{\geq mn}$ and consequently $\Sym_A^m(M)\in\Mod_A^{\geq mn}$ as well.
Thus if $f:\Sym_A M\to A$ is the projection to $A\simeq\Sym^0_A M$, it is straightforward to show that $\fib(\epsilon_f)\in\Mod_A^{\geq 2n}$.
The general case is obtained from this special case by connectivity and induction arguments.
\end{example}

\begin{proposition}{\em \cite[Corollary 7.4.3.2]{HA}}
Let $f:A\to B$ be a map of connective commutative algebra spectra such that $\cof(f)$ is
$n$-connective for some $n\geq 0$. Then the relative cotangent complex $L_{B/A}$ is $n$-connective, and the converse holds provided that $\pi_0 f:\pi_0 A\to\pi_0 B$ is an isomorphism.
\end{proposition}

\begin{remark}
The absolute cotangent complex of a connective commutative algebra spectrum is
itself connective.
This follows immediately from the previous proposition since the cofiber of the unit map is connective.
\end{remark}
\begin{corollary}
 Let $f : A\to B$ be a map of connective commutative algebra spectra. Then $f$ is an equivalence if
and only if $\pi_0 f:\pi_0 A\to\pi_0 B$ is an isomorphism and the relative cotangent complex $L_{B/A}$ vanishes.
\end{corollary}
\begin{remark}
Let $f : A\to B$ be a map of connective commutative algebra spectra such that $\cof(f)$ is $n$-connective for some $n \geq 0$.
The induced map $L_A\to L_B$ factors as the composite
\[
L_A\overset{g}{\too} B\otimes_A L_A\overset{h}{\too} L_B
\]
and the equivalence $\cof(g)\simeq\cof(f)\otimes_A L_A$, together with the connectivity of $A$ and $L_A$, imply that $\cof(g)$ is $n$-connective.
We also have an exact triangle
\[
B\otimes_A\cof(f)\too L_{B/A}\too\cof(\epsilon_f)
\]
in which $B\otimes_A\cof(f)$ and $\cof(\epsilon_f)$ are $n$-connective, so that $L_{B/A}$ is $n$-connective as well.
It follows that the cofiber of $L_A\to L_B$ is $n$-connective..
\end{remark}

\begin{proposition}{\em \cite[Lemma 7.4.3.8]{HA}}
Let A be a connective commutative algebra spectrum.
There are canonical isomorphisms of $\pi_0$-modules
\[
\pi_0 L_A\simeq\pi_0 L_{\pi_0 A}\simeq\Omega_{\pi_0 A}.
\]
\end{proposition}

\begin{remark}
Let $f:A\to B$ be a map of connective commutative algebra spectra such that $\cof(f)$ is $n$-connective.
Tensoring the exact triangle $A\to B\to\cof(f)$ with the $A$-module $\cof(f)$, we obtain an exact triangle
\[
\cof(f)\overset{\delta}{\too} B\otimes_A\cof(f)\too\cof(f)\otimes_A\cof(f)
\]
which exhibits $\delta$ as a $(2n-1)$-connective map.
Composing with $\epsilon_f:B\otimes_A\cof(f)\to L_{B/A}$, we obtain a $(2n-1)$-connective map $\cof(f)\to L_{B/A}$.
\end{remark}

\begin{proposition}{\em \cite[Proposition 7.4.3.9]{HA}}\label{prop:pi0L}
Let $f:A\to B$ be a morphism in $\CAlg^{\cn}$.
Then $L_{B/A}$ is connective and $\pi_0L_{B/A}\cong\Omega_{\pi_0 B/\pi_0 A}$.
\end{proposition}

\begin{theorem}{\em \cite[Theorem 7.4.3.18]{HA}}\label{prop:perfection}
Let $f:A\to B$ be a morphism in $\CAlg^{\cn}$.
If $B$ is locally of finite presentation over $A$, $L_{B/A}$ is a perfect $B$-module.
The converse holds provided $\pi_0 B$ is of finite presentation over $\pi_0 A$.
\end{theorem}

\subsection{Classification of \'etale maps}

Recall that a map  of discrete commutative rings $f:A\to B$ is said to be \'etale if $B$ is a finitely presented flat commutative $A$-algebra such that the multiplication map $B\otimes_A B\to B$ is the projection onto a summand. Geometrically, this is the algebraic analogue of a (not necessarily surjective) covering space: there exists a commutative ring $C$ and a cartesian square of schemes of the form 
\[
\xymatrix{
\Spec(B)\coprod\Spec(C)\ar[r]\ar[d] & \Spec(B)\ar[d]_f\\
\Spec (B)\ar[r]^f & \Spec(A)}.
\]
\'Etale maps of commutative rings $A\to B$ are smooth of relative dimension zero \cite[Lemma 10.141.2]{stacks}. In particular, there exists a presentation of $B$ as a commutative $A$-algebra of the form $B\cong A[x_1,\ldots,x_n]/(f_1,\ldots,f_n)$, where the $\{f_i\}_{1\leq i\leq n}$ are a sequence of elements of $A[x_1,\ldots,x_n]$ such that the image of the Jacobian matrix of partial derivatives $\{\partial f_i/\partial x_j\}_{1\leq i,j\leq n}$ is invertible in $B$.

\begin{definition}
A map $f:A\to B$ of commutative algebra spectra is {\'etale}
\index{\'etale!map}
if $f$ exhibits $B$ as a finitely presented commutative $A$-algebra such that $B$ is a flat $A$-module and $\pi_0 f:\pi_0 A\to\pi_0 B$ is \'etale.
\end{definition}

Ideally, we would like to be able to calculate the space of \'etale maps between commutative algebra spectra $A$ and $B$ in terms of the set of \'etale maps between their underlying discrete commutative algebras $\pi_0 A$ and $\pi_0 B$.
Since an \'etale map $f:A\to B$ induces an \'etale map $\pi_0 f:\pi_0 A\to\pi_0 B$, we can address this question directly by studying the space of \'etale lifts $f:A\to B$ of an \'etale morphism $f_0:\pi_0 A\to\pi_0 B$.
The main result of this section, one of the major results of higher algebra, is that the space of such lifts is contractible.

\begin{proposition}{\em \cite[Remark 7.5.1.7]{HA}}\label{prop:mbemie}
Given a commutative triangle
\[
\xymatrix{
& A\ar[rd]^g\ar[ld]_f &\\
B\ar[rr]^h & & C}
\]
 in $\CAlg$, if $f$ is \'etale, then $g$ is \'etale if and only if $h$ is \'etale.
 In particular, any map between \'etale commutative $A$-algebras is automatically \'etale.
\end{proposition}

\begin{remark} Let $f : A\to B$ be a map of connective commutative ring spectra. Then $f$ is \'etale if and only if each  $\tau_{\leq n}f:\tau_{\leq n} A\to\tau_{\leq n} B$ is \'etale.
\end{remark}

\begin{remark}
Since an \'etale map of discrete rings is smooth of relative dimension zero, its module of relative K\"ahler differentials vanishes, and one can show using ordinary algebraic methods that its relative cotangent complex also vanishes.
This begs the question of whether or not the relative cotangent complex $L_{B/A}$ of an \'etale map of commutative algebra spectra $f:A\to B$ vanishes.
Using flatness, one reduces to the connective case, so that $L_{B/A}$ is also connective, and if $L_{B/A}\neq 0$  there's a least $n\in\NN$ with $\pi_n L_{B/A}\neq 0$.
This is a contradiction: by connectivity considerations as in the previous section, 
\[
0\cong\pi_n L_{\pi_0 B/\pi_0 A}\cong\pi_n(L_{B/A}\otimes_B\pi_0 B)\cong\pi_n L_{B/A}.
\]
Hence $L_{B/A}\simeq 0$ for any \'etale map $f:A\to B$.
\end{remark}

\begin{remark}
Said differently, this means that the absolute cotangent complex functor $L:\CAlg\to\Mod$ satisfies \'etale basechange:
\index{\'etale!basechange}
if $f:A\to B$ is an \'etale map of commutative algebra spectra, then the exact triangle
\[
B\otimes_A L_A\too L_B\too L_{B/A}
\]
together with the vanishing of the relative cotangent complex $L_{B/A}$ implies that $B\otimes_A L_A\simeq L_B$.
\end{remark}

Let $\Der^{\et}\subset\Der$ denote the (not full) subcategory consisting of the connective derivations $A\to\Sigma M$ and those morphisms of connective derivations
\[
\xymatrix{A\ar[r]\ar[d] &\Sigma M\ar[d]\\
B\ar[r] & \Sigma N}
\]
such that the induced map $B\otimes_A M\to N$ is an equivalence and the commutative algebra map $A\to B$ is \'etale.
Let $\CAlg^{\et}\subset\CAlg^{\cn}$ denote the (not full) subcategory consisting of the (connective) commutative algebra spectra and the \'etale maps.
\begin{remark}
The forgetful functor $\Der^{\et}\to\CAlg^{\et}$ is a cocartesian fibration such that, for each $A\in\CAlg^{\et}$, the fiber $\Der^{\et}_A$ is an $\i$-groupoid (i.e., it is a left fibration).
This is because an \'etale morphism of derivations is cocartesian: given connective derivations $\eta:A\to\Sigma M$ and $\eta':B\to\Sigma N$ and a map $f:A\to B$, we have an equivalence $N\simeq B\otimes_A M$, so if $f=\id_A$ we obtain an equivalence $N\simeq M$.
In particular, for any connective derivation $\eta:A\to\Sigma M$, the forgetful functor induces an equivalence
$
\Der^{\et}_{\eta/}\overset{\simeq}{\too}\CAlg_{A}^{\et}.
$
\end{remark}

\begin{remark}
Consider the functor $\Phi:\Der\to\Fun(\Delta^1,\CAlg)$ which sends the derivation $\eta:A\to\Sigma M$ to the square-zero extension $A^\eta\to A$.
Compositing with the restrictions to the set of vertices $\{0,1\}$ of $\Delta^1$, we obtain functors $\Phi_0:\Der\to\CAlg$ and $\Phi_1:\Der\to\CAlg$ such that $\Phi_1$ restricts to the left fibration $\Der^{\et}\to\CAlg^{\et}$.
For a given connective derivation $\eta:A\to\Sigma M$, $\Phi_0$ and $\Phi_1$ induce functors
\[
\CAlg^{\et}_{A^\eta}\overset{\Phi'_0}{\from}\Der^{\et}_{\eta/}\overset{\Phi'_1}{\too}\CAlg^{\et}_A
\]
such that $\Phi'_0$ is an equivalence by \autoref{prop:key} and $\Phi'_1\simeq\Phi'_0\otimes_{A^\eta} A$ and $\Phi'_0$ is an equivalence by the remark above.
We obtain the following corollary.
\end{remark}
\begin{corollary} Let $f:A'\to A$ be a square-zero extension of connective commutative ring spectra. The relative
tensor product
\[
(-)\otimes_{A'} A:\CAlg_{A'}\too\CAlg_A
\]
induces an equivalence from the $\infty$-category of \'etale commutative $A'$-algebras to the $\infty$-category of \'etale commutative $A$-algebras.
\end{corollary}

\begin{theorem}{\em \cite[Corollary 7.5.4.3]{HA}}\label{thm:et}
For any commutative algebra spectrum $A$, 
$
\pi_0 : \CAlg_A\to\CAlg_{\pi_0 A}
$
induces an equivalence  $\CAlg_A^{\et}\simeq\CAlg_{\pi_0 A}^{\et}$.
\end{theorem}

\begin{remark}
Note that any \'etale commutative $\pi_0 A$-algebra $B$ is automatically in $\CAlg_A^\heartsuit$ since the flatness condition implies that $B$ must be discrete: $\pi_n B\cong\pi_n A\otimes_{\pi_0 A}\pi_0 B\cong 0$ if $n\neq 0$.
Hence $\CAlg_{\pi_0 A}^{\et}\simeq\CAlg_{\pi_0 A}^{\heartsuit\et}$, and the theorem asserts that the $\infty$-category of \'etale commutative $A$-algebras is equivalent to the {\em ordinary} category of \'etale commutative $\pi_0 A$-algebras.
\end{remark}

\begin{remark}
If $R$ is a commutative ring spectrum, the structure of the $\i$-category $\CAlg_R^{\et}$ of \'etale commutative $R$-algebras implies that the small \'etale site of $R$ is equivalent to the small \'etale site of the discreet ring $\pi_0 R$, and the analogous result holds for the small Zariski sites of \autoref{eg:zar}.
These facts form the cornerstones of {\em spectral algebraic geometry}, as treated in \cite{SAG} and \cite{TV08}.
\end{remark}
\begin{remark}
There are robust notions of spectral scheme and Deligne-Mumford stack.
Versions of the Artin representability
\index{representability theorems}
theorem for these higher categorical objects are formulated and proved as \cite[Theorems 18.1.0.1 and 18.1.0.2]{SAG}, providing necessary and sufficient conditions for a functor $F:\CAlg^{\cn}\to\S$ to be represented by such an object.
These conditions are surprisingly straightforward: the ordinary stack $F|_{\CAlg^\heartsuit}\to\S$ must be represented by an ordinary scheme or Deligne-Mumford stack, $F$ must admit a cotangent complex, $F$ must preserve limits of Postnikov towers, and $F$ must preserve pullbacks of diagrams of the form $A\to C\leftarrow B$ in $\CAlg^{\cn}$ in which both of the maps $A\to C$ and $B\to C$ surjective on $\pi_0$ with nilpotent kernel.
\end{remark}

\printindex

{\small
\bibliographystyle{plain}
\bibliography{bib}}

\vspace{30pt}
\small
\noindent
David Gepner\\
School of Mathematics and Statistics\\
The University of Melbourne\\
Parkville VIC 3010 Australia\\
\texttt{david.gepner@unimelb.edu.au}

\end{document}